\documentclass[oneside]{amsart}
\usepackage{amsmath,amssymb,color,amsthm,graphicx,cite}
\usepackage{color,bm}
\graphicspath{{./img_pdf/}}
\usepackage{enumerate} 
\usepackage[backref]{hyperref}
\RequirePackage{luatex85}
\usepackage[all]{xy}
\newtheorem{theorem}{Theorem}[section]
\newtheorem{proposition}[theorem]{Proposition}
\newtheorem{lemma}[theorem]{Lemma}
\newtheorem{corollary}[theorem]{Corollary}
\newtheorem{claim}{Claim}

\theoremstyle{definition}
\newtheorem{definition}{Definition}
\newtheorem{main}{Theorem}

\newtheorem{main_cor}[main]{Corollary}

\def\Sv{\mathop{\mathrm{Sing}}(v)} 
\def\Pv{\mathop{\mathrm{Per}}(v)} 
 
\def\Z{\mathbb{Z} }
\def\R{\mathbb{R} }
\def\A{\mathbb{A} }
\def\T{\mathbb{T} }
\def\bp{\partial^{+} }
\def\nbd{neighborhood }
\def\nbds{neighborhoods }

\author{Tomoo Yokoyama}
\date{\today}
\address{Department of Mathematics, Faculty of Science, Saitama University, Shimo-Okubo 255, Sakura-ku, Saitama-shi, 338-8570 Japan\\}
\email{tyokoyama@rimath.saitama-u.ac.jp}

\subjclass[2010]{}

\makeatletter
\@namedef{subjclassname@2020}{%
  \textup{2020} Mathematics Subject Classification}
\makeatother

\title[A Poincar\'{e}-Bendixson theorem for flows]{A Poincar\'{e}-Bendixson theorem for flows with arbitrarily many singular points}


\begin{document}

\baselineskip=17pt

\maketitle

\begin{abstract}
The Poincar\'{e}-Bendixson theorem is one of the most fundamental tools to capture the limit behaviors of orbits of flows. It was generalized and applied to various phenomena in dynamical systems, differential equations, foliations, group actions, translation lines, and semi-dynamical systems. On the other hand, though the no-slip boundary condition is a fundamental condition in differential equations and appears in various fluid phenomena, and Lakes of Wada attractors naturally occur in discrete and continuous real dynamical systems and complex dynamics, no generalizations of the Poincar\'{e}-Bendixson theorem can be applied to any differential equations with no-slip boundary condition on surfaces with boundary and flows with Lakes of Wada attractors. To analyze them, we generalize the Poincar\'{e}-Bendixson theorem into one for flows with arbitrarily many singular points on possibly non-compact surfaces by introducing some concepts to describe limit behaviors and using methods of foliation theory and general topology. 
\end{abstract}

\section{Introduction}

The Poincar\'{e}-Bendixson theorem is one of the most fundamental tools to capture the limit behaviors of orbits of flows and was applied to various phenomena (e.g. \cite{bhatia1966application,du2021traveling,koropecki2019poincare,hajek1968dynamical,roussarie2020topological,roussarie2021some,roussarie2021some02,pokrovskii2009corollary}). 
The theorem is stated by Poincar\'{e} in \cite{poincare1882courbes} for analytic vector fields on the plane and proved by Bendixson in \cite{bendixson1901courbes} for $C^1$ vector fields on the plane. 
In \cite{birkhoff1927dynamical}, Birkhoff introduced the concepts of $\omega$-limit set and $\alpha$-limit set of a point. 
Using these concepts, one can simply describe the limit behaviors of orbits stated in the works of Poincar\'{e} and Bendixson in detail. 
In fact, the Poincar\'{e}-Bendixson theorem was generalized for flows on surfaces in various ways  \cite{andronov1966qualitative,aranson1996introduction,buendia2019top,ciesielski1994poincare,demuner2009poincare,gardiner1985structure,gutierrez1986smoothing,Levitt1982foliation,lopez2004accumulation,lopez2007topological,markley1969poincare,marzougui1996,marzougui1998structure,nikolaev1999flows,schwartz1963generalization,vanderschoot2003local,yano1985asymptotic}, and also for foliations \cite{levitt1987differentiability,plante1973generalization}, translation lines on the sphere \cite{koropecki2019poincare}, geodesics for a meromorphic connection on Riemann surfaces \cite{abate2011poincare,abate2016poincare}, group actions \cite{hounie1981minimal}, and semidynamical systems \cite{bonotto2008limit}. 
For instance, the following statement holds (cf. \cite{nikolaev1999flows}):
The $\omega$-limit set of an orbit of a flow with finitely many fixed points on a compact surface is either a closed orbit, a semi-attracting limit circuit, or a Q-set, where a Q-set is the closure of non-closed recurrent orbit. 
A part of this classification is based on the following Ma\v{i}er's result \cite{mayer1943trajectories} (see  \cite[Theorem 2.4.4 p.32]{nikolaev1999flows}, \cite[Theorem 4.2]{aranson1996maier} for general cases, and \cite[Theorem 19]{markley2008recurrent} for orientable hyperbolic cases for details): Any point contained in an $\omega$-limit set of some point for a flow on a compact surface whose $\omega$-limit set contains non-closed orbits is positively recurrent. 
Furthermore, in \cite{vinograd1952limit,lopez2004accumulation}, it is shown that the $\omega$-limit set of a non-recurrent point is the boundary of an open annulus. 
Moreover, in \cite{marzougui1996}, the Poincar\'{e}-Bendixson theorem is generalized for flows with totally disconnected singular point sets on closed orientable surfaces as follows: The $\omega$-limit set of an orbit of such a flow is either a closed orbit, a union of singular points and of connecting separatrices, or a Q-set. 
Here a non-singular orbit is a connecting separatrix if each of the $\omega$-limit and the $\alpha$-limit sets is a singular point. 
However, any differential equations with no-slip boundary condition on surfaces with boundary can be applied no such generalizations of the Poincar\'{e}-Bendixson theorem to, because of the degeneracy of singular points.

%

On the other hand, the $\omega$-limit sets of orbits of analytic flows on the plane, the sphere, and the projective plane are studied in \cite{lopez2007topological}. 
However, Buend{\'i}a and L{\'o}pez pointed out a gap in a crucial lemma in the work by constructing a counterexample in the case of the sphere minus two points \cite{buendia2017rem}. 
Though the lemma is not generally true, they fixed the gap in the case of the sphere, the plane, the projective plane, and the projective plane minus one point. 
They also showed that $\omega$-limit sets of orbits of analytic flows on open connected subsets of the sphere are essentially the boundaries of simply connected Peano subcontinua \cite{buendia2019top}.


\subsection{Statements of main results}

To describe the Poincar\'{e}-Bendixson theorem for flows with arbitrarily many singular points on surfaces, we recall some concepts and introduce two concepts (quasi-circuit and quasi-Q-set) as follows. 
A closed connected invariant subset is a {\bf non-trivial quasi-circuit} if it is a boundary component of an open annulus, contains a non-recurrent point, and consists of non-recurrent points and singular points. 
The $\omega$-limit set $\omega(x)$ of a point $x$ is a {\bf quasi-semi-attracting limit quasi-circuit} with respect to a positive invariant small collar $\A_{-1}$ if $\omega(x)$ is a non-trivial quasi-circuit and the collar $\A_{-1}$ contains a quasi-semi-attracting collar basin of $\omega(x)$ (see Definition~\ref{def:qst_lqc} for details). 
A {\bf Q-set} is the closure of a non-closed recurrent orbit. 
A Q-set is a {\bf transversely Cantor} Q-set if there is a small neighborhood $U$ of a non-closed recurrent point of the Q-set $\mathcal{M}$ such that $\mathcal{M} \cap U$ is a product of an open interval and a Cantor set. 
An $\omega$-limit set of a point is a {\bf quasi-Q-set} if it intersects an essential closed transversal infinitely many times.
A non-recurrent orbit is a {\bf connecting quasi-separatrix} if each of the $\omega$-limit and the $\alpha$-limit sets is contained in a boundary component of the singular point set.
A subset is {\bf locally dense} if its closure has a nonempty interior. 
A topological space is {\bf locally connected} at a point $x$ if there is a small connected \nbd of $x$. 
A subset is {\bf locally connected} if it is locally connected at any point in it. 

We generalize the Poincar\'{e}-Bendixson theorem for flows with the totally disconnected singular point set on a closed orientable surface (see \cite[Theorem~3.1]{marzougui1996} (cf. \cite[Theorem~2.1]{marzougui2000flows})) as follows. 

\begin{main}\label{main:a}
The following statements hold for a flow with arbitrarily many singular points on a compact surface: 
\\
{\rm(a)} The $\omega$-limit set of any non-closed orbit is one of the following exclusively:
\begin{quote}
\setlength{\leftskip}{-25pt}$(1)$ A nowhere dense subset of singular points.
\\
$(2)$ A semi-attracting limit cycle.
\\
$(3)$ A quasi-semi-attracting limit quasi-circuit that is the image of a circle.
\\
$(4)$ A quasi-semi-attracting limit quasi-circuit that is not locally connected. 
\\
$(5)$ A locally dense Q-set. 
\\
$(6)$ A transversely Cantor Q-set. 
\\
$(7)$ A quasi-Q-set that consists of singular points and non-recurrent points.  
\end{quote}
{\rm(b)} Every non-recurrent orbit in the $\omega$-limit set of a point is a connecting quasi-separatrix. 
\\
{\rm(c)} If the singular point set is totally disconnected, then any non-recurrent orbits in the $\omega$-limit set of a point are connecting separatrices.  
\\
{\rm(d)} If the $\omega$-limit set of a point is a Q-set, then the Q-set corresponds to the orbit closure of any non-closed recurrent point in the Q-set. 
\end{main}

The previous theorem implies a generalization of the Poincar\'{e}-Bendixson theorem for a flow with countably many singular points on a compact surface (see Theorem~\ref{cor:PB_countable}). 

Using the classification of $\omega$-limit sets, we can show the non-existence of quasi-Q-sets on a sphere and a projective plane (see Corollary~\ref{cor:nonex_qQset}). 
Moreover, we demonstrate the dependency between the $\omega$-limit set and the $\alpha$-limit set (i.e. the dual concept of the $\omega$-limit set) of a point \cite{yokoyama2022omega}. 
For instance, the $\omega$-limit set of a point whose $\alpha$-limit set is a locally dense Q-set either is a locally dense Q-set or consists of singular points. 
Moreover, it is known that the number of Q-sets for flows on compact surfaces is finite \cite{markley1969poincare,markley1970number} and that Q-sets are essential. 
Therefore it is naturally asked whether the number of quasi-Q-sets for flows on compact surfaces is finite and whether quasi-Q-sets are essential \footnote{These questions are suggested by an  anonymous referee}. 
We show the finiteness and the essential property of quasi-Q-sets (see Lemma~\ref{lem:ess_qqset}, Proposition~\ref{prop:num_qqset}, and Proposition~\ref{prop:num_qqset_02}). 
In addition, we topologically characterize quasi-Q-sets (see Proposition~\ref{prop:ch_qqset} and Proposition~\ref{prop:ch_qqset_nontrivial}). 

The previous theorem can be applied to capture a time-reversal symmetric condition for limit sets. 
In fact, the characterization of a flow that is either irrational or Denjoy by Athanassopoulos \cite{athanassopoulos1992characterization} is refined as an application of the previous theorem in \cite{yokoyama2021flows}. 
From the construction in \cite[Example~3]{yokoyama2021flows}, notice that a subsets of singular points can become Lakes of Wada attractors (or Plykin type attractors) and that quasi-semi-attracting limit quasi-circuits that is not locally connected can contain a Wada-Lakes-like structure. 
By this construction, cutting closed transversals and collapsing the new boundary components, we can construct flows on a sphere with Lakes of Wada attractors and with an arbitrarily large number of complementary domains,  which are flow versions of such attractors of spherical homeomorphisms constructed by Boro\'{n}ski, \v{C}in\v{c}, and Liu \cite{boronski2020prime} and such an attractor of a transcendental entire function constructed by Mart{\'\i}-Pete, Rempe and Waterman \cite{marti2021wandering}. 
In particular, such constructions of homeomorphisms are motivated by a generalization of the Poincar\'{e}-Bendixson theorem \cite{koropecki2019poincare}. 
Moreover, using the previous theorem, for any Hamiltonian flow with arbitrarily many singular points on a compact surface, it can be shown that the $\omega$-limit set of any non-closed orbit consists of singular points \cite{yokoyama2022omega}.
In addition, the previous theorem implies a generalization of the Ma\v{i}er's description of recurrence as follows. 

\begin{main}\label{main:rec}
Let $v$ be a flow on a compact surface $S$ and $\mathop{\mathrm{Cl}}(v)$ the union of closed orbits. 
The following statements hold for a point $x \in \omega(z)$ for some point $z \in S$: 
\\
{\rm (1)} $\omega(x) \setminus \mathop{\mathrm{Cl}}(v) \neq \emptyset$ if and only if $x$ is non-closed positively recurrent. 
\\
{\rm (2)} $\alpha(x) \setminus \mathop{\mathrm{Cl}}(v) \neq \emptyset$ if and only if $x$ is non-closed negatively recurrent. 
\\
{\rm (3)} $(\omega(x) \cup \alpha(x)) \setminus \mathop{\mathrm{Cl}}(v) \neq \emptyset$ if and only if $x$ is non-closed recurrent.
\\
{\rm (4)} $\omega(x) \setminus \mathop{\mathrm{Cl}}(v) \neq \emptyset$ and $\alpha(x) \setminus \mathop{\mathrm{Cl}}(v) \neq \emptyset$ if and only if $x$ is non-closed Poisson stable.
\end{main} 

The previous theorem implies the following topological characterizations of non-closed positive recurrence. 

\begin{main_cor}\label{main:rec02}
The following are equivalent for a point $x$ for a flow on a compact surface $S$: 
\\
{\rm (1)} The point $x$ is non-closed positively recurrent. 
\\
{\rm (2)} $\omega(x) \setminus \mathop{\mathrm{Cl}}(v) \neq \emptyset$ and there is a point $z \in S$ with $x \in \omega(z)$. 
\\
{\rm (3)} $\omega(x) \setminus \mathop{\mathrm{Cl}}(v) \neq \emptyset$ and there is a point $z \in S$ with $x \in \alpha(z)$. 
\end{main_cor} 

Using the end completion of surfaces of finite genus and finitely many boundary components, we can obtain analogous results for such surfaces (see Theorem~\ref{main:b}, Theorem~\ref{main:c}, ant Corollary~\ref{cor:rec02} for details). 
Furthermore, applying Theorem~\ref{main:c}, minimal flows on compact surfaces are characterized and the Poincar\'{e} recurrence theorem for flows on surfaces is generalized \cite{yokoyama2024topological}.

Recall that a flow $v \colon \R \times Z \to Z$ is topologically semi-conjugate to a flow $w \colon \R \times Y \to Y$ via $h \colon Y \to Z$ if $h$ is a continuous surjection such that $v(t, h(y)) = h(w(t,y))$ for any $(t,y) \in \R \times Y$. 
We also introduce the following blow-up operation, which can modify any limit circuits into quasi-circuits that are not circuits (and more generally modify any non-locally-dense $\omega$-limit sets into $\omega$-limit sets that are not arcwise-connected).

\begin{main}\label{main:deformation_02}
Let $v$ be a flow on a surface $S$ with an $\omega$-limit set $\omega$ of a point containing a non-singular point  $x$ and with a point $y \in S - \omega$ satisfying $\omega(y) = \omega$. 
Then there is a flow $\widetilde{v}$ on $S$ satisfying the following properties: 
\\
{\rm(1)} The $\omega$-limit set $\omega_{\widetilde{v}}(y)$ is not arcwise-connected. 
\\
{\rm(2)} The restriction $v\vert_{S - \omega}$ of $v$ is topologically equivalent to the restriction $\widetilde{v}\vert_{S - \omega_{\widetilde{v}}(y)}$. 
\\
{\rm(3)} The flow $v_x$ is topologically equivalent to some flow $v'$ which is topologically semi-conjugate to the flow $\widetilde{v}$, where $v_x$ is the resulting flow of $v$ by replacing $x$ with a singular point (see Lemma~\ref{lem:deformation_01} for details of the definition of $v_x$). 
\\
{\rm(4)} The topological semi-conjugacy from $\widetilde{v}$ to $v'$ can be obtained by collapsing a connected closed invariant subset of $\omega_{\widetilde{v}}(y)$ into a singleton. 
\end{main}

The present paper consists of seven sections.
In the next section, as preliminaries, we introduce fundamental concepts.
In \S 3, we generalize the Poincar\'{e}-Bendixson theorem to one for a flow with arbitrarily many singular points on a compact surface. 
Moreover, we demonstrate the finiteness and the essential property of quasi-Q-sets and topologically characterize quasi-Q-sets. 
In \S 4, the Poincar\'{e}-Bendixson theorem is generalized to one for a flow with arbitrarily many singular points on a surface of finite genus and finitely many boundary components. 
Moreover, we characterize the recurrence, which is a generalization of the Ma\v{i}er's description of recurrence. 
In \S 5, the total disconnectivity of singular points implies that a limit quasi-circuit is the image of a circle, and the finiteness of singular points implies that a limit quasi-circuit is a limit circuit. 
Morevoer, the countability of singular points implies that a quasi-Q-set is a Q-set. 
These reductions imply a proof of a generalization of the Poincar\'{e}-Bendixson theorem for a flow with finitely many singular points on a compact surface. 
In addition, we obtain a generalization of the Poincar\'{e}-Bendixson theorem for a flow with countably many singular points on a compact surface. 
In \S 6, we introduce an operation that makes $\omega$-limit sets not arcwise-connected by constructing flow boxes with non-arcwise-connected invariant subsets. 
In the final section, we state some examples with $\omega$-limit sets which are non-locally-connected subsets of singular points, quasi-circuits that are non-circuits, and quasi-Q-sets that are non-Q-sets respectively.

\section{Preliminaries}

We recall the topological notion and the notion of dynamical systems. 

\subsection{Topological notion}
Denote by $\overline{A}$ the closure of a subset $A$ of a topological space, by $\mathop{\mathrm{int}}A$ the interior of $A$, and by $\partial A := \overline{A} - \mathop{\mathrm{int}}A$ the boundary of $A$, where $B - C$ is used instead of the set difference $B \setminus C$ when $C \subseteq B$.
We define the {\bf coborder} $\bm{\bp A}$ of $A$ by $\overline{A} - A$ and the {\bf border} $\bm{\partial^-  A}$ by $A - \mathop{\mathrm{int}}A$ of $A$. 
Then $\partial A = \partial^-  A \sqcup \bp A$, where $\sqcup$ denotes a disjoint union.
A boundary component of a subset $A$ is a connected component of the boundary of $A$. 
A subset is {\bf locally dense} if its closure has a nonempty interior. 

\subsubsection{Curves and loops}
A {\bf curve} is a continuous mapping $C: I \to Y$ where $I$ is a non-degenerate connected subset of a circle $\mathbb{S}^1$.
A curve is simple if it is injective.
We also denote by $C$ the image of a curve $C$.
Denote by $\partial C := C(\partial I)$ the boundary of a curve $C$ if $C$ can be extended into a continuous map whose domain is $I \cup \partial I$, where $\partial I$ is the boundary of $I \subset \mathbb{S}^1$. 
Put $\mathrm{int} C := C \setminus \partial C$ if $\partial C$ is defined. 
A simple curve is a simple closed curve if its domain is $\mathbb{S}^1$ (i.e. $I = \mathbb{S}^1$).
A simple closed curve is also called a {\bf loop}. 
An {\bf arc} is a simple curve whose domain is an interval. 

\subsubsection{Essential property}
A subset $A$ of a compact surface $S$ is {\bf inessential} (cf. \cite[2.4]{koropecki2014strictly}) if there is an open disk in $S^*$ which is a \nbd of $A^*$, where $S^*$ is the resulting closed surface from $S$ by collapsing all boundary components into singletons and $A^*$ is the resulting subset from $A$.
A subset $A$ of a compact surface $S$ is {\bf essential} if it is not inessential. 
Note that a loop in the interior $S - \partial S$ of a compact surface $S$ is essential if and only if it is not null homotopic in $S^*$. 

\subsection{Notion of dynamical systems}
By a {\bf surface}, we mean a paracompact two dimensional manifold, that does not need to be orientable.
A {\bf flow} is a continuous $\R$-action on a manifold.
From now on, we suppose that flows are on surfaces unless otherwise stated.
Let $v : \R \times S \to S$ be a flow on a surface $S$.
For $t \in \R$, define $v_t : S \to S$ by $v_t := v(t, \cdot )$.
For a point $x$ of $S$, we denote by $O(x)$ the orbit of $x$ (i.e. $O(x) := \{ v_t(x) \mid t \in \R \}$), $O^+(x)$ the positive orbit (i.e. $O^+(x) := \{ v_t(x) \mid t > 0 \}$), and $O^-(x)$ the negative orbit (i.e. $O^-(x) := \{ v_t(x) \mid t < 0 \}$).
A point $x$ of $S$ is {\bf singular} if $x = v_t(x)$ for any $t \in \R$ and is {\bf periodic} if there is a positive number $T > 0$ such that $x = v_T(x)$ and $x \neq v_t(x)$ for any $t \in (0, T)$.
A point is {\bf closed} if it is singular or periodic.
An orbit is singular (resp. periodic, closed) if it contains a singular (resp. periodic, closed) point. 
Denote by $\mathop{\mathrm{Sing}}(v)$ the set of singular points and by $\mathop{\mathrm{Per}}(v)$ (resp. $\mathop{\mathrm{Cl}}(v)$) the union of periodic (resp. closed) orbits. 
The subset $\mathop{\mathrm{Sing}}(v)$ is called the singular point set. 
The $\omega$-limit (resp. $\alpha$-limit) set of a point $x$ is $\omega(x) := \bigcap_{n\in \mathbb{R}}\overline{\{v_t(x) \mid t > n\}}$ (resp.  $\alpha(x) := \bigcap_{n\in \mathbb{R}}\overline{\{v_t(x) \mid t < n\}}$).
An $\omega$-limit set of a point is locally dense if it has a nonempty interior. 
Similarly, an $\alpha$-limit set of a point is locally dense if it has a nonempty interior. 
For an orbit $O$, define $\omega(O) := \omega(x)$ and $\alpha(O) := \alpha(x)$ for some point $x \in O$.
Note that an $\omega$-limit (resp. $\alpha$-limit) set of an orbit is independent of the choice of point in the orbit.

\subsubsection{Topological equivalence}
A flow $v$ on a surface $M$ is {\bf topologically equivalent} to a flow $w$ on a surface $N$ if there is a homeomorphism $h \colon M \to N$ such that the images of any orbits of $v$ are orbits of $w$ with preservation of the direction in time. 
Then the homeomorphism $h \colon M \to N$ is called the {\bf topologically equivalent homeomorphism}. 

\subsubsection{Separatrices}
A {\bf separatrix} is a non-singular orbit whose $\alpha$-limit or $\omega$-limit set is a singular point.
A separatrix is {\bf connecting} if each of its $\omega$-limit set and $\alpha$-limit sets is a singular point.
Note that any connecting separatrices are non-recurrent.

\subsubsection{Orbit arcs}
An arc is an {\bf orbit arc} if it is contained in an orbit. 
An orbit arc in a subset $A$ is {\bf maximal} in $A$ if it is an orbit arc in $A$ which is maximal with respect to the inclusion order. 
In other words, for an orbit $O$ and for an orbit arc $C \subseteq O$ in a subset $A$ is maximal in $A$ if and only if the orbit arc is a path component of $A \cap O$. 
Note that a maximal orbit arc of an orbit $O$ in a subset $A$ need not be a connected component of $O \cap A$. 
In fact, let $v$ be a minimal flow on a torus $\mathbb{T}^2$, $x \in \mathbb{T}^2$ a point, and $A := \mathbb{T}^2 - \{ x \}$ the complement of $x$. 
Then the positive orbit $O^+(x)$ of $x$ is a maximal orbit arc of $O(x) \cap A$, but the set difference $O(x) - \{ x \}$ is the connected component of $O(x) \cap A$. 
Indeed, assume that $O(x) - \{ x \}$ is not a connected component of $O(x) \cap A$. 
Take disjoint nonempty open subsets $U$ and $V$ of $O(x) - \{ x \}$ with $O(x) - \{ x \} \subset U \sqcup V$. 
Since the positive orbit $O^+(x)$ is connected in either $O^+(x) \subset U$ or $O^+(x) \subset V$. 
By renaming if necessary, we may assume that $O^+(x) \subset U$. 
By $\overline{O^+(x)} = \mathbb{T}^2$, we have that $\emptyset \neq O^+(x) \cap V \subset U \cap V$, which contradicts $U \cap V = \emptyset$. 
This means that $O(x) - \{ x \}$ is the connected component of $O(x) \cap A$.

%
%
%
%
%

\subsubsection{Recurrence and invariance}
A point $x$ is {\bf positively recurrent} (or positively Poisson stable) if $x \in \omega(x)$. 
A point $x$ is {\bf negatively recurrent} (or negatively Poisson stable) if $x \in \alpha(x)$, 
A point $x$ is {\bf recurrent} if $x \in \omega(x) \cup \alpha(x)$. 
A point $x$ of $S$ is {\bf Poisson stable} (or strongly recurrent) if $x \in \omega(x) \cap \alpha(x)$.  
Denote by $\mathrm{R}(v)$ the set of non-closed recurrent points.
A point is {\bf wandering} if there are its neighborhood $U$ and a positive number $N$ such that $v_t(U) \cap U = \emptyset$ for any $t > N$. 
A point is non-wandering if it is not wandering. 
Note that a recurrent point is non-wandering. 
An orbit is recurrent (resp. Poisson stable, wandering, non-wandering) if it contains a recurrent (resp. Poisson stable, wandering, non-wandering)  point. 
A non-closed recurrent orbit is also called a non-trivial recurrent orbit. 

\begin{definition}
The closure of a non-closed recurrent orbit is called a {\bf Q-set} (or quasi-minimal set). 
\end{definition}

A subset is {\bf invariant} (or {\bf saturated}) if it is a union of orbits. 
The {\bf saturation} of a subset is the union of orbits intersecting it.
A nonempty closed invariant subset is {\bf minimal} if it contains no proper nonempty closed invariant subsets. 
A subset $A$ is {\bf positive invariant} if $v(t,A) \subseteq A$ for any $t \in \R_{\geq 0}$. 
A subset $A$ is {\bf negative invariant} if $v(t,A) \subseteq A$ for any $t \in \R_{\leq 0}$. 
Recall that the {\bf (orbit) class} $\hat{O}$ of an orbit $O$ is the union of orbits each of whose orbit closure equals $\overline{O}$ (i.e. $\hat{O} = \{ y \in S \mid \overline{O(y)} = \overline{O} \} $).

\subsubsection{Topological properties of orbits}

The following properness, local density, and exceptional properties of orbits are analogous concepts of codimension one foliation theory (cf. \cite{Hector1981,Candel2000foliation}). 

An orbit $O$ is {\bf proper} if there is its \nbd $U$ with $\overline{O} \cap U = O$. 
Note that an orbit $O$ is proper if and only if it is an embedded submanifold. 
Moreover, any closed orbit is proper. 
Recall that an orbit is locally dense if and only if the closure of the orbit has a nonempty interior. 
%
An orbit is exceptional if it is neither proper nor locally dense. 
A point is proper (resp. locally dense) if its orbit is proper (resp. locally dense).
Denote by $\mathrm{LD}(v)$ (resp. $\mathrm{E}(v)$, $\mathrm{P}(v)$) the union of locally dense orbits (resp. exceptional orbits, non-closed proper orbits). 
We have the following observation. 

\begin{lemma}{\cite[Lemma~2.1]{yokoyama2024topological}}\label{lem:decomp}
The following statements hold for a flow $v$ on a paracompact manifold $M$:
\\
{\rm(1)} A point of $M$ is non-proper if and only if it is non-closed recurrent. 
\\
{\rm(2)} $M = \mathop{\mathrm{Cl}}(v) \sqcup \mathrm{P}(v) \sqcup \mathrm{R}(v) = \mathop{\mathrm{Sing}}(v) \sqcup \mathop{\mathrm{Per}}(v) \sqcup \mathrm{P}(v) \sqcup \mathrm{LD}(v) \sqcup \mathrm{E}(v)$. 
\\
{\rm(3)} The union $\mathrm{P}(v)$ is the set of non-recurrent points. 
\\
{\rm(4)} The union $\mathrm{R}(v) = \mathrm{LD}(v) \sqcup \mathrm{E}(v)$ is the set of non-proper points. 
\end{lemma}

For the self-containedness, we prove the previous lemma in the different way from the proof of \cite[Lemma~2.1]{yokoyama2024topological} as follows. 

\begin{proof}
By definitions, note that a closed orbit is proper and recurrent. 
Fix a point $x \in M$.  
By \cite[Corollary~3.4]{yokoyama2019properness}, the orbit $O(x)$ is proper if and only if $O(x) = \hat{O}(x)$. 
From \cite[Theorem~VI]{cherry1937topological}, the closure of a non-closed recurrent orbit $O$ of a flow on a manifold contains uncountably many non-closed recurrent orbits whose closures are $\overline{O}$. 
Therefore the orbit $O(x)$ is non-closed recurrent if and only if $\hat{O}(x)$ consists of uncountably many orbits. 
Then any non-closed proper orbit is not recurrent. 

We claim that assertion {\rm(1)} holds. 
Indeed, suppose that $x$ is non-closed recurrent. 
Since any non-closed proper orbit is not recurrent, the point $x$ is not proper. 
Conversely, suppose that $x$ is not proper. 
Since an orbit $O$ is proper if and only if $O = \hat{O}$, we have $O(x) \subsetneq \hat{O}(x)$. 
Then $x$ is not closed and there is a point $y \in \hat{O}(x) - O(x)$ whose orbit does not contain $x$ such that $\overline{O(y)} = \overline{O(x)}$. 
This implies that $y \in \overline{O(x)} - O(x) \subseteq \alpha(x) \cup \omega(x)$. 
From the closedness and the invariance of $\alpha$-limit sets and $\omega$-limit sets, we have $x \in \overline{O(x)} = \overline{O(y)} \subseteq \alpha(x) \cup \omega(x)$. 
This means that $x$ is recurrent. 

Since $\mathop{\mathrm{Cl}}(v)$ is both the set of closed proper points and the set of closed recurrent points, assertion {\rm(1)} implies assertions {\rm(2)}--{\rm(4)}. 
\end{proof}



\subsubsection{Flow boxes}

We define a trivial flow box as follows. 

\begin{figure}
\begin{center}
\includegraphics[scale=0.35]{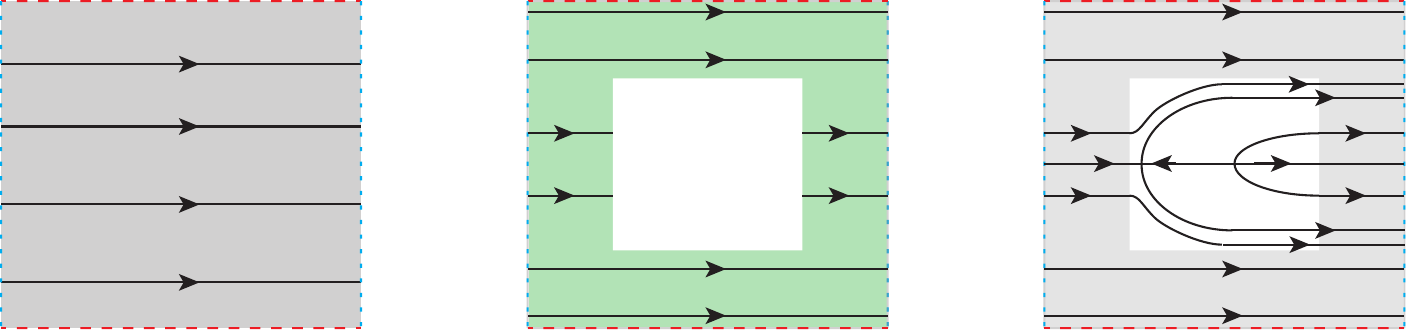}
\end{center}
\caption{Left, a trivial flow box; middle, the intersection of a flow box and its small neighborhood of the boundary, right; a Cherry flow box.}
\label{fig:flowbox}
\end{figure}

\begin{definition}
A disk $B$ on a surface $S$ is a {\bf trivial flow box} with respect to a flow $v$ on $S$ if there are nondegenerate intervals $I, J \subset \R$ and a homeomorphism $f \colon B \to I \times J \subset \R^2$ which carries the maximal orbit arc in $B$ to the maximal orbit arc in $I \times J$ with respect to the flow $v_X$ generated by a vector field $X = \partial/\partial x_1 = (1,0)$ on the plane $\R^2$ as in the left of Figure~\ref{fig:flowbox}. 
\end{definition}

Note that any orbit arcs on the square $I \times J$ are horizontal intervals and satisfy $\dot{x}_1 = 1$ and $\dot{x}_2 = 0$ with respect to the chart $(x_1, x_2)$ as in the left of Figure~\ref{fig:flowbox}. 
In other words, a disk $B$ on a surface $S$ is a trivial flow box with respect to a flow $v$ on $S$ if and only if there are a continuous vector field $X$ on a surface $T$, a topologically equivalent homeomorphism $h \colon S \to T$, and a chart $(x_1, x_2) \colon h(B) \to I \times J \subset \R^2$ such that $v$ is topologically equivalent to the flow $v_X$ generated by $X$ via $h$ and that $X$ is given by the expression $\partial/\partial x_1$ in the chart $(x_1, x_2)$. 

%
%

\begin{definition}
A disk $B$ on a surface $S$ is a {\bf flow box} with respect to a flow $v$ on $S$ if there are intervals $I, J$ which are either $(-1,1)$, $(-1,1]$, $[-1,1)$, or $[-1,1]$, and there are a closed disk $D \subset \mathop{\mathrm{int}} B$, a continuous vector field $X = \partial/\partial x_1$ on the plane $\R^2$, and a homeomorphism $f \colon B - D \to \mathbb{A} \subset \R^2$ which carries the maximal orbit arc in the annulus $B - D$ to the maximal orbit arc in $\mathbb{A}$, where $\mathbb{A} := (I \times J) - [-1/2, 1/2]^2$ is an annulus, as in the middle of Figure~\ref{fig:flowbox}. 
\end{definition}

In other words, a disk $B$ on a surface $S$ is a flow box with respect to a flow $v$ on $S$ if and only if there are a continuous vector field $X$ on a surface $T$, a topologically equivalent homeomorphism $h \colon S \to T$, a closed disk $D \subset \mathop{\mathrm{int}} B$, and a chart $(x_1, x_2) \colon h(B - D) \to \mathbb{A} \subset \R^2$ such that $v$ is topologically equivalent to the flow $v_X$ generated by $X$ via $h$ and that $X$ is given by the expression $\partial/\partial x_1 = (1,0)$ in the chart $(x_1, x_2)$. 

By definition, trivial flow boxes are flow boxes. 
The orbit arcs on the annulus $\mathbb{A}$ in the flow box $I \times J$ are horizontal intervals as in the middle of Figure~\ref{fig:flowbox}. 
Notice that Cherry flow boxes (cf. \cite{cherry1937topological} and Figure~1 in \cite{gardiner1985structure}) are flow foxes as in the right of Figure~\ref{fig:flowbox}. 

\subsubsection{Topological properties of Q-sets}

We observe the following characterization of local density. 

\begin{lemma}\label{lem:locally_dense_ch}
The following statements are equivalent for a non-closed recurrent orbit $O$
\\
{\rm(1)} The Q-set $\overline{O}$ is locally dense. 
\\
{\rm(2)} The orbit $O$ is locally dense. 
\\
{\rm(3)} The Q-set $\overline{O}$ contains locally dense orbits. 
\\
{\rm(4)} The Q-set $\overline{O}$ is a \nbd of $O$. 
\end{lemma}

\begin{proof}
By definition of local density for subsets, assertions {\rm(1)} and {\rm(2)} are equivalent. 
Obviously, assertion {\rm(4)} implies assertion {\rm(2)}, assertion {\rm(2)} implies assertion {\rm(3)}, and assertion {\rm(3)} implies assertion {\rm(1)}. 

Suppose $O$ is locally dense. 
The closure $\overline{O}$ has a nonempty interior. 
Fix a point $x \in \mathop{\mathrm{int}}\overline{O} =: U$. 
Since $x \in \overline{O}$, there is a point $y \in O \cap U$. 
Then $U$ is a \nbd of $y \in O$. 
Because $v(t, \cdot)$ is a homeomorphism, the image $v(t, U)$ is an open \nbd of $v(t, y)$ for any $t \in \R$. 
Therefore $O = \bigcup_{t \in \R} v(t, y) \subseteq \bigcup_{t \in \R} v(t, U) = v(U) = v(\mathop{\mathrm{int}}\overline{O}) \subseteq \overline{O}$ because of the invariance of $\overline{O}$. 
This means that $\overline{O}$ is a \nbd of $O$. 
\end{proof}

A Q-set is {\bf exceptional} if it is not locally dense. 
A Q-set is {\bf transversely Cantor} if there is a small neighborhood $U$ of a non-closed recurrent point of the Q-set $\mathcal{M}$ such that $\mathcal{M} \cap U$ is a product of an open interval and a Cantor set. 
We observe the following characterization of exceptional property. 

\begin{lemma}\label{lem:excptional_ch}
The following statements are equivalent for a non-closed recurrent orbit $O$
\\
{\rm(1)} The Q-set $\overline{O}$ is exceptional. 
\\
{\rm(2)} The orbit $O$ is exceptional. 
\\
{\rm(3)} The Q-set $\overline{O}$ contains exceptional orbits. 
\\
{\rm(4)} The Q-set $\overline{O}$ is a transversely Cantor Q-set. 
\end{lemma}

\begin{proof}
Obviously, assertion {\rm(2)} implies assertion {\rm(3)}. 
Suppose that $\overline{O}$ is exceptional. 
By Lemma~\ref{lem:locally_dense_ch}, the Q-set $\overline{O}$ contains no locally dense orbits. 
Since any non-closed recurrent orbit is not proper, the non-closed recurrent orbit $O$ is not locally dense and so is exceptional. 
This implies that assertion {\rm(1)} implies assertion {\rm(2)}.

Suppose that the Q-set $\overline{O}$ contains exceptional orbits. 
By \cite[Proposition~2.2]{yokoyama2016topological}, any exceptional Q-set contains no locally dense orbits.
Lemma~\ref{lem:locally_dense_ch} implies that $\overline{O}$ is not locally dense and so is exceptional. 
This means that assertion {\rm(3)} implies assertion {\rm(1)}.
 
Suppose that $\overline{O}$ is a transversely Cantor Q-set. 
Then $O$ is non-closed recurrent. 
Moreover, there are a non-closed recurrent point $x \in \overline{O}$ and a small neighborhood $U$ of $x$ such that $\overline{O} \cap U$ is a product of an open interval and a Cantor set. 
Then $\overline{O}$ is not a \nbd of $x$ and is not a \nbd of $O$. 
Lemma~\ref{lem:decomp} and Lemma~\ref{lem:locally_dense_ch} imply that the non-closed recurrent orbit $O$ is not locally dense and so exceptional. 
This implies that assertion {\rm(4)} implies assertion {\rm(2)}.

Suppose that $\overline{O}$ is exceptional. 
By \cite[Proposition~2.2]{yokoyama2016topological}, we have that $\hat{O} = \overline{O} \setminus (\mathop{\mathrm{Sing}}(v) \sqcup \mathrm{P}(v)) \subseteq \mathrm{R}(v)$. 
Since $O$ is exceptional and so not locally dense, we obtain $\mathop{\mathrm{int}}\overline{O} = \emptyset$. 
For any $y \in \hat{O}$, we have that $\mathop{\mathrm{int}}\overline{O(y)} = \mathop{\mathrm{int}}\overline{O} = \emptyset$ and so that $O(y)$ is not locally dense. 
Therefore $\hat{O} \subseteq \mathrm{R}(v) - \mathrm{LD}(v) = \mathrm{E}(v)$. 
Fix a non-closed recurrent point $x \in O$. 
Since $\overline{O}$ is not locally dense, there is a transverse closed arc $T$ whose interior contains $x$ such that $\overline{O} \cap \partial T = \emptyset$. 
Then there is a closed trivial flow box $U$ which is a \nbd of $x$ and contains no singular points such that $T \subset U$ and $\partial T \subset \partial U$. 
Then the intersection $\overline{O} \cap \partial T$ is a compact metrizable space. 
By time reversion if necessary, we may assume that $\overline{O} = \omega(x)$. 
By $\overline{O} \cap \partial T = \emptyset$, we obtain $\overline{O} \cap T \subset \mathop{\mathrm{int}}T$. 
Therefore the intersection $\overline{O} \cap T =  \omega(x) \cap \mathop{\mathrm{int}} T$ is perfect. 
Since $\overline{O}$ is not locally dense, the intersection $\overline{O} \cap \partial T$ contained in a closed interval is totally disconnected. 
Because a Cantor set is characterized as a compact metrizable perfect totally disconnected space, the intersection $\overline{O} \cap \partial T$ is a Cantor set. 
Therefore the intersection $\overline{O} \cap \mathop{\mathrm{int}} U$ is a product of an open interval and a Cantor set. 
This implies that $\overline{O}$ is a transversely Cantor Q-set. 
This means that assertion {\rm(1)} implies assertion {\rm(4)}.
\end{proof}

\subsubsection{Circuits}
An {\bf annular} subset is homeomorphic to an annulus. 
An open annular subset $\mathbb{A}$ of a surface is a {\bf collar} of a singular point $x$ if the union $\mathbb{A} \sqcup \{ x \}$ is a neighborhood of $x$. 
By a {\bf cycle} or a periodic circuit, we mean a periodic orbit.

\begin{definition}
A {\bf circuit} is one of the following subsets: 
\\
{\rm(1)} A singular point.
\\
{\rm(2)} A cycle. 
\\
{\rm(3)} An image of an oriented circle by a continuous orientation-preserving mapping which is a directed graph but not a singleton and which is the union of separatrices and finitely many singular points.
\end{definition}

\begin{figure}
\begin{center}
\includegraphics[scale=0.2]{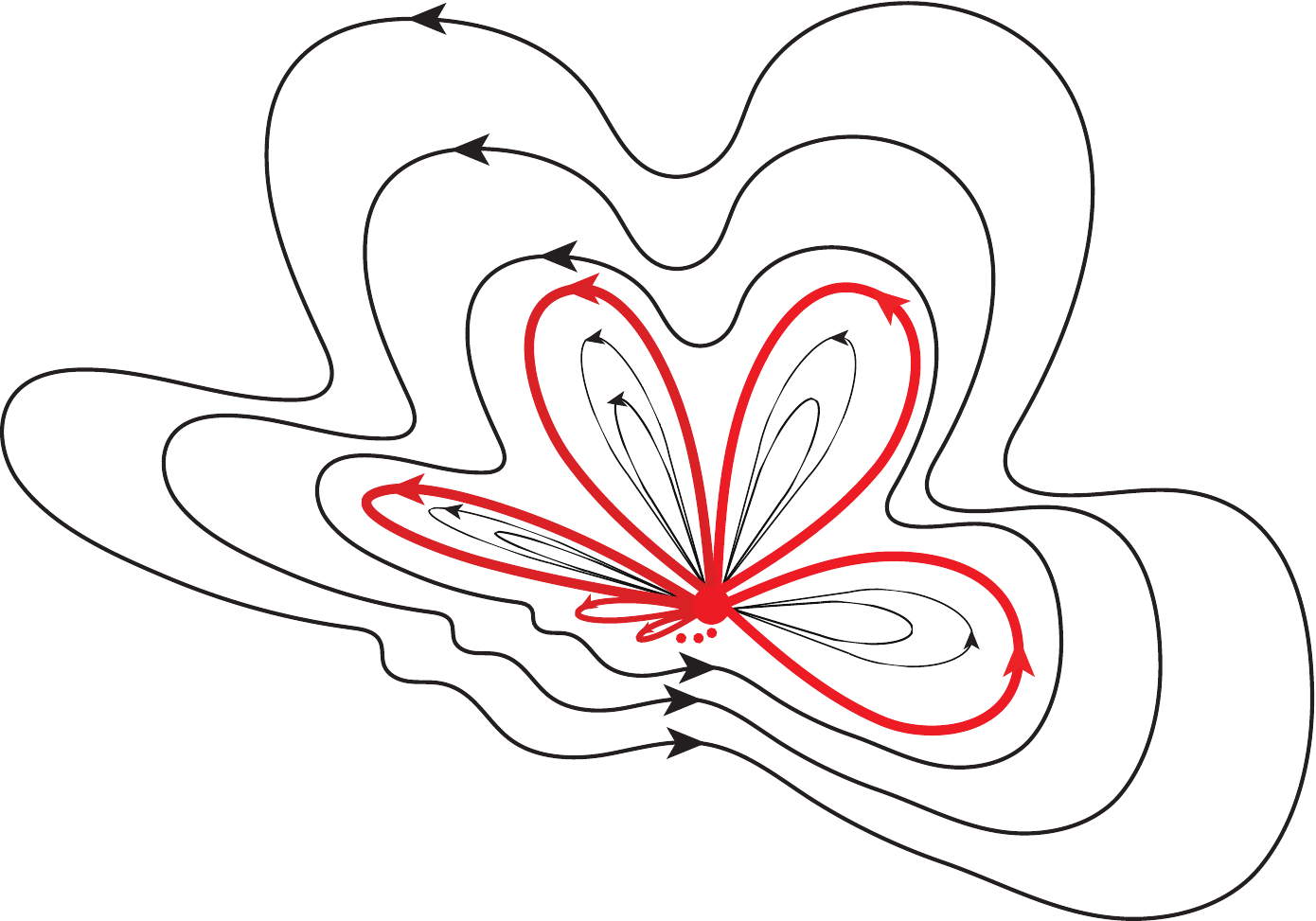}
\end{center}
\caption{A circuit that consists of a degenerate singular point and infinitely many connecting separatrices, and its \nbd which consists of a singular point, periodic orbits, and non-recurrent orbits.}
\label{NAC02}
\end{figure}

A circuit is {\bf trivial} if it is a singular point. 
A circuit is {\bf nontrivial} if it is not trivial. 
%
Note that there are non-trivial circuits with infinitely many edges, and that any non-trivial non-periodic circuit contains non-recurrent orbits as in Figure~\ref{NAC02}.
\begin{definition}
An open annular subset $\mathbb{A}$ of a surface is a {\bf collar} of a non-trivial circuit $\gamma$ if $\gamma$ is a boundary component of $\mathbb{A}$ and there is a neighborhood $U$ of $\gamma$ such that $\mathbb{A}$ is a connected component of the complement $U - \gamma$. 
\end{definition}
A nontrivial circuit $\gamma$ is a {\bf semi-attracting} (resp. {\bf semi-repelling}) circuit with respect to a positive invariant (resp. negative invariant) small collar $\A$ if $\omega(x)  = \gamma$ (resp. $\alpha(x) = \gamma$) for any point $x \in \A$. 
Then $\A$ is called a semi-attracting (resp. a semi-repelling) collar basin of the non-trivial circuit $\gamma$.
A nontrivial circuit $\gamma$ is a {\bf quasi-semi-attracting} (resp. {\bf quasi-semi-repelling}) circuit with respect to a positive invariant (resp. negative invariant) small collar $\A$ if there is a point $x \in \A$ such that $\omega(x)  = \gamma$ (resp. $\alpha(x) = \gamma$). 

\begin{definition}
A cycle is a {\bf limit cycle} if it is a quasi-semi-attracting or quasi-semi-repelling circuit. 
\end{definition}

We will show that any limit cycles have small semi-attracting or semi-repelling collar basins (see Lemma~\ref{lem:limit_cycle}).

\begin{definition}
A non-trivial circuit is a {\bf limit circuit} if it is a quasi-semi-attracting or quasi-semi-repelling circuit. 
\end{definition}

We will show that any limit circuits have small semi-attracting or semi-repelling collar basins under the finiteness of singular points (see Lemma~\ref{non-q-lim-circuit}).

\subsubsection{Transversality}
Notice that we can define transversality using tangent spaces of surfaces because each flow on a compact surface is topologically equivalent to a $C^1$-flow by Gutierrez's smoothing theorem~\cite{gutierrez1986smoothing}.

A $C^1$ simple curve $C$ is {\bf transverse} at a point $p \in S$ to the flow $v_X$ generated by a vector field $X$ if $T_p S = T_p C \oplus T_p O_{v_X}(p)$, where $T_p C$ is the tangent space of $C$ at $p$ and $T_p O_{v_X}(p)$ is the tangent space of the orbit of $p$ with respect to the flow $v_X$. 

\begin{definition}
A simple curve $C$ is {\bf transverse} to $v$ at a point $p \in S$ if there are a vector field $X$ on a surface $T$ and a topologically equivalent homeomorphism $h \colon S \to T$ such that $v$ is topologically equivalent to the flow $v_X$ generated by $X$ via $h$ and that the image $h \circ C$ is a $C^1$ simple curve which is transverse at the point $h(p) \in T$ to the flow $v_X$. 
\end{definition}

A simple curve $C$ is {\bf transverse} to $v$ if so is it at any point in $C$.  
An arc $C$ transverse to $v$ is called a {\bf transverse arc}.
\begin{definition}
A simple closed curve is a {\bf closed transversal} (cf. \cite[Definition~3.4.7 p.41]{Hector1981} and \cite[Definition~3.3.6 p.86]{Candel2000foliation}) if it is transverse to $v$.
\end{definition}

Notice that the closed transversal is one of the fundamental tools in foliation theory to analyze transverse relations among leaves of codimension one foliations (cf. \cite{Hector1981,Candel2000foliation}).  

A simple curve $C$ is {\bf tangent} to $v$ at a point $p \in S$ if it is not transverse at $p$. 
Then the point $p$ is called a tangency of $C$ to $v$. 

\subsubsection{Quasi-circuits}

We introduce a quasi-circuit as follows. 

\begin{definition}
A {\bf quasi-circuit} is one of the following subsets: 
\\
{\rm(1)} A singular point.
\\
{\rm(2)} A cycle. 
\\
{\rm(3)} A closed connected invariant subset which is a boundary component of an open annulus, contains a non-recurrent orbit, and consists of non-recurrent orbits and singular points. 
\end{definition}

A quasi-circuit is {\bf trivial} if it is either a cycle or a singular point. 
A quasi-circuit is {\bf nontrivial} if it is not trivial. 
In case {\rm(3)} in the previous definition, the open annulus is called a {\bf collar} of the nontrivial quasi-circuit.

\begin{definition}\label{def:qst_lqc}
Let $x$ be a point whose $\omega$-limit set is a quasi-circuit which is not a singular point. 
A positive invariant collar $\A$ of $\omega(x)$ is a {\bf quasi-semi-attracting collar basin} of $\omega(x)$ if $O^+(x) \cap \A \neq \emptyset$ and there is a positive invariant collar $\A_{-1}$ of $\omega(x)$ with $\A \subseteq \A_{-1}$ satisfying the following conditions hold: 
\\
\textrm{(1)} The $\omega$-limit set $\omega(x)$ is a boundary component of $\A$
\\
\textrm{(2)} For any non-singular point $y' \in \omega(x)$, there is a transverse closed arc $I_y \subseteq \A_{-1} \sqcup \omega(x)$ whose boundary contains $y'$ such that $I_{y'} - \{y' \} \subset \A$. 
\\
\textrm{(3)} There are a non-singular point $y \in \omega(x)$ and a transverse closed arc $I \subseteq \A_{-1} \sqcup \omega(x)$ whose boundary $\partial I$ consists of $y$ and a point $x_0 \in O^+(x) \cap \A_{-1}$ and which intersects $O^+(x)$ infinitely many times, as shown in Figure~\ref{fig:return_map_01}, 
\begin{figure}
\begin{center}
\includegraphics[scale=0.5]{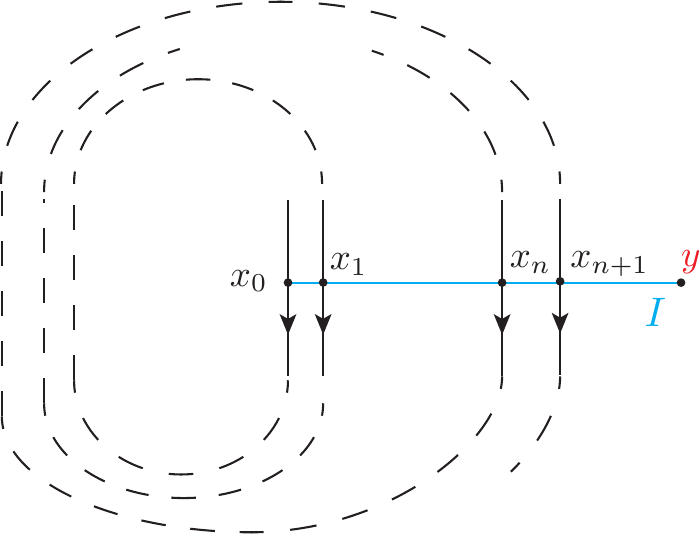}
\end{center}
\caption{A transverse closed arc $I$ and the first return map $f_I$ on $I$.}
\label{fig:return_map_01}
\end{figure}
and there are flow boxes $D_i$ each of whose boundaries $\partial D_i$ is a loop $I_i \cup C_i \cup I_{i+1} \cup C_{i+1}$, as in Figure~\ref{fig:rectangle_01}, 
\begin{figure}
\begin{center}
\includegraphics[scale=0.5]{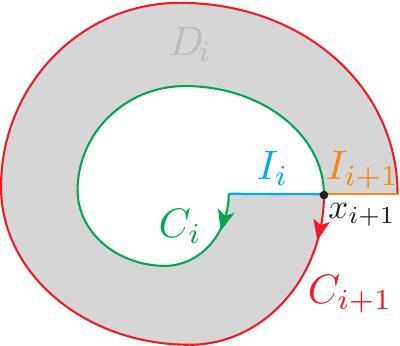}
\end{center}
\caption{A flow box $D_i$ whose boundary is the union of $I_i \cup C_i \cup I_{i+1} \cup C_{i+1}$.}
\label{fig:rectangle_01}
\end{figure}
such that the union $\mathbb{A}_k := (\bigcup_{i =k}^{\infty} (D_i \cup I_{i+1} \cup C_{i+1})) - \{ x_{k+1} \}$ for any $k \in \Z_{\geq 0}$ is a positive invariant open annulus one of whose boundary component is $\omega(x)$, and such that $\A = \A_0$, where $f_I$ is the first return map on $I$, $x_i$ is the $i$-th return image of $x_0$ by $f_I$, $I_{i}$ is the closed sub-arc of $I$ whose boundary consists of $x_{i}$ and $x_{i+1}$, $C_{i}$ is the closed orbit arc in $O^+(x)$ whose boundary consists of $x_{i}$ and $x_{i+1}$.
\end{definition}

\begin{definition}\label{def:qst_lqc_02}
The $\omega$-limit set $\omega(x)$ of a point $x$ is a {\bf quasi-semi-attracting} limit quasi-circuit with respect to a positive invariant small collar $\A_{-1}$ if $\omega(x)$ is a non-trivial quasi-circuit and the collar $\A_{-1}$ contains a quasi-semi-attracting collar basin of $\omega(x)$.
\end{definition}
Using the time reversion, we can define a {\bf quasi-semi-repelling} limit quasi-circuit with respect to a negative invariant small collar $\A$ and its {\bf quasi-semi-repelling collar basin}. 
We introduce the concept of a limit quasi-circuit, which is a generalized concept of a limit circuit. 

\begin{definition}
A non-trivial quasi-circuit is a {\bf limit quasi-circuit} if it is a quasi-semi-attracting or quasi-semi-repelling limit quasi-circuit.
\end{definition}
Note that any circuit with a collar is a quasi-circuit with a collar. 
On the other hand, a quasi-circuit is not a circuit in general (see an example in \S~\ref{ex02}). 

\subsubsection{Quasi-Q-set}
Recall the quasi-Q-set as follows. 
\begin{definition}
The $\omega$-limit (resp. $\alpha$-limit) set of a point is a {\bf quasi-Q-set} if it intersects an essential closed transversal infinitely many times. 
\end{definition}

\begin{figure}
\begin{center}
\includegraphics[scale=0.25]{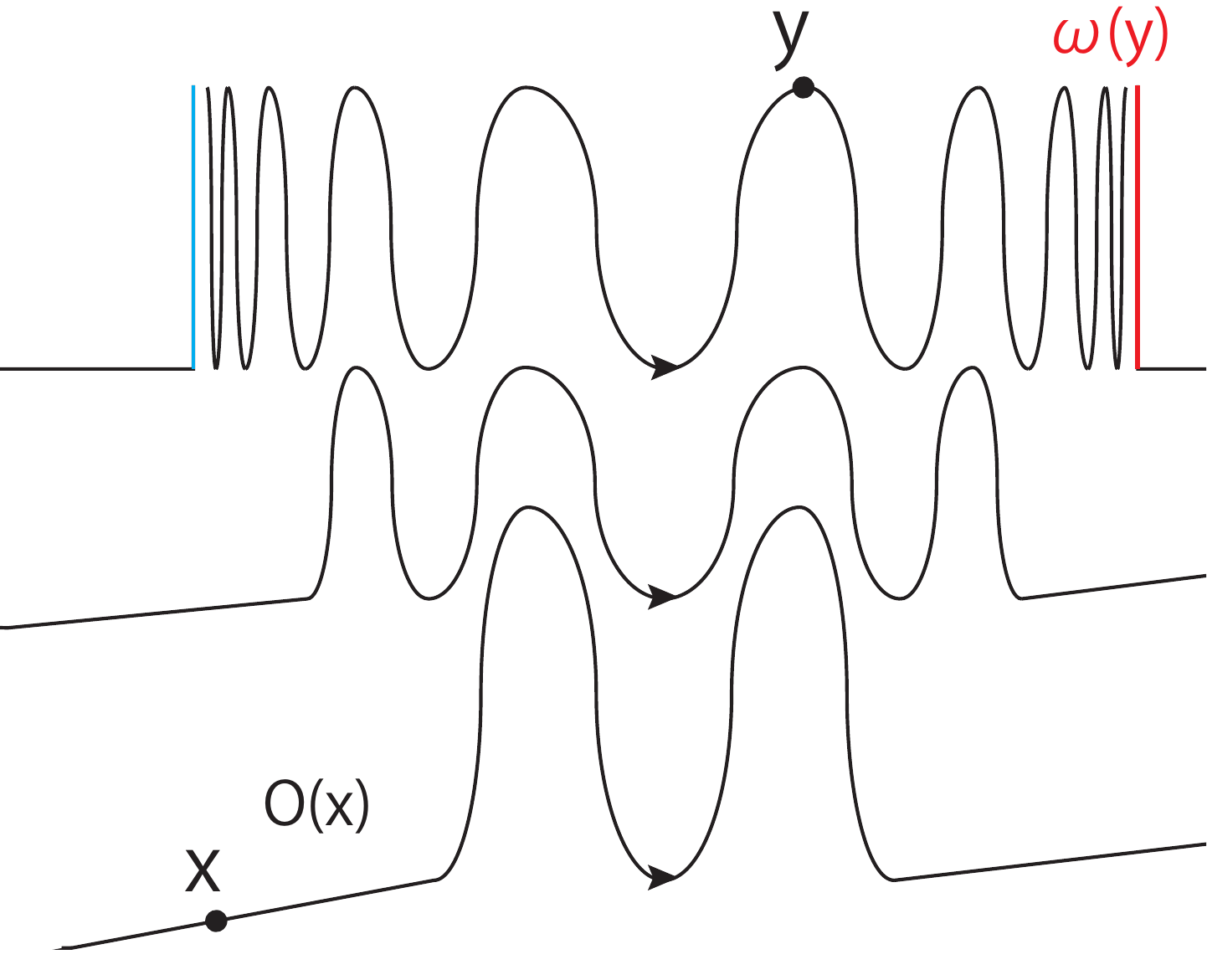}
\end{center}
\caption{A flow box with a non-arcwise-connected connected invariant subset.}
\label{NAC}
\end{figure}

Quasi-Q-sets are topologically characterized in Proposition~\ref{prop:ch_qqset}. 
Note that a quasi-Q-set need not be arcwise-connected. 
In fact, Hastings constructed an attractor of a flow on $\R^2$ which is homeomorphic to a Warsaw circle but is not an $\omega$-limit set \cite[Example 3.3]{Hastings1979PB}. 
By modifying the construction of a Warsaw circle, one can obtain a quasi-Q-set that is not arcwise-connected by replacing a trivial flow box of a point of an exceptional minimal set with a flow box as in Figure~\ref{NAC} (see details \S~\ref{ex00}). 
Moreover, quasi-Q-set is not a Q-set in general (see an example in \S~\ref{ex01}). 
Conversely, we will show that a Q-set is a quasi-Q-set(see Lemma~\ref{preQ}).
On the other hand, if $v$ has at most countably many singular points, then a quasi-Q-set is a Q-set (see Lemma~\ref{q-Q}).
In addition, we have the following observation.

\begin{lemma}\label{lq_q}
Any locally dense quasi-Q-set of a flow on a surface is a Q-set. 
In particular, any quasi-Q-set which intersecting a locally dense orbit is a Q-set. 
\end{lemma}

\begin{proof}
If a quasi-Q-set intersects a locally dense orbit, then the quasi-Q-set is locally dense. 
Fix a locally dense quasi-Q-set $\mathcal{M}$. 
By definition of quasi-Q-set, there is a point $x$ with $\omega(x) = \mathcal{M}$ or $\alpha(x) = \mathcal{M}$. 
By time reversion if necessary, we may assume that $\omega(x) = \mathcal{M}$. 
The locally density implies that $O(x) \cap \mathop{\mathrm{int}} \, \omega(x) = O(x) \cap \mathop{\mathrm{int}} \, \mathcal{M} \neq \emptyset$ and so that 
$x \in \mathop{\mathrm{int}} \,\omega(x)$. 
This means that $x$ is non-closed recurrent and $\overline{O(x)} = \omega(x) = \mathcal{M}$. 
Lemma~\ref{lem:locally_dense_ch} implies that $\mathcal{M} = \overline{O(x)}$ is a locally dense Q-set. 
\end{proof}

We call a quasi-Q-set {\bf non-trivial} if it is not a Q-set. 
Non-trivial quasi-Q-sets are topologically characterized in Proposition~\ref{prop:ch_qqset_nontrivial}. 
Notice that any non-trivial quasi-Q-set of a flow on a compact surface intersects uncountably many connected components of the singular point set (see Lemma~\ref{lem:q-Q_uncountable}).

\section{Poincar\'{e}-Bendixson theorem for surfaces on flows with arbitrarily many singular points}

%
%
%
%
%

\subsection{A generalization of Poincar\'{e}-Bendixson theorem}

This subsection has the following key lemma to generalize the Poincar\'{e}-Bendixson theorem. 

\begin{lemma}\label{lemma:a}
Let $v$ be a flow on a compact surface $S$.
The $\omega$-limit set of any non-closed orbit $O$ is one of the following exclusively:
\\
$(1)$ A nowhere dense subset of singular points.
\\
$(2)$ A semi-attracting limit cycle.
\\
$(3)$ A quasi-semi-attracting limit quasi-circuit.
\\
$(4)$ A locally dense Q-set.
\\
$(5)$ A quasi-Q-set that is not locally dense.

Moreover, in case $(3)$, the orbit $O$ is wandering, and $O \cap \overline{O'} = \emptyset$ for any orbit $O' \neq O$. 
\end{lemma}

To show the previous lemma, we show some technical lemmas. 
The proof methods use techniques derived from foliation theory (cf. \cite{Hector1981,Candel2000foliation}).
First, we show the existence of closed transversals near infinite intersections of transverse arcs and orbits. 

\begin{lemma}\label{loops}
Let $I$ be a transverse arc and $x \in I$ such that $\vert  I \cap O(x)\vert  = \infty$. 
Then there are an orbit arc $C$ in $O(x)$ and a transverse closed arc $J \subseteq I$ such that the union $\mu := J \cup C$ is a loop with $C \cap J = \partial C = \partial J$ and that the return map along $C$ is orientation-preserving between neighborhoods of $\partial C$ in $I$. 
Moreover, for any small number $\varepsilon > 0$, there is a closed transversal $\gamma \subset B_{\varepsilon}(\mu)$, where $B_{\varepsilon}(\mu) := \{ y \in S \mid \min_{z \in \mu} d(y,z) < \varepsilon \}$ is the $\varepsilon$-\nbd of $\mu$ with respect to the Riemannian distance for a Riemannian metric on $S$. 
\end{lemma}

\begin{proof}
By time reversion if necessary, we may assume that $O^+(x) \cap \mathop{\mathrm{int}} I = \infty$. 
Fix a point $x_0 \in O^+(x) \cap \mathop{\mathrm{int}} I$. 
Let $f_v: I' \to I$ be the first return map on $I$ induced by $v$ with the maximal domain $I' \subseteq I$, $x_i := (f_v)^i(x_0)$ the $i$-th return of $x_0$, $C_{a,b} \subset O^+(x)$ the orbit arc from $a$ to $b$, and $I_{a,b} \subset I$ the subinterval between $a$ and $b$ of $I$.  
We may assume that $x_0 < x_1$. 

Suppose that the restriction of $f_v$ to a neighborhood of $x_i$ for some $i \in \Z_{\geq 0}$ is orientation-preserving. 
Then put $C := C_{x_i, x_{i+1}}$ and $J := I_{x_i, x_{i+1}}$. 
By the waterfall construction (cf. \cite[Lemma~3.3.7 p.86]{Candel2000foliation}) to the loop $\mu := C \cup J$ (see Figure~\ref{wf}), there is a closed transversal $\gamma$ intersecting $O(x)$ near $\mu$. 
\begin{figure}
\begin{center}
\includegraphics[scale=0.4]{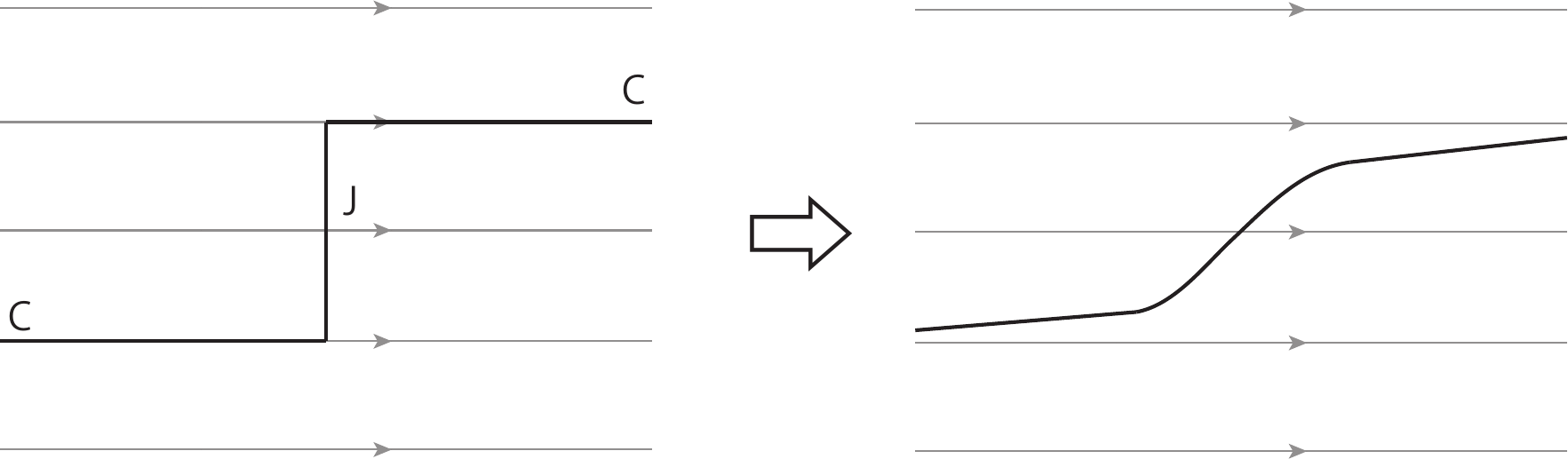}
\end{center}
\caption{The waterfall construction}
\label{wf}
\end{figure}
Thus we may assume that that the restriction of $f_v$ to a neighborhood of $x_i$ for any $i \in \Z_{\geq 0}$ is orientation-reversing. 

We claim that there is a natural number $i$ such that $x_{i+1} < x_i$. 
Indeed, otherwise $x_i < x_{i+1}$ for any $i \in \Z_{\geq 0}$. 
Then each pair of loops $\gamma_i := C_{x_{2i}, x_{2i+1}} \cup I_{x_{2i+1}, x_{2i}}$ has disjoint neighborhoods each of which is a M\"obius band. 
This means that $S$ has infinite non-orientable genus, which contradicts the compactness of $S$.

By renumbering, we may assume that $x_2 < x_1$.  
From $x_0 < x_1$, the first return map for $I_{x_2, x_0}$ along $C_{x_0, x_2}$ is orientation-preserving. 
Put $C := C_{x_0, x_2}$ and $J := I_{x_2, x_0}$. 
As above, the waterfall construction to the loop $\mu := C \cup J$ completes the assertion. 
\end{proof}

We show that the infinite intersection of a transverse closed arc implies the existence of a quasi-circuit under properness. 

\begin{lemma}\label{lem-quasi-limit-cricuits}
If there is a transverse closed arc $J: [-1,0] \to S$ with $\{ J(0) \} = J([-1,0]) \cap \omega(J(-1)) \subset \overline{J([-1,0]) \cap O^+(J(-1))}$, then the following properties hold: 
\\
{\rm(1)} The $\omega$-limit set $\omega(J(-1))$ is either a limit cycle or a limit quasi-circuit. 
\\
{\rm(2)} The point $J(-1)$ is wandering and $J(-1) \not\in \overline{O(p)}$ for any point $p \in S - O(J(-1))$. 

Moreover, we can choose a small quasi-semi-attracting collar basin $\A$ of $\omega(J(-1))$ such that $\partial \A - \omega(J(-1))$ is a loop consisting of a closed orbit arc and a transverse closed interval.
\end{lemma}

\begin{proof}
Let $J: [-1,0] \to S$ be a transverse closed arc, $x_0 := J(-1)$ a point, $y := J(0)$ a point, and $I := J([-1,0])$ a closed interval with $\{ y \} = I \cap \omega(x_0) \subset \overline{I \cap O^+(x_0)}$, as shown in Figure~\ref{fig:return_map}.
\begin{figure}
\begin{center}
\includegraphics[scale=0.5]{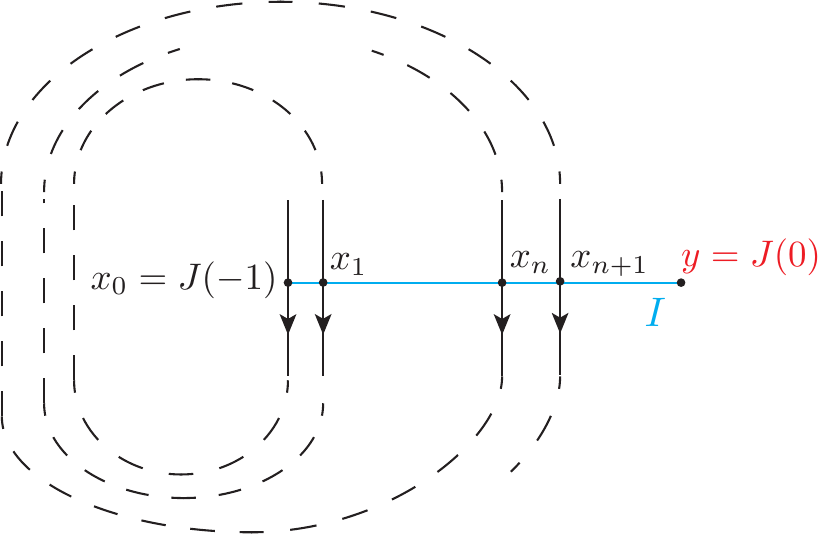}
\end{center}
\caption{A return map on $I$.}
\label{fig:return_map}
\end{figure}
Since $y \in \overline{I \cap O^+(x_0)} \cap \omega(x_0)$, we obtain that $O^+(x_0) \cap I$ is infinite and  $\{ y \} =  I \cap \omega(x_0) \subset \partial I = \{ x_0, y \}$. 
Denote by $x_i$ the $i$-th return image of $x_0$ on $I$.
By $x_i \not\in \{ J(0) \} = J([-1,0]) \cap \omega(x_0)$ for any $i \in \Z_{\geq 0}$, the point $x_i$ for any $i \in \Z_{> 0}$ is not positively recurrent and so is $x_0$. 
Since $\{ y \} =  I \cap \omega(x_0)$, any closed interval contained in $I - \{y\}$ intersects $O^+(x_0)$ at most finitely many times.

\begin{claim}
The sequence $(x_i)_{i \in \Z_{\geq N }}$ in $I$ is strictly increasing and converges to $y$ for some $N > 0$.  
\end{claim}
\begin{proof}
It suffices to show that there is a large number $N > 0$ such that $x_{i} < x_{i+1}$ in $I$ for any natural number $i \geq N$ with respect to a natural total order on the sub-arc $I$. 
Otherwise $x_{i} > x_{i+1}$ holds for infinitely many natural numbers $i$.
Since each closed sub-arc of $\mathop{\mathrm{int}} I$ intersects at most finitely many points of $O^+(x_0)$, there are infinitely many triples $i_{k}^-:= i_{k}^0-1 <  i_{k}^0 <  i_{k}^0+1:=i_{k}^+$ of natural numbers with $i_{k}^+ < i_{k+1}^-$ such that either $x_{i_{k}^-} < x_{i_{k}^+} < x_{i_{k}^0} < x_{i}$ or $x_{i_{k}^+} < x_{i_{k}^-} < x_{i_{k}^0} < x_{i}$  for any $i \geq i_{k+1}^-$ as in Figure~\ref{fig:return_map02}.
Denote by $I_{i_{k}^{-0}}$ (resp. $I_{i_{k}^{0+}}$) by the sub-arc of $I$ whose boundary consists of $x_{i_{k}^-}$ and $x_{i_{k}^0}$ (resp. $x_{i_{k}^0}$ and $x_{i_{k}^+}$), and by $C_{i_{k}^{-0}}$ (resp. $C_{i_{k}^{0+}}$) the curve contained in $O^+(x_0)$ whose boundary consists of $x_{i_{k}^-}$ and $x_{i_{k}^0}$ (resp.  $x_{i_{k}^0}$ and $x_{i_{k}^+}$).
Then the unions $T_{i_{k}^{-0}} := C_{i_{k}^{-0}} \cup I_{i_{k}^{-0}}$ and $T_{i_{k}^{0+}} := C_{i_{k}^{0+}} \cup I_{i_{k}^{0+}}$ are simple closed curves whose intersection is a closed arc $I_{i_{k}^{0+}}$ or $I_{i_{k}^{-0}}$ as on the upper of Figure~\ref{fig:return_map02}.
\begin{figure}
\begin{center}
\includegraphics[scale=0.48]{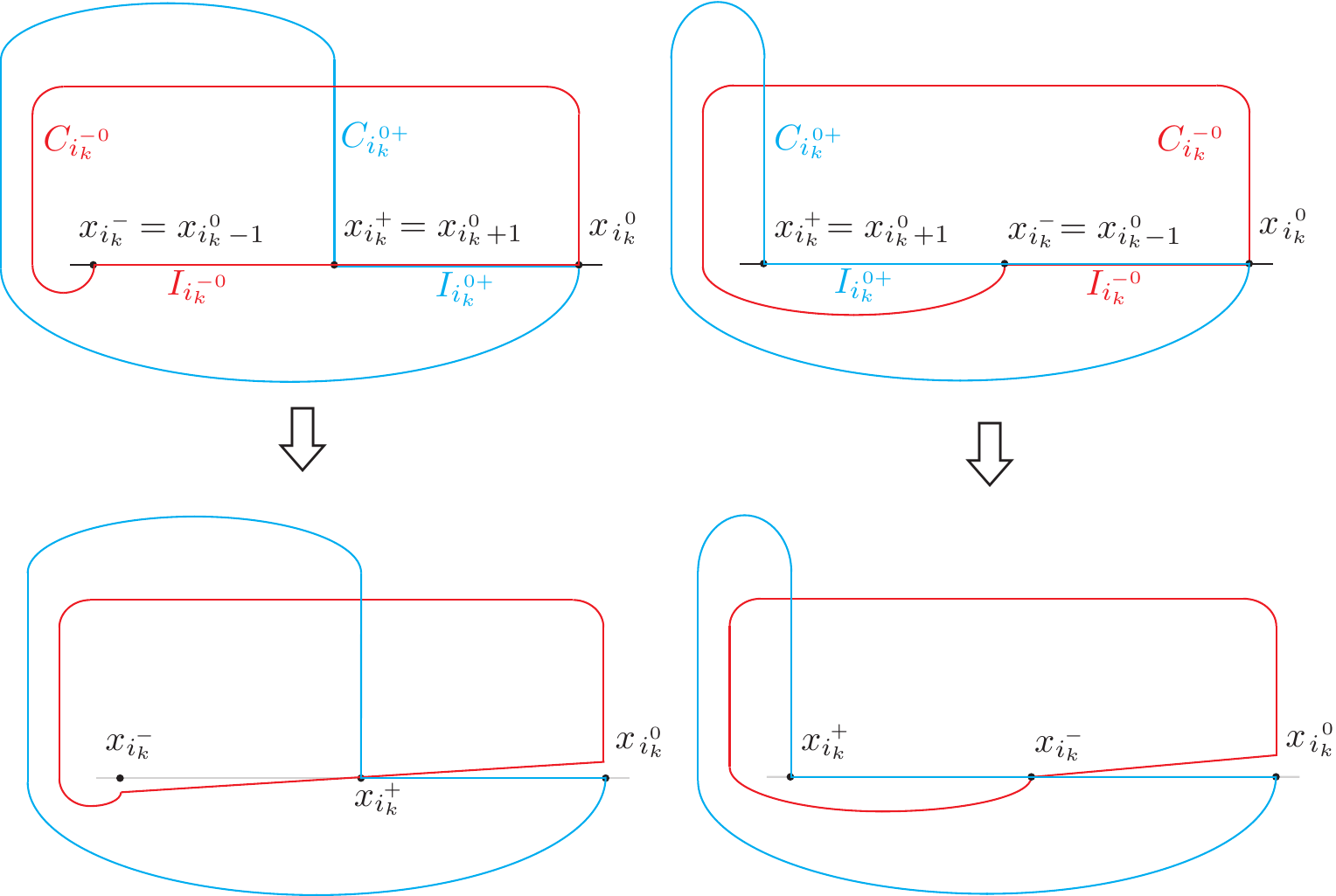}
\end{center}
\caption{(Case where $(x_i)$ is non-increasing) Upper, Two simple closed curves $T_{i_{k}^{-0}} = C_{i_{k}^{-0}} \cup I_{i_{k}^{-0}}$ and $T_{i_{k}^{0+}} = C_{i_{k}^{0+}} \cup I_{i_{k}^{0+}}$; lower, deformed two transverse simple closed curves $T'_{i_{k}^{-0}}$ and $T'_{i_{k}^{0+}}$.}
\label{fig:return_map02}
\end{figure}
By a deformation like a waterfall construction as on the lower of Figure~\ref{fig:return_map02}, we obtain two simple closed curves $T'_{i_{k}^{-0}}$ and $T'_{i_{k}^{0+}}$ whose intersection is either $x_{i_{k}^-}$ or $x_{i_{k}^+}$ and which are close to the original simple closed curves $T_{i_{k}^{-0}}$ and $T_{i_{k}^{0+}}$ respectively with respect to the Riemannian distance for a Riemannian metric on $S$.
(Note that if $S$ is orientable then we can choose $T'_{i_{k}^{-0}}$ and $T'_{i_{k}^{0+}}$ to be closed transversals.)
Since $\max \{x_{i_{k}^-},x_{i_{k}^0},x_{i_{k}^+} \} < x_{i}$  for any $i \geq i_{k+1}^-$, we have $O^+(x_{i_{k+1}^-}) \cap (I_{i_{k}^{-0}} \cup I_{i_{k}^{0+}}) = \emptyset$ and so $O^+(x_{i_{k+1}^-}) \cap (T'_{i_{k}^{-0}} \cup T'_{i_{k}^{0+}}) = \emptyset$.
Hence $(T'_{i_{k}^{-0}} \cup T'_{i_{k}^{0+}}) \cap (T'_{i_{l}^{-0}} \cup T'_{i_{l}^{0+}}) = \emptyset$ for any $k \neq l$.
Since these simple closed curves intersect at exactly one point, they are essential.
Cutting $T'_{i_{1}^{0+}}$ and collapsing new boundary components into singletons, we obtain the resulting surface whose genus is the genus of $S$ minus one.
Since there are infinitely many disjoint bouquets $T'_{i_{k}^{-0}} \cup T'_{i_{k}^{0+}}$,
the genus of $S$ is not finite, which contradicts the compactness of $S$.
Thus there is a large number $N > 0$ such that $x_{i} < x_{i+1}$ in $I$ for any natural number $i \geq N$. 
\end{proof}

For any $i \in \Z_{\geq N}$, denote by $I_{i}$ by the sub-arc of $I$ whose boundary consists of $x_{i}$ and $x_{i+1}$, and by $C_{i}$ the curve contained in $O^+(x_0)$ whose boundary consists of $x_{i}$ and $x_{i+1}$.
Fix a connected component $D_i$ of $S - \bigcup_{j =N}^\infty (I_j \cup C_j)$ one of whose boundary components is the union $I_i \cup C_i \cup I_{i+1} \cup C_{i+1}$.
Since there are at most finite genus, by renumbering $x_i$, we may assume that $D_i$ is a rectangle for any $i \in \Z_{\geq N}$ as in Figure~\ref{fig:rectangle_01}.  
Define a union $\mathbb{A}_k := (\bigcup_{i =k}^{\infty} (D_i \cup I_{i+1} \cup C_{i+1})) - \{ x_{k+1} \}$   for any $k \in \Z_{\geq N}$, which is an open annulus. 
By the monotonicity of $x_i$ in $I$, each $\mathbb{A}_k$ is a positive invariant open annulus homotopic to $\mathbb{A}_N$ for any $k \in \Z_{\geq N}$.
Since $\omega(x_0)$ is connected, by construction, the closure $\overline{D_i}$ of any rectangle $D_i$ does not intersect $\omega(x_0)$ and so $\mathbb{A}_k \cap \omega(x_0) = \emptyset$.
Therefore $\omega(x_0) = \bigcap_{n \in \mathbb{R}_{\geq N}}\overline{\{v_t(x_0) \mid t > n\}} = \bigcap_{k \in \mathbb{Z}_{\geq N}}\overline{\{v_t(x_0) \mid t > k\}} \subseteq \bigcap_{k \in \mathbb{Z}_{\geq N}}\overline{\mathbb{A}_k} = \partial \overline{\mathbb{A}_N} - (I_N \cup C_N)$.
This means that $\omega(x_0)$ is contained in a boundary component $\omega$ of the annuli $\mathbb{A}_N$ and so $\mathbb{A}_k$ for any $k \in \mathbb{Z}_{\geq N}$.
Moreover, the boundary component of $\partial \mathbb{A}_N$ which does not intersect $\omega$ is a loop consisting of the orbit arc $C_N$ of $O(x)$ and the transverse closed arc $I_N$. 

\begin{claim}\label{claim:bd_ann}
$\omega(x_0)$ is the boundary component $\omega$ of the annulus $\mathbb{A}_N$. 
\end{claim}
\begin{proof}
By the previous claim, from the existence of a closed flow box containing $I$, there is a positive number $\varepsilon > 0$ such that the length of any $C^1$-arc in $\mathbb{A}_N - O^+(x_0)$ from a point in $I_i \setminus O^+(x_0)$ to a point in $I_{i+1} \setminus O^+(x_0)$ is at least $\varepsilon$. 

Assume that there is a point $x \in \omega - \omega(x_0) \subset \partial \mathbb{A}_N$. 
Then there is a small \nbd $U_x$ of $x$ with $U_x \cap (O^+(x_0) \cup I_N \cup C_N) = \emptyset$. 
Since $x \in \omega = \partial  \mathbb{A}_N - (I_N \cup C_N)$, there is a point $a \in U_x \cap \mathbb{A}_N$. 
Take a closed $C^1$-arc $\gamma \subset U_x$ from $x$ to $a$ which has a finite length. 
Since $U_x \cap (O^+(x_0) \cup I_N \cup C_N) = \emptyset$, we have $\gamma \cap (O^+(x_0) \cup I_N \cup C_N) = \emptyset$. 
By $\mathbb{A}_N = (\bigcup_{i =N}^{\infty} (D_i \cup I_{i+1} \cup C_{i+1})) - \{ x_{N+1} \}$, 
we obtain $\mathbb{A}_N \setminus O^+(x_0) \subseteq (\bigcup_{i =N}^{\infty} (D_i \cup I_{i+1})) \setminus O^+(x_0)$. 
There is an integer $N' \geq N$ such that $a \in D_{N'} \cup I_{N'+1}$. 
Since  $\gamma \subset U_x \subset S - (O^+(x_0) \cup (\partial \mathbb{A}_N - \omega))$ is a closed arc from $x \in \partial \mathbb{A}_N = \overline{\mathbb{A}_N} - \mathbb{A}_N$ to $a \in D_{N'} \cup I_{N'+1}$, the closed arc $\gamma$ intersects $I_n \setminus O^+(x_0)$ for any $n > N'$. 
Since $\gamma$ contains a $C^1$-arc in $\mathbb{A}_N - O^+(x_0)$ from a point in $I_{N'+1} \setminus O^+(x_0)$ to a point in $I_{N'+1+k} \setminus O^+(x_0)$ for any $k \in \Z_{>0}$ whose length is at least $k \varepsilon$, the length of $\gamma$ is infinite, which contradicts the finite length of $\gamma$. 
\end{proof}

\begin{claim}
For any non-singular point $y' \in \omega(x)$, there is a transverse closed arc $I_{y'} \subseteq \A_{N} \sqcup \omega(x)$ whose boundary contains $y'$ and which intersects $O^+(x)$ infinitely many times such that $I_{y'} - \{y'\} \subset \A_N$. 
\end{claim}
\begin{proof}
By $y' \in \omega(x)$, since $y'$ is non-singular, Claim~\ref{claim:bd_ann} implies that there is a transverse closed arc $I_{y'} \subseteq \A_{N} \sqcup \omega(x)$ whose boundary contains $y'$ and which intersects $O^+(x)$ infinitely many times such that $I_{y'} - \{ y'\} \subset \A_N$. 
\end{proof}

Therefore $\omega(x_0)$ is either a limit cycle or a limit quasi-circuit with its quasi-semi-attracting collar basin $\A := \A_N$ such that  $\partial \mathbb{A} - \omega(x_0) = C_N \cup I_N$ .

\begin{claim}
The point $x_0$ is wandering. 
\end{claim}
\begin{proof}
It suffices to show that $x_{N+1}$ is wandering. 
By definition, we have $x_{N+1} \in \mathop{\mathrm{int}} (I_N \cup I_{N+1})$ and $\partial I_{N+k} = \{ x_{N+k}, x_{N+k+1} \}$ for any $k \in \Z_{\geq 0}$.  
Let $f_v \colon I' \to I$ be the first return map on $I$ with the maximal domain $I' \subseteq I$. 
Since the point $x_k$ for any $k \in \Z_{\geq 0}$ is contained in the domain of $f_v$, the flow box theorem (cf. \cite[Theorem 1.1, p.45]{aranson1996introduction}) implies that there is a small closed interval $J$ whose interior contains $x_{N+1}$ such that $J$ and $f_v(J)$ are contained in the domain of $f_v$.  
Then $f_v^2(J)$ is contained in $\mathbb{A}_{N+2}$. 
Taking $J$ short, we may assume that there is a positive number $T>0$ such that $x_{N+3} \in v_T(J) \subset \mathbb{A}_{N+2}$. 
Then there is a small positive number $\varepsilon >0$ such that $v_T(\bigsqcup_{t \in (-\varepsilon, \varepsilon)} v_t(J)) \subset \mathbb{A}_{N+2}$. 
The open subset $U := \bigsqcup_{t \in (-\varepsilon, \varepsilon)} v_t(J) \subset \mathbb{A}_{N-1} - \mathbb{A}_{N+2}$ is an open \nbd of $x_{N+1}$ and $v_T(U) \subset \mathbb{A}_{N+2}$. 
By the positive invariance of $\mathbb{A}_{N+2}$, we have $v_t(U) \subset \mathbb{A}_{N+2}$ and so $U \cap v_t(U) = \emptyset$ for any $t >T$. 
This means that $x_{N+1}$ is wandering. 
\end{proof}

Since any wandering point $w$ is not contained in the orbit closure of points outside of $O(w)$, we obtain $x_0 \not\in \overline{O(z)}$ for any point $z \in S - O(x_0)$. 
\end{proof}

We have the following observations. 

\begin{lemma}\label{lem:limit_cycle}
The $\omega$-limit set of a non-periodic point $x$ that intersects periodic points is a limit cycle with its semi-attracting collar basin $\A$. 
\end{lemma}

\begin{proof}
Fix a periodic point $y \in \omega(x)$. 
The flow box theorem implies that the limit cycle $O(y)$ is covered by finitely many trivial flow boxes $B_1, \ldots , B_k$ with $y \in B_1$. 
By $O(y)  \subseteq \Pv \cap \omega(x)$, there is a small transverse closed arc $J: [-1,0] \to B_1 \subset S$ with $J(-1) \in O^+(x)$, $J(0) = y$, and $\bigcup_{z \in J([-1,0])} O^+(z) \subset \bigcup_{i=1}^k B_i$ such that the first return map $f_v \colon J([-1,0]) \to J([-1,0])$ to the transverse closed arc $J([-1,0])$ is an attracting map with $\bigcap_{n \in \Z_{>0}} f_v^n(J([-1,0])) = \{ y \}$. 
Then the union $\A := \bigcup_{z \in J((-1,0))} O^+(z) \subset \bigcup_{i=1}^k B_i$ is a semi-attracting collar basin of $\omega(x)$ such that $\bigcap_{t > 0} v(t, \A) = O(y) \subseteq \Pv \cap \bigcup_{i=1}^k B_i$. 
Therefore $\omega(x') = O(y) = \omega(x)$ for any $x' \in \A$. 
\end{proof}

\begin{corollary}\label{cor:limit_cycle}
Any locally dense Q-sets and any quasi-Q-sets intersect no periodic points.
\end{corollary}

\begin{proof}
The properness of periodic orbits implies that each periodic orbit intersects any closed transversal at most finitely many times. 
If a quasi-Q-set $Q$ intersects periodic orbits, then Lemma~\ref{lem:limit_cycle} implies that $Q$ is a periodic orbit that intersects any closed transversal at most finitely many times, which contradicts the infinite intersection of a closed transversal and $Q$. 
If a locally dense Q-set $Q$ intersects periodic points, then Lemma~\ref{lem:limit_cycle} implies that $Q$ is a periodic orbit, which contradicts that $Q$ contains non-closed recurrent points. 
\end{proof}

We show that the infinite intersection of an essential closed transversal implies the existence of either a quasi-Q-set, an essential limit cycle, or a quasi-circuit. 

\begin{lemma}\label{lem4-07}
Let $x$ be a point contained in a closed transversal $\gamma$ such that $O^+(x)$ intersects $\gamma$ infinitely many times.
Then $\omega(x)$ is either a quasi-Q-set, an essential limit cycle, or a limit quasi-circuit. 
Moreover, if $\omega(x)$ is a quasi-circuit, then $x$ is wandering and is not contained in the orbit closures of points outside of $O(x)$.
\end{lemma}

\begin{proof}
Suppose that a positive orbit $O^+(x)$ intersects a closed transversal $\gamma$ infinitely many times. 
Then $x$ is not periodic. 
Since any orbit and any inessential closed transversal intersects at most once, the closed transverse $\gamma$ is essential such that $\gamma \cap \omega(x) \neq \emptyset$ and so that $\omega(x) \not\subseteq \mathop{\mathrm{Sing}}(v)$. 
If $\omega(x)$ contains periodic orbits, then Lemma~\ref{lem:limit_cycle} implies that it is an essential limit cycle because the positive orbit $O^+(x)$ intersects the basin of the semi-attracting limit cycle.
Thus we may assume that $\omega(x) \cap \mathop{\mathrm{Per}}(v) = \emptyset$.
If $\omega(x)$ is a quasi-Q-set, then the assertion holds. 
Thus we may assume that $\omega(x)$ is not a quasi-Q-set.
Then $\gamma \cap \omega(x)$ is nonempty and finite. 

We claim that there is a transverse closed arc $J:[-1, 0] \to \gamma$ with $J(-1 ) \in O^+(x)$ and $\{ J(0) \} = J([-1,0]) \cap \omega(J(-1)) \subset \overline{O^+(J(-1)) \cap J([-1,0])}$. 
Indeed, fix a point $y \in \gamma \cap \omega(x)$. 
Since $\gamma \cap \omega(x)$ is finite, there are a point $x_0 \in \gamma \cap O^+(x)$ and a transverse closed arc $J:[-1, 0] \to \gamma$ with $x_0 = J(-1)$ and $y = J(0)$ such that $\{ y \}  = J([-1,0]) \cap \omega(x_0) \subset \overline{O^+(x_0) \cap J([-1,0])}$. 

By Lemma~\ref{lem-quasi-limit-cricuits}, the $\omega$-limit set $\omega(x)$ is a limit quasi-circuit such that $x$ is wandering and is not contained in the orbit closures of points outside of $O(x)$.
\end{proof}

We have the following observation. 

\begin{lemma}\label{lem:transverse_bdry}
Let $\A$ be a quasi-semi-attracting collar basin of a quasi-semi-attracting limit quasi-circuit $\omega(x)$.
Then every $C^1$ transverse closed arc whose interior is contained in $\A$ and whose boundary is contained in $\partial \A$ and contains a point in $\omega(x)$ connects the boundary components $\omega$ and $\A - \omega(x)$. 
\end{lemma}

\begin{proof}
Let $\gamma$ be a $C^1$ transverse closed arc with $\mathop{\mathrm{int}}\gamma \subset \A$ such that $\partial \gamma$ contains a point $y' \in \omega(x)$. 
Fix any Riemannian metric on $S$ which induces the Riemannian distance.
Then the length of $\gamma$ is finite. 
By renumbering, we may assume that $\mathbb{A} := (\bigcup_{i =0}^{\infty} (D_i \cup I_{i+1} \cup C_{i+1})) - \{ x_{1} \}$, where $D_j$, $C_j$, and $x_j$ as in Definition~\ref{def:qst_lqc}. 
Let $I$ be a transverse closed arc as in Definition~\ref{def:qst_lqc}.
Replacing $\A$ with $\A_k$ for some large $k$, we may assume that the set difference $D := \A \setminus I$ is a rectangle. 

\begin{claim}
We may assume that there is a closed sub-arc $J$ of $\gamma$ from $y' \in \omega(x)$ with $J - \{ y' \} \subset D$ and $J \cap I = \emptyset$, by deforming $\gamma$ near $y'$ in $\A$.  
\end{claim}
\begin{proof}
From the existence of a closed trivial flow box containing $I$, there is a positive number $\varepsilon > 0$ such that, for any $i  \in \Z_{\geq 0}$, the length of any $C^1$-arc in $\A - O^+(x)$ from a point in $I_i \setminus O^+(x)$ to a point in $I_{i+1} \setminus O^+(x)$ is at least $\varepsilon$. 
Considering the universal covering of $\A$, the finite length of $\gamma$ implies that the lift $\widetilde{\gamma}$ of $\gamma \cap \A$ intersect at most finitely many lifts of $I \cap \A$. 
Therefore the existence of a closed trivial flow box containing $I$ implies that we may assume that $\gamma$ contains a closed sub-arc $J$ from $y' \in \omega(x)$ with $J - \{ y' \} \subset \A \setminus I = D$ and $J \cap I = \emptyset$, by perturbing $\gamma$ into a $C^1$ transverse closed arc with $\mathop{\mathrm{int}}\gamma \subset \A$. 
\end{proof}


\begin{claim}
The sub-arc $J$ intersects the open interval $C_j$ at most once for any $j \in \Z_{\geq 0}$. 
\end{claim}
\begin{proof}
Assume that $J \subset D$ intersects the open interval $C_j$ at least twice for some $j > N$. 
Since $D - C_i$ is the disjoint union of two open disks, the orientability of the open disk $D$ implies the incompatibility of the direction of the orbit arc $C_j$ at a pair of points $c_j, c'_j \in C_j \cap \gamma$ as in Figure~\ref{fig:orientation_reversion_loop}, which is a contradiction. 
\begin{figure}
\begin{center}
\includegraphics[scale=0.5]{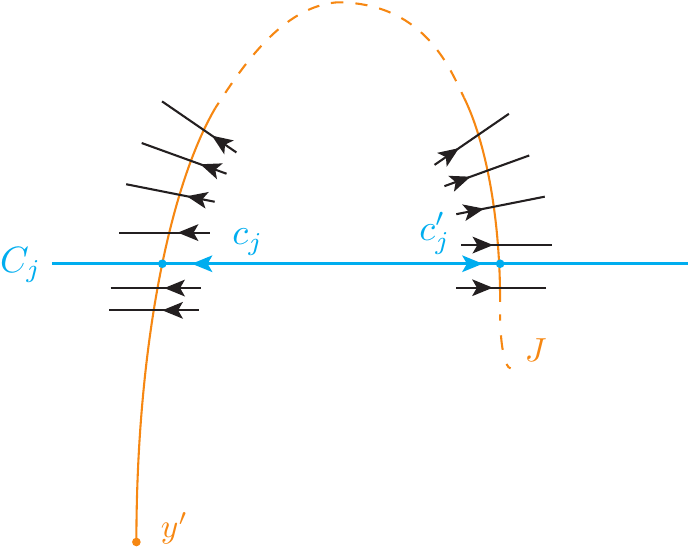}
\end{center}
\caption{An open sub-arc $J$ of $\gamma$ from $y' \in \omega(x)$ to a point in $B_r(y') \cap \A \subset D$.}
\label{fig:orientation_reversion_loop}
\end{figure}
\end{proof}

From $J \cap \mathop{\mathrm{int}}\A \neq \emptyset$, there is a point $z \in J \cap (D_{j_0} \cup C_{j_0})$ for some $j_0  \in \Z_{\geq 0}$. 
By the previous claim, we have $\vert J \cap C_j \vert = 1$ for any $j > j_0$. 
Therefore $\vert \gamma \cap C_j \vert \geq 1$ for any $j > j_0$. 
Moreover, the set difference $D' := \A_{j_0 +1} \setminus J \subset \A \setminus J$ is a rectangle, where $\mathbb{A}_{j_0 +1} := (\bigcup_{i =j_0 +1}^{\infty} (D_i \cup I_{i+1} \cup C_{i+1})) - \{ x_{j_0 +2} \}$. 

\begin{claim}
The transverse closed arc $\gamma$ intersects the open interval $C_j$ exactly once for any $j > j_0 + 1$.
\end{claim}
\begin{proof}
Assume that $\gamma \subset \A - \bigcup_{n \geq j_0} I_n$ intersects the open interval $C_j$ at least twice for some $j > j_0 + 1$.
The simplicity of $\gamma$ implies that $\gamma - J \subset \A \setminus J$. 
Since $\A  - (C_j \cup J)$ is the disjoint union of an open annuls and an open disk, the orientability of the annulus $\A$ implies the incompatibility of the direction of the orbit arc $C_j$ at a pair of points $c_j, c'_j \in C_j \cap \gamma$ as in Figure~\ref{fig:orientation_reversion_loop}, which is a contradiction. 
\end{proof}
The previous claim implies the assertion. 
\end{proof}

We state the non-existence of transversely accumulation non-singular points in limit quasi-circuits as follows. 

\begin{lemma}\label{lem:acc_lc}
For any limit quasi-circuit $C$, there is no transverse closed arc $\gamma$ such that $C \cap \gamma$ has accumulation points. 
\end{lemma}

\begin{proof}
Let $C$ be a limit quasi-circuit.  
By time reversion if necessary, we may assume that there is a point $x$ with $\omega(x) = C$. 
Since $\omega(x)$ is a limit quasi-circuit, there is a small open annulus $\mathbb{A}$ which is a quasi-semi-attracting collar basin and of which $\omega(x)$ is a boundary component such that the boundary $\partial \mathbb{A}$ consists of two connected components. 
Denote by $\partial_1 := \partial \A - \omega(x)$ another boundary component of $\mathbb{A}$.
In other words, we have $\partial \mathbb{A} = \omega(x) \sqcup \partial_1$.
%
Fix any distance function $d$ on $S$ induced by a Riemannian metric.
Since the boundary components of $\mathbb{A}$ are compact and disjoint, there is a positive number $d_{\mathbb{A}}$ such that $d_{\mathbb{A}} = \min \{ d(y, z) \mid y \in \omega(x), z \in \partial_1 \}$.

Assume that there is a transverse closed arc $\gamma$ such that $C \cap \gamma$ has accumulation points $x_\infty \in \overline{\gamma \cap \omega(x)}$. 
Extending the transverse closed arc $\gamma$ if necessary, we may assume that $x_\infty \in  \mathop{\mathrm{int}} \gamma$. 
By definition of transverse, by taking a topologically equivalent homeomorphism if necessary, we may assume that $\gamma$ is $C^1$. 
From the flow box theorem applying to the compact subset $\gamma$, the fact that $\gamma$ contains no singular point implies that there are a trivial flow box $U$ with $\gamma \subset U
$ and a homeomorphism $h \colon [0,1]^2 \to U$ such that the images $h(\{p_1 \} \times [0,1])$ for any $p_1 \in [0,1]$ are orbit arcs, and that there is a small number $\delta \in (0, d_{\mathbb{A}})$ such that $B_{\delta}(\gamma) := \{ y \in S \mid d(y, \gamma ) < \delta \} \subset U$ is an open disk
as in Figure~\ref{fig:arc}.
\begin{figure}
\begin{center}
\includegraphics[scale=0.5]{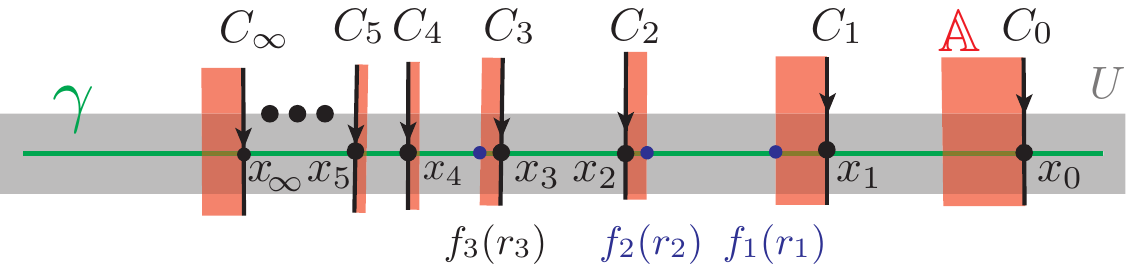}
\end{center}
\caption{An example of a transverse closed arc $\gamma$, its neighborhood $U$, and a subset of an open annulus $\mathbb{A}$}
\label{fig:arc}
\end{figure}
Moreover, there are an orbit arc $C_{\infty}$ with $x_\infty \in C_\infty$ and a sequence $(C_i)_{i \in \Z_{\geq 0}}$ of orbit arcs which intersect $\gamma$ exactly once, are contained in $\omega(x)$, and converge to $C_{\infty}$ such that $\partial C_i \cap U = \emptyset$.
For any $i \in \Z_{\geq 0}$, put $C'_i := C_i \cap U \subset \omega (x)$ and let $x_i \in C'_i \cap \gamma$ be the intersection point.
By Lemma~\ref{lem:transverse_bdry}, any transverse closed arc whose interior is contained in $\mathbb{A}$ and whose boundary is contained in $\partial \mathbb{A}$ and contains a point in $C'_i \subset \omega(x)$, connects $\omega(x)$ and $\partial_1$. 
For any $i \in \Z_{> 0}$, let $f_{i}:[0,l_{i}] \to \gamma \subset S$ be the closed arc parameterized by arc length from $x_i$ to either $x_{i-1}$ or $x_{i+1}$ such that $f_{i}((0,r_i)) \subset \mathbb{A}$ and $f_{i}(r_i) \in \partial_1$ for some real number $r_i \in (0, l_i)$.
By definition, we obtain $f_i(l_i) \in \omega(x)$.
Since $\gamma$ is compact, we have $\lim_{i \to \infty} l_i = 0$. 
Fix a large integer $N$ such that $l_N < \delta < d_{\mathbb{A}}$. 
Then $f_{N}([0, l_{N}]) \subset U$. 
%
Any closed arc $\mu_i \colon [0,\min\{l_i, \delta\}] \to S$ from $x_i$ parameterized by arc length which does not intersect $C'_i$ except the starting point is contained in $\mathbb{A}$ (i.e. $\mu_i((0,\min\{l_i, \delta\}]) \subset \mathbb{A}$) because $\min\{l_i, \delta\} \leq \delta <  d_{\mathbb{A}}$. 
In particular, since $l_N < \delta$, the closed arc $f_{N} : [0, l_{N}] \to S$ is contained in $\mathbb{A}$ except the starting point $x_i$ (i.e. $f_{N}((0, l_{N}]) \subset \A$). 
Therefore $f_N(l_N) \in \A$, which contradicts $f_{N}(l_N) \in \omega(x) \subset \partial \mathbb{\A} = \overline{\A} - \A$.
Thus the assertion holds. 
\end{proof}

We show that each limit quasi-circuit is not a quasi-Q-set, and that each quasi-Q-set is not a limit quasi-circuit. 

\begin{lemma}\label{lem4-07a}
There are no limit quasi-circuits that are also quasi-Q-sets.
\end{lemma}

\begin{proof}
Let $Q$ be a quasi-Q-set. 
Then there is a closed transversal $\gamma$ which intersects $Q$ infinitely many times.
Since the intersection $\gamma \cap Q$ is closed,  the intersection $\gamma \cap Q$ has an accumulation point. 
Lemma~\ref{lem:acc_lc} implies that $Q$ is not a limit quasi-circuit. 
\end{proof}


We have the following equivalence. 

\begin{lemma}\label{lem:equiv_locally_Q}
The following statements are equivalent for a point $x$:
\\
{\rm(1)} $\omega(x) \cap \mathrm{LD}(v) \neq \emptyset$. 
\\
{\rm(2)} The $\omega$-limit set $\omega(x)$ is locally dense.
\\
{\rm(3)} The $\omega$-limit set $\omega(x)$ is a locally dense Q-set which is not transversely Cantor.
\\
In any case, we have that $x \in \mathrm{LD}(v)$ and $\overline{O(x)}= \omega(x)$. 
\end{lemma}

\begin{proof}
Trivially, assertion {\rm(3)} implies assertions {\rm(1)} and {\rm(2)}. 
Suppose that $\omega(x) \cap \mathrm{LD}(v) \neq \emptyset$. 
\cite[Theorem VI]{cherry1937topological} implies that there is a Poisson stable point $y \in \omega(x) \cap \mathrm{LD}(v)$ such that $\emptyset \neq O(x) \cap \mathop{\mathrm{int}} (\overline{O(y)})$ and so that $x \in \overline{O(y)}$. 
Since $y \in \omega(x)$, this means that $\omega(x) = \overline{O(x)} = \omega(y) = \overline{O(y)}$ is a locally dense Q-set. 
By \cite[Lemma~2.3]{yokoyama2016topological}, we have $\omega(y) \cap \mathrm{E}(v) = \overline{O(y)} \cap \mathrm{E}(v) \subseteq \overline{\mathrm{LD}(v)} \cap \mathrm{E}(v) =  \emptyset$.
From Lemma~\ref{lem:excptional_ch}, the $\omega$-limit set $\omega(y) = \omega(x) = \overline{O(x)}$ is not transversely Cantor.

Suppose that $\omega(x)$ is locally dense.
Then $\emptyset \neq O(x) \cap \mathop{\mathrm{int}} (\omega(x)) \subseteq O(x) \cap \mathop{\mathrm{int}} (\overline{O(x)})$. 
Therefore $x \in \mathop{\mathrm{int}} (\omega(x))$ and so $\overline{O(x)}= \omega(x)$. 
This means that $\omega(x)$ is a locally dense Q-set and so that $x \in \omega(x) \cap \mathrm{LD}(v)$. 
%
\end{proof}

We have the following observation. 

\begin{lemma}\label{lem:equiv_exceptional_Q}
If the $\omega$-limit set $\omega(x)$ of a point $x$ intersects $\mathrm{E}(v)$, then $\omega(x)$ is a quasi-Q-set that is not locally dense.
\end{lemma}

\begin{proof}
Suppose $\omega(x)$ contains a non-closed recurrent orbit $O \subset \mathrm{E}(v)$. 
We claim that $\omega(x) \cap \mathrm{LD}(v) = \emptyset$. 
Indeed, assume that $\omega(x) \cap \mathrm{LD}(v) \neq \emptyset$. 
Lemma~\ref{lem:equiv_locally_Q} implies that $x \in \mathrm{LD}(v)$ and that $\overline{O(x)}= \omega(x)$. 
\cite[Lemma~2.3]{yokoyama2016topological} implies that $O \subseteq \omega(x) \cap \mathrm{E} = \overline{O(x)} \cap \mathrm{E}(v) \subseteq \overline{\mathrm{LD}(v)} \cap \mathrm{E}(v) =  \emptyset$, which is a contradiction. 


By Lemma~\ref{lem:equiv_locally_Q}, the $\omega$-limit set $\omega(x)$ is not locally dense. 
Since $O$ is non-closed recurrent, by Lemma~\ref{loops}, taking a small transverse arc, the waterfall construction implies that there is a closed transversal $\gamma$ intersecting $O$ infinitely many times. 
This means that  $\omega(x)$ is a quasi-Q-set that is not locally dense.
\end{proof}

We show the key lemma as follows. 

\begin{proof}[Proof of Lemma~\ref{lemma:a}]
Lemma~\ref{lem:decomp} implies that $S = \mathop{\mathrm{Cl}}(v) \sqcup \mathrm{P}(v) \sqcup \mathrm{R}(v) = \mathop{\mathrm{Sing}}(v) \sqcup \mathop{\mathrm{Per}}(v) \sqcup \mathrm{P}(v) \sqcup \mathrm{LD}(v) \sqcup \mathrm{E}(v)$. 
Lemma~\ref{cor:limit_cycle} and Lemma~\ref{lem4-07a} imply that the five possible invariant subsets in the lemma are exclusive. 
Let $x$ be a point whose orbit is not closed. 
We may assume that $\omega(x)$ is not contained in $\mathop{\mathrm{Sing}}(v)$.

\begin{claim}\label{claim:non-rec_sing_nonproper}
We may assume that $\omega(x)$ contains a non-recurrent orbit $O$ and $\omega(x) \subset \mathop{\mathrm{Sing}}(v) \sqcup \mathrm{P}(v)$. 
\end{claim}
\begin{proof}
From Lemma~\ref{lem:equiv_locally_Q}, we may assume that $\omega(x) \cap \mathrm{LD}(v) = \emptyset$ and that $\omega(x)$ is not locally dense. 
Lemma~\ref{lem:equiv_exceptional_Q} implies that we may assume that $\omega(x) \cap \mathrm{E}(v) = \emptyset$. 
Then $\omega(x) \cap \mathrm{R}(v) = \emptyset$ and so $\omega(x) \subseteq \mathop{\mathrm{Cl}}(v) \sqcup \mathrm{P}(v)$.
If $\omega(x)$ contains a periodic orbit, then $\omega(x)$ is a semi-attracting limit cycle, because of Lemma~\ref{lem:limit_cycle}.
Thus we may assume that $\omega(x)$ contains neither periodic orbits. 
Then $\omega(x) \subseteq \mathop{\mathrm{Sing}}(v) \sqcup \mathrm{P}(v)$. 
Since $\omega(x) \not\subseteq \Sv$, the $\omega$-limit set $\omega(x)$ contains a non-recurrent orbit $O$. 
\end{proof}

\begin{claim}\label{claim:infinite_ct}
We may assume that $O^+(x)$ has no closed transversal intersecting it infinitely many times. 
\end{claim}
\begin{proof}
If there is a closed transversal $\gamma$ that intersects $O^+(x)$ infinitely many times, then Lemma~\ref{lem4-07} implies $\omega(x) = \omega(x')$ is either a quasi-Q-set, an essential limit cycle, or a quasi-circuit for any point $x' \in O^+(x) \cap \gamma$. 
This means that the assertion of Lemma~\ref{lemma:a} holds.
\end{proof}


%
By Claim~\ref{claim:non-rec_sing_nonproper}, take a non-recurrent point $y \in O \subset \omega(x)$ and a transverse closed arc $J: [-1,0] \to S$ with $y = J(0)$ and $x_0 := J(-1) \in O^+(x)$ such that $O^+(x_0)$ intersects $J((-1,0))$ infinitely many times. 
Write $I := J([-1,0])$ and $\mathrm{int} I := J((-1,0))$.
Let $f_v : J' \to I$ be the first return map with the maximal domain $J' \subseteq I$ and $x_n :=(f_v)^n(x_0)$ $n$-th return of $x_0$. 
Denote by $C_n$ the orbit arc from $x_n$ to $x_{n+1}$ and by $I_n \subset I$ the closed arc with $\partial C_n = \partial I_n = \{ x_n, x_{n+1} \}$. 

\begin{claim}\label{claim:ori_pre_loop}
We may assume that $f_v\vert_{J' \cap \mathrm{int} I}$ is orientation-preserving by shortening the transverse intervals $I$. 
\end{claim}
\begin{proof}
Otherwise there are a subsequence $(x_{k_n})_{n \in \Z_{>0}}$ of $(x_{n})_{n \in \Z_{>0}}$ in $J$ converging to $J(0)$ and small neighborhoods $U_{k_n}$ of the unions $I_{k_n} \cup C_{k_n}$ each of which is an open M\"obius band such that $U_{k_n} \cap U_{k_m} = \emptyset$ for any $n \neq  m \in \Z_{>0}$, which contradicts that $S$ has finite non-orientable genus. 
\end{proof}

\begin{claim}
We may assume that 
\[
\{ y \} = I \cap \omega(x_0)  \subset \overline{I \cap O^+(x_0)}
\]
 by shortening $I$. 
\end{claim}
\begin{proof}
Otherwise there is a point  ${x_N} \in J' \cap \mathrm{int} I \cap O^+(x_0)$ with $I_{N} \cap \omega(x_0) \neq \emptyset$ such that $I_{N}$ intersects $O^+(x_0)$ infinitely many times. 
Then the union of $I_{N} \cup C_{N}$ is a loop. 
By Claim~\ref{claim:ori_pre_loop}, from the waterfall construction to the loop $C_{N} \cup I_{N}$, there is a closed transversal $T_{x_N}$ near the loop $C_{N} \cup I_{N}$ such that $T_{x_N}$ intersects $O^+(x_0) \subset O^+(x)$ infinitely many times, which contradicts Claim~\ref{claim:infinite_ct}. 
\end{proof}

Then $\{ J(0) \} = J([-1,0]) \cap \omega(J(-1)) \subset \overline{J([-1,0]) \cap O^+(J(-1))}$. 
By Lemma~\ref{lem-quasi-limit-cricuits}, the $\omega$-limit set $\omega(x)$ is a limit quasi-circuit such that $x$ is wandering and is not contained in the orbit closures of points outside of $O(x)$. 
\end{proof}

\subsection{Classification of quasi-circuits}

We have the following dichotomy. 

\begin{lemma}\label{lem:qc_ch}
A quasi-semi-attracting limit quasi-circuit either is the image of a circle or is not locally connected exclusively. 
\end{lemma}

\begin{proof}
%
Let $\gamma$ be a quasi-semi-attracting limit quasi-circuit. 
By definition of non-trivial quasi-circuit, the quasi-circuit $\gamma$ consists of singular points and non-recurrent points, and there is an open annulus $\A$ such that $\gamma$ is a boundary component of $\A$. 

\begin{claim}
$\gamma$ can be realized as a boundary component of an open disk. 
\end{claim}
\begin{proof}
Taking $\A$ small and using a bump function whose inverse image of the zero is another boundary component $\mu$ of $\A$, cutting $\mu$ into two boundary components, and collapsing the two boundary components into two singletons, the resulting space $S'$ is a surface, the resulting subset from $\A$ is an open disk $U \subset S'$, and the quasi-circuit $\gamma$ can be realized as the boundary component of the open disk $U$. 
\end{proof}

Thus we may assume that there is an open disk $U$ such that $\gamma$ is a boundary component of $U$. 
By Riemann mapping theorem, there is a Riemann mapping $\varphi : U \to D$ from an open disk $U$ to a unit open disk $D$ in a complex plane. 
By Carath\'eodory-Torhorst theorem (see \cite{torhorst1921rand} and also  \cite[Theorem~16.6]{milnor2011dynamics} for detail), the boundary $\partial U = \gamma$ is locally connected if and only if the inverse mapping $\varphi^{-1} : D \to U$ extends continuously to a map, also denoted by  $\varphi^{-1}$, from the closed disk $\overline{D}$ onto $\overline{U}$. 

We claim that $\gamma$ is locally connected if and only if $\gamma$ is the image of a circle. 
Indeed, if $\gamma$ is locally connected, then $\gamma = \varphi^{-1}(\partial D)$ is the image of a circle. 
Conversely, suppose that $\gamma$ is the image of a circle. 
Let $p \colon \mathbb{S}^1 \to \gamma$ be the continuous surjection from a circle $\mathbb{S}^1$. 
Since the circle $\mathbb{S}^1$ is compact and $\gamma$ is Hausdorff, the continuous surjection $p$ is closed and so is a quotient mapping. 
Because each quotient space of a locally connected space is locally connected (cf. \cite[Proposition~12 (p.112)]{bourbaki2013general}), the image $p(\mathbb{S}^1) = \gamma$ is locally connected. 
\end{proof}

\subsection{Classification of quasi-Q-sets}

In this subsection, we classify quasi-Q-sets. 

\subsubsection{Inclusion relation between Q-sets and quasi-Q-sets}

We describe a property between non-closed recurrent orbits and closed transversals. 

\begin{lemma}\label{lem3-02a}
For any non-closed recurrent point $x \in S$, there is a closed transversal $\gamma$  through $O(x)$ such that the intersection $\gamma \cap O(x)$ is infinite.
Moreover, each closed transversal through a non-closed recurrent point $y$ is essential and intersects $O(y)$ infinitely many times.
\end{lemma}

\begin{proof}
Fix a point $x \in \mathrm{R}(v)$ and a transverse arc $I \subset U$ such that $x$ is the interior point of $I$. 
Then $\vert  I \cap O(x)\vert  = \infty$. 
By Lemma~\ref{loops}, there are an orbit arc $C$ in $O(x)$ and a transverse closed arc $J \subseteq I$ such that the union $\mu := J \cup C$ is a loop with $C \cap J = \partial C = \partial J$ and that the return map along $C$ is orientation-preserving between neighborhoods of $\partial C$ in $I$. 
By the waterfall construction to the loop $\mu$, there is a closed transversal $\gamma$ intersecting $O(x)$ near $\mu$.
Since $x$ is non-closed recurrent, the intersection $\gamma \cap O(x)$ is infinite.

Let $\nu$ be a closed transversal through a non-closed recurrent point $y$. 
By time reversion if necessary, we may assume that $y \in \omega(y)$. 
From the recurrence of $y$, we have that $y \in \overline{\nu \cap O^+(y)}$ (i.e. the point $y$ is an accumulation point of $\nu \cap O^+(y)$). 
Assume that $\nu$ is inessential.
Let $S^*$ be the resulting closed surface from the compact surface $S$ by collapsing all boundary components into singletons, and $v^*$ the resulting flow on $S^*$ from $v$.
Then $\nu$ is null homotopic in $S^*$ and so $\nu$ is the boundary of an open disk $D$ with either $D \cap O_{v^*}(y) = O_{v^*}^+(y)$ or $D \cap O_{v^*}(y) = O_{v^*}^-(y)$. 
Therefore $\nu \cap O_{v^*}(y) = \{ y \}$ and so $y$ is not recurrent with respect to $v^*$. 
By construction of $v^*$, the point $y$ is also not recurrent with respect to $v$, which contradicts the recurrence of $y$.
\end{proof}

This implies the following corollary.

\begin{corollary}\label{preQ}
A Q-set is a quasi-Q-set. 
\end{corollary}

\begin{proof}
Let $\gamma$ be a Q-set.
Then there is a non-closed recurrent orbit $O \subset \gamma$.
By Lemma~\ref{lem3-02a}, there is an essential closed transversal intersecting $O$ infinitely many times.
\end{proof}

\subsubsection{Non-recurrent orbits in the $\omega$-limit sets of points}

We recall the following Ma\v{i}er's result \cite{mayer1943trajectories} (cf.  \cite[Theorem 2.4.4 p.32]{nikolaev1999flows},  \cite[Theorem 4.2]{aranson1996maier}, \cite[Lemma 3.7]{yokoyama2017decompositions}). 

\begin{lemma}[Ma\v{i}er]\label{pos_rec}
Let $v$ be a flow on a compact surface $S$. 
A point $x \in \omega(z)$ for some point $z \in S$ with $\omega(x) \setminus  \mathop{\mathrm{Cl}}(v) \neq \emptyset$ is non-closed positively recurrent $(\mathrm{i.e. } \,\, x \in \omega(x) \cap \mathrm{R}(v) \, )$. 
\end{lemma}

The similar argument of the proof of the previous result implies the following result.

\begin{lemma}\label{neg_rec}
Let $v$ be a flow on a compact surface $S$. 
A point $x \in \omega(z)$ for some point $z \in S$ with $\alpha(x) \setminus \mathop{\mathrm{Cl}}(v) \neq \emptyset$ is non-closed negatively recurrent. 
In particular, we have $x \in \alpha(x) \cap \mathrm{R}(v)$. 
\end{lemma}

\begin{proof}
Since $\alpha(x) \setminus  \mathop{\mathrm{Cl}}(v) \neq \emptyset$, the point $x$ is not closed. 

\begin{claim}
If $O(x) = O(z)$, then $x$ is non-closed negatively recurrent. 
\end{claim}

\begin{proof}
Suppose that $O(x) = O(z)$. 
Since $x$ is not closed, so is $z$. 
By $x \in \omega(z) = \omega(x)$, the point is positively recurrent. 
\cite[Theorem~VI]{cherry1937topological} implies that there is a Poinsson stable point $z' \in S$ with 
$x \in \omega(z) = \overline{O(x)} = \overline{O(z')} = \alpha(z')$. 
Applying the dual statement of Lemma~\ref{pos_rec} to $x \in \alpha(z')$ with $\alpha(x) \setminus \mathop{\mathrm{Cl}}(v) \neq \emptyset$, the point $x$ is non-closed negatively recurrent. 
\end{proof}

Thus we may assume that $O(x) \neq O(z)$. 
Fix a non-closed point $y \in \alpha(x) \setminus \mathop{\mathrm{Cl}}(v)$. 
Then there is a transverse closed arc $I_{[-1,1]}: [-1,1] \to S$ with $y = I_{[-1,1]}(0)$ such that the negative orbit $O^-(x)$ intersects $I_{[-1,1]}([-1,0])$ infinitely many times. 
Denote by $I := I_{[-1,1]}([-1,0])$ a directed closed interval. 
Therefore there is a sequence $(x_i)_{i \in \Z_{\geq 0}}$ of points in $O^-(x) \cap I$ with $x_{i+1} \in O^-(x_i)$ which converges to $y$ monotonically from one side. 
Denote by $I_{a,b}$ the sub-arc in $I$ whose boundary consists of $a$ and $b$ for any points $a, b \in I$ and by $C_{a,b}$ the orbit arc in an orbit $O$ from $a$ to $b$ for any points $a, b \in O \cap I$. 

Assume that $x$ is not negatively recurrent (i.e. $x \notin \alpha(x)$). 
Then there is an open sub-arc $J$ in $I$ with $\{ x_2 \} = J \cap O^-(x)$. 
By $x_2 \in \omega(z)$, the first return map $f_{v,J}$ on $J$ induced by $v$ is well-defined and injective. 
From the finiteness of genus of $S$, by replacing $x$ with a point of $O^-(x)$, 
we may assume that the restriction of the first return map $f_{v,I}$ to the transverse closed arc $I$ induced by $v$, restricted to a \nbd of $f_{v,I}^{-1}(O^-(x) \cap I)$ in $I$ is orientation-preserving. 
Therefore $I$ and $O^-(x)$ intersect in a same orientation infinitely many times. 
\begin{claim}
We can define a strictly increasing subsequence $(n_i)_{i \in \Z_{\geq 0}}$ of $\Z_{\geq 0}$ with $n_i +3 \leq n_{i+1}$ and a sequence $(z_i)_{i \in \Z_{\geq 0}}$ of $J \cap O^+(z)$ with $z_{i+1} \in O^+(z_i)$ converging to $x_2$ monotonically from one side in $J$ such that $C_{z_{i-1},z_{i}} \cap \mathop{\mathrm{int}} I_{z_{i-1}, z_{i}} \neq \emptyset$, $C_{z'_{i}, z_{i}} \cap \mathop{\mathrm{int}} I_{x_{n_{i}},x_{n_{i}+1}} \neq \emptyset$, and $C_{z_0, z'_{i}} \cap I_{x_{n_{i}}, y} = \emptyset$ for any $i \in \Z_{>0}$, where $z'_{i} \in \mathop{\mathrm{int}} C_{z_{i-1},z_{i}} \cap I_{z_{i-1},z_{i}}$ is the first return image of $z_{i}$ into $I_{z_{i-1},z_i}$ induced by the time reversed flow of $v$. 
\end{claim}
\begin{proof}
By induction, 
fix a point $z_0 \in J \cap O^+(z)$ and $n_0 = 0$ such that $O^+(z)$ intersects $I_{z_0, x_2}$ infinitely many times. 
Since the sequence $(x_k)_{k \in \Z_{\geq 0}}$ converges to $y$ monotonically from one side, by $O(x) \neq O(z)$, for any $i \in \Z_{\geq 0}$, there are an integer $k_i \geq 3$ 
and a point $z_{i+1} \in I_{z_i, x_2} \cap O^+(z_i)$ with $C_{z_{i},z_{i+1}} \cap \mathop{\mathrm{int}} I_{z_{i}, z_{i+1}} \neq \emptyset$, $C_{z_0,z'_{i+1}} \cap \mathop{\mathrm{int}} I_{x_{n_{i}+k_i}, y} = \emptyset$, and $C_{z'_{i+1}, z_{i+1}} \cap I_{x_{n_{i}+k_i}, y} \neq \emptyset$. 
Fix an integer $n_{i+1} \geq n_i + k_i \geq n_i +3$ such that $C_{z'_{i+1}, z_{i+1}} \cap \mathop{\mathrm{int}} I_{x_{n_{i+1}}, x_{n_{i+1}+1}} \neq \emptyset$ and $C_{z_0, z'_{i+1}} \cap I_{x_{n_{i+1}}, y} = \emptyset$. 
Then $C_{z_{i-1},z_{i}} \cap \mathop{\mathrm{int}} I_{z_{i-1}, z_{i}} \neq \emptyset$, $C_{z'_{i}, z_{i}} \cap \mathop{\mathrm{int}} I_{x_{n_{i}}, x_{n_{i}+1}} \neq \emptyset$, and $C_{z_0, z'_{i}} \cap I_{x_{n_{i}}, y} = \emptyset$ for any $i \in \Z_{>0}$. 
\end{proof}

Fix a Riemannian metric $g$ on $S$ which induces the Riemannian distance $d_g$. 
Since the sequence $(z_i)_{i \in \Z_{\geq 0}}$ of $J \cap O^+(z)$ converging to $x_2$ monotonically from one side, the sequence of the lengths of $I_{z_{i+1},z_i}$ converges to zero. 
For any $i \in \Z_{\geq 0}$, let $f_{v,I_{z_i, x_2}}$ be the first return map from $I_{z_i, x_2}$ to $I_{z_{i+1}, x_2} \subset I_{z_i, x_2}$ induced by $v$. 
Then $C_{z'_{i+1}, z_{i+1}} \cap I_{z_{i}, x_2} = \{ z'_{i+1}, z_{i+1} \}$ and $C_{z'_{i+1}, z_{i+1}} \cap \mathop{\mathrm{int}} I_{x_{n_{i+1}},x_{n_{i+1}+1}} \neq \emptyset$. 
Since $\mathop{\mathrm{int}} C_{z_i,z_{i+1}} \cap I_{z_i, z_{i+1}} \neq \emptyset$, we have that $z_i \neq z'_{i+1}$ and so that the closed intervals $I_{z'_{i+1}, z_{i+1}} \subset J$ are pairwise disjoint. 
Therefore the unions $\gamma_i := C_{z'_{i+1}, z_{i+1}} \cup I_{z'_{i+1}, z_{i+1}} \subset O^+(z) \cap J$ are pairwise disjoint loops intersecting $\mathop{\mathrm{int}} I_{x_{n_{i+1}},x_{n_{i+1}+1}}$. 
Let $\A_i$ be the connected component of $S - \bigcup_{k \in \Z_{\geq 0}} \gamma_k$ intersecting $I_{z'_{i+1},z_i}$. 

\begin{claim}
We may assume that $\A_i$ is a closed annulus whose boundary is a disjoint union $\gamma_i \sqcup \gamma_{i+1} \subset O^+(z_0) \cup J$ such that the pairwise disjoint loops $\gamma_i$ are homotopic to each other.
\end{claim}
\begin{proof}
Then the boundary of any domain $\A_i$ is contained in $(O^+(z_0) \cup J) \sqcup \partial S$. 
Since there are at most finitely many boundary components and finite genus, by renumbering, we may assume that each domain $\A_i$ is annular and that the restriction of $f_{v,I_{z_i, x_2}}$ whose domain is a small neighborhood of $z'_{i+1} \in I_{z_i, z_{i+1}}$ and codomain is a small neighborhood of $z_{i+1}$ is orientation-preserving. 
Then $\A_i$ is a closed annulus whose boundary is a disjoint union $\gamma_i \sqcup \gamma_{i+1} \subset O^+(z_0) \cup J$. 
Since $S$ is compact, by renumbering, we may assume that the pairwise disjoint loops $\gamma_i$ are homotopic to each other.
\end{proof}

Then the union $\A_{i-1} \cup \A_{i}$ is also a closed annulus with $\A_{i-1} \cap \A_{i} = \gamma_i$. 
Denote by $d_0 > 0$ the distance between $\gamma_0$ and $\gamma_1$ in $\A_0$ (i.e. $d_0 := d_g(\gamma_0,\gamma_1)$, where $d_g(A,B) := \min_{a \in A, b \in B}d_g(q_0,q_1)$).

Fix a large integer $N \in \Z_{>2}$ such that the length of $I_{x_{n_{i-1}},x_{n_{i+1}}}$ is less than $d_0/2$ for any $i \geq N$. 
Then $x_2 \notin O^-(x_{n_N +1})$. 
Since $\gamma_i \cap O^-(x_2) \subset (O^+(z) \cup J) \cap O^-(x_2) = \emptyset$ for any $i \geq N$, put $D := \min \{ d_g(x_{n_N +1}, \partial \A_N), d_g(x_{n_{N+1} +1}, \partial \A_{N+1}) \} = \min \{ d_g(\{ x_{n_N +1}\} , \gamma_N \sqcup \gamma_{N+1}), d_g(\{ x_{n_{N+1} +1}\} , \gamma_{N+1} \sqcup \gamma_{N+2}) \} >0$. 

For any $i \in \Z_{\geq 0}$, applying the waterfall construction to the loop $\gamma_i$, there is a closed transversal $T_i$ isotopic to $\gamma_i$ with $x_{n_{i}},x_{n_{i}+1} \notin T_i$ such that $T_i$ intersects $\mathop{\mathrm{int}} I_{x_{n_{i}},x_{n_{i}+1}}$ transversely and $d_{H}(T_i, \gamma_i) < \min \{D, d_0, d_g(\gamma_{i-1}, \gamma_{i}), d_g(\gamma_{i}, \gamma_{i+1}) \}/4$, where $d_H$ is the Hausdorff distance. 
For any $i \in \Z_{\geq 0}$, denote by $\A'_i$ the closed annulus whose boundary is $T_i \sqcup T_{i+1}$ and which is near $\A_i$. 
Then the union $\A'_{i} \cup \A'_{i+1}$ is also a closed annulus with $\A'_{i} \cap \A'_{i+1} = T_{i+1}$ and $\partial (\A'_{i} \cup \A'_{i+1}) = T_{i} \sqcup T_{i+2}$. 

\begin{claim}\label{claim:ineq_001}
$d_g(x_{n_N +1}, \partial \A'_N) \geq 3D/4$. 
\end{claim}
\begin{proof}
We have the following inequality: 
\[
\begin{split}
& \,\, d_g(x_{n_N +1}, \partial \A'_N) =  d_g(x_{n_N +1}, T_N \sqcup T_{N+1}) 
\\
= & \, \min \{ d_g(x_{n_N +1}, T_N), d_g(x_{n_N +1}, T_{N+1})\} \\
\geq & \,  \min \{ d_g(x_{n_N +1}, \gamma_N) - d_H(T_N,\gamma_{N}), d_g(x_{n_N +1}, \gamma_{N+2}) - d_H(T_{N+2},\gamma_{N+2})\} 
\\
\geq & \,  D - D/4 = 3D/4
\end{split}
\]
\end{proof}

\begin{claim}
The closed transversal $T_i$ intersects $I_{x_{n_{i}},x_{n_{i}+1}}$ exactly once for any $i \in \Z_{\geq N}$. 
\end{claim}
\begin{proof}
Assume that $T_i$ intersects $I_{x_{n_{i}},x_{n_{i}+1}}$ at least twice. 
Since $\A'_i$ is a closed annulus with $\partial \A'_i = T_i \sqcup T_{i+1}$ such that $T_i$ and $T_{i+1}$ are closed transversals, the transverse closed arc $I_{x_{n_{i}},x_{n_{i}+1}}$ goes outside of $\A'_i$ and goes into $\A'_i$ from $\gamma_i$ with respect to the positive or negative direction. 
The fact that the union $\bigcup_{k=0}^i \A'_k$ is a closed annulus whose boundary components are closed transversals implies that  $I_{x_{n_{i}},x_{n_{i}+1}} \cap T_k \neq \emptyset$ for any $k =0, 1, \ldots , i$. 
Since the transverse closed arc $I_{x_{n_{i}},x_{n_{i}+1}}$ goes through $\A'_0$, it contains a sub-arc in $\A'_0$ whose boundary component consists of a point in $T_0$ and a point in $T_1$. 
Then the length of $I_{x_{n_{i}},x_{n_{i}+1}}$ is more than $d_0/2$, which contradicts that the length is less than $d_0/2$. 
\end{proof}

By the previous claim, we have that $x_{n_{N}}< T_N \cap I_{x_{n_{N}},x_{n_{N}+1}} < x_{n_{N}+1} < x_{n_{N+1}} < T_{N+1} \cap I_{x_{n_{N+1}},x_{n_{N+1}+1}} < x_{n_{N+1}+1} < x_{n_{N+2}}$ in the closed interval $I$. 

\begin{claim}
$x_{n_N +1} \in \A_N$. 
\end{claim}
\begin{proof}
Since $\A'_{N}$ is a closed annulus with $\partial \A'_{N} = T_{N} \sqcup T_{N+1}$, by $T_N \cap I_{x_{n_{N}},x_{n_{N}+1}} < x_{n_{N+1}} < T_{N+1} \cap I_{x_{n_{N+1}},x_{n_{N+1}+1}}$, we obtain that  $x_{n_N +1} \in \A'_N$.
By Claim~\ref{claim:ineq_001}, we have that $d_g(x_{n_N +1}, \partial \A'_N) \geq 3D/4 > D/4 > \max \{ d_{H}(T_N, \gamma_N), d_{H}(T_{N+1}, \gamma_{N+1})\} \geq d_H(\partial \A'_N, \partial \A_N) \geq d_H(\A'_N, \A_N)$. 
Since the boundary $\partial \A'_N = T_N \sqcup T_{N+1}$ is isotopic to $\partial \A_N = \gamma_N \sqcup \gamma_{N+1}$,  the annulus $\A'_N$ is isotopic to the annulus $\A_N$ with $d_H(\partial \A'_N, \partial \A_N)< D/4 < d_g(x_{n_N +1}, \partial \A'_N)$ and so $x_{n_N +1} \in \A_N$.
\end{proof}
By the same argument of the proof of the previous claim, we have $x_{n_{N+1} +1} \in \A_{N+1}$. 
From $O^-(x_{n_N +1}) \subseteq O^-(x)$ and $x_{n_{N+1} +1} \in O^-(x_{n_N +1}) \setminus \A_N$, the negative orbit $O^-(x_{n_N +1})$ intersects $\A_N$ but is not contained in $\A_N$. 
By $O^-(x_{n_N +1}) \cap (\bigcup_{k} C_{z'_{k+1}, z_{k+1}}) \subseteq O(x) \cap O(z) = \emptyset$ and $\partial \A_N = \gamma_N \sqcup \gamma_{N+1} \subset O^+(z) \cup (I_{z'_{N}, z_{N}} \sqcup I_{z'_{N+1}, z_{N+1}})$, we have $\emptyset \neq O^-(x_{n_N +1}) \cap \partial \A_N = O^-(x_{n_N +1}) \cap (\gamma_N \sqcup \gamma_{N+1})  = O^-(x_{n_N +1}) \cap (I_{z'_{N}, z_{N}} \sqcup I_{z'_{N+1}, z_{N+1}}) \subset O^-(x) \cap J = \{ x_2 \}$, which contradicts $x_2 \notin O^-(x_{n_N +1})$. 
Thus, the point $x$ is negatively recurrent. 
\end{proof}

\begin{proposition}\label{preQ03}
Let $v$ be a flow on a compact surface $S$. 
An orbit in the $\omega$-limit set of a point in $S$ is non-recurrent if and only if it is a connecting quasi-separatrix. 
\end{proposition}

\begin{proof}
Let $x \in S$ be a point with $x \in \omega(z)$ for some point $z \in S$.
If $x$ is closed, then $O(x)$ is recurrent and is not a connecting quasi-separatrix.
Thus we may assume that $x$ is non-closed. 

\begin{claim}
$(\omega(x) \cup \alpha(x)) \cap \mathop{\mathrm{Per}}(v) = \emptyset$. 
\end{claim}
\begin{proof}
Assume that $\omega(x) \cap \mathop{\mathrm{Per}}(v) \neq \emptyset$. 
By $\omega(x) \subseteq \omega(z)$, we have $\emptyset \neq \omega(x) \cap \mathop{\mathrm{Per}}(v) \subseteq \omega(z) \cap \mathop{\mathrm{Per}}(v)$. 
Lemma~\ref{lemma:a} implies that $\omega(z)$ is a limit cycle and so $x \in \omega(z)  \subseteq \mathop{\mathrm{Per}}(v)$, which contradicts the non-closedness of $x$. 
Thus $\omega(x) \cap \mathop{\mathrm{Per}}(v) = \emptyset$. 

Assume that $\alpha(x) \cap \mathop{\mathrm{Per}}(v) \neq \emptyset$. 
By $\alpha(x) \subseteq \omega(z)$, we obtain $\emptyset \neq \alpha(x) \cap \mathop{\mathrm{Per}}(v) \subseteq \omega(z) \cap \mathop{\mathrm{Per}}(v)$. 
Lemma~\ref{lemma:a} implies that $\omega(z)$ is a limit cycle and so $x \in \omega(z)  \subseteq \mathop{\mathrm{Per}}(v)$, which contradicts the non-closedness of $x$. 
Thus $\alpha(x) \cap \mathop{\mathrm{Per}}(v) = \emptyset$. 
\end{proof}

Suppose that $O(x)$ is a connecting quasi-separatrix. 
Then $x$ is not-closed recurrent. 
Conversely, suppose that $x$ is not-closed non-recurrent. 
If $\omega(x) \setminus \mathop{\mathrm{Cl}}(v) \neq \emptyset$, then Lemma~\ref{pos_rec} implies that $x$ is positive recurrent, which contradicts the non-recurrence of $x$. 
If $\alpha(x) \setminus \mathop{\mathrm{Cl}}(v) \neq \emptyset$, then Lemma~\ref{neg_rec} implies that $x$ is negative recurrent, which contradicts the non-recurrence of $x$. 
Thus $\omega(x) \cup \alpha(x) \subset \mathop{\mathrm{Cl}}(v)$. 
By $(\omega(x) \cup \alpha(x)) \cap \mathop{\mathrm{Per}}(v) = \emptyset$, we obtain $\omega(x) \cup \alpha(x) \subset \mathop{\mathrm{Sing}}(v)$. 
This means that $O(x)$ is a connecting quasi-separatrix. 
\end{proof}

Proposition~\ref{preQ03} implies the following reduction. 

\begin{corollary}\label{cor:td}
Let $v$ be a flow with totally disconnected singular points on a compact surface $S$. 
Then each orbit in the $\omega$-limit set of a point is non-recurrent if and only if it is a connecting separatrix. 
\end{corollary}

\subsubsection{Properties of the resulting flows by collapsing connected components of singular points into singletons}

Consider a flow $v$ on a surface $S$. 
Denote by $\bm{S_{\mathrm{me}}}$ the metric completion of the difference $S - \mathop{\mathrm{Sing}}(v)$ and by $\bm{v_{\mathrm{me}}}$ the resulting $\R$-action such that the new points are singular points.
Let $p_{\mathrm{me}}: S_{\mathrm{me}} \to S$ be the canonical projection. 
Then $\mathop{\mathrm{Sing}}(v_{\mathrm{me}}) = p_{\mathrm{me}}^{-1}(\mathop{\mathrm{Sing}}(v))$.
Let $\bm{S_{\mathrm{col}}}$ be the resulting space from $S_{\mathrm{me}}$ by collapsing any connected components of $\mathop{\mathrm{Sing}}(v_{\mathrm{me}})$ into singletons. 
By construction, the resulting space $S_{\mathrm{col}}$ is a disjoint union of closed surfaces. 
Let $\bm{v_{\mathrm{col}}$} the resulting $\R$-action on $S_{\mathrm{col}}$, and $p_{\mathrm{col}}: S_{\mathrm{me}} \to S_{\mathrm{col}}$ the canonical projection (see Figure~\ref{Fig:blowdown}). 
\begin{figure}
\[
\xymatrix@=18pt{
S \ar@{}[d]| {\bigcup} & & S_{\mathrm{me}}\ar[ll]_{p_{\mathrm{me}}} \ar[rr]^{p_{\mathrm{col}}} \ar@{}[d]| {\bigcup} & &  S_{\mathrm{col}} \ar@{}[d]| {\bigcup} \\
S - \mathop{\mathrm{Sing}}(v) & &  S_{\mathrm{me}} - \mathop{\mathrm{Sing}}(v_{\mathrm{me}}) \ar@{=}[ll]_{p_{\mathrm{me}}| }
 \ar@{=}[rr]^{p_{\mathrm{col}}| } & &  S_{\mathrm{col}} - \mathop{\mathrm{Sing}}(v_{\mathrm{col}})
 }
\]
\caption{Canonical quotient mappings induced by the metric completion and the collapse}
\label{Fig:blowdown}
\end{figure}
Then $\mathop{\mathrm{Sing}}(v_{\mathrm{col}}) = p_{\mathrm{col}}(\mathop{\mathrm{Sing}}(v_{\mathrm{me}})) = p_{\mathrm{col}}(p_{\mathrm{me}}^{-1}(\mathop{\mathrm{Sing}}(v)))$.
By construction, we obtain $S - \mathop{\mathrm{Sing}}(v) = S_{\mathrm{me}} - \mathop{\mathrm{Sing}}(v_{\mathrm{me}}) = S_{\mathrm{col}} - \mathop{\mathrm{Sing}}(v_{\mathrm{col}})$. 
We have the following continuity. 

\begin{lemma}\label{lem:collapse_equiv}
Let $v : \R \times S \to S$ be a flow on a surface $S$.
The resulting $\R$-actions $v_{\mathrm{me}}$ and $v_{\mathrm{col}}$ are flows {\rm(i.e.} continuous $\R$-actions{\rm)}. 
\end{lemma}

To demonstrate the previous lemma, we state the following observation. 

\begin{lemma}\label{lem:collapse}
Let $(a_n)_{n \in \Z_{\geq 0}}$ be a sequence of points in $S$ and $(t_n)_{n \in \Z_{\geq 0}}$ be a sequence of points in $\R$ such that the sequence $((-t_n,b_n))_{n \in \Z_{\geq 0}}$ in $\R \times S$ converges a point $(-t_\infty, b_\infty) \in \R \times S$ with $b_\infty \in \Sv$, where $b_n := v(t_n,a_n)$. 
Then $\lim_{n \to \infty} a_n = b_\infty = \lim_{n \to \infty} b_n$. 
\end{lemma}

\begin{proof}
By definition, we have that $a_n = v(-t_n, b_n)$ for any nonnegative integer $n \in \Z_{\geq 0}$, and that $b_\infty = \lim_{n \to \infty} b_n$. 
Since the sequence $(t_n)_{n \in \Z_{\geq 0}}$ in $\R$ converges a point $t_\infty \in \R$, there is a number $T_0 >0$ with $\{ t_n \mid n \in \Z_{\geq 0} \} \subset [-T_0,T_0]$. 
For any positive integer $n \in \Z_{>0}$, there is a \nbd $B_n$ of $b_\infty$ such that $\sup \{ d(b_\infty, v(t, y)) \mid (t,y) \in [-T_0,T_0] \times B_n \} < 1/n$. 
Taking a subsequence of $((-t_n,b_n))_{n \in \Z_{\geq 0}}$, we may assume that $v(-t_n,b_n) \in B_n$. 
Then $d(b_\infty, a_n) = d(b_\infty, v(-t_n,b_n)) < 1/n$ for any positive integer $n \in \Z_{>0}$. 
This means that $b_\infty = \lim_{n \to \infty} a_n$. 
\end{proof}

Using the previous observation, we show the following continuity of specific flows. 

\begin{lemma}\label{lem:continuity_extension}
Let $v$ be a $\R$-action on a surface $S$ and $U$ an open subset of $S$ to which the restriction of $v$ is a continuous action such that the set difference $S - U$ consists of singular points. 
Then $v$ is continuous {\rm (i.e.} a flow{\rm)}. 
\end{lemma}

\begin{proof}
Fix any closed subset $A \subseteq S$. 
By the closedness of $A$ and the openness of $U$, from $A = (A \cap U) \sqcup (A \setminus U)$, the set difference $A \setminus U$ is closed and $A \cap U = \overline{A \cap U} \cap U$.
Since the restriction $v \vert_{\R \times U} \colon \R \times U \to U$ is continuous, the inverse image $v \vert_{\R \times U}^{-1}(A \cap U) = v^{-1}(A \cap U)$ is closed with respect to $\R \times U$. 
Then $\overline{v^{-1}(A \cap U)} \cap (\R \times U) = v^{-1}(A \cap U)$. 
From the closedness of $\mathop{\mathrm{Sing}}(v)$ and definition of product topology, the inverse image $v^{-1}(\mathop{\mathrm{Sing}}(v)) = \R \times \mathop{\mathrm{Sing}}(v)$ is closed. 
%
By $S - U \subseteq \mathop{\mathrm{Sing}}(v)$, we have that $\partial^+ ({v^{-1}(A \cap U)})  = \overline{v^{-1}(A \cap U)} - v^{-1}(A \cap U) \subseteq \R \times (S - U) \subseteq \R \times \mathop{\mathrm{Sing}}(v) = v^{-1}(\mathop{\mathrm{Sing}}(v))$. 

We claim that $\partial^+ ({v^{-1}(A \cap U)}) \subseteq \R \times (A \setminus U)$. 
Indeed, assume $\partial^+ ({v^{-1}(A \cap U)}) \not\subseteq \R \times (A \setminus U)$. 
Fix a point $(-t',p') \in \partial^+ ({v^{-1}(A \cap U)}) \setminus (\R \times (A \setminus U))  \subseteq (\R \times (S - U) ) \setminus (\R \times (A \setminus U))  \subseteq \R \times ((S - U) \setminus A) \subseteq \R \times (\mathop{\mathrm{Sing}}(v) \setminus A)$. 
Then $p' \in \mathop{\mathrm{Sing}}(v) \setminus A$. 
Since $(-t',p') \in \partial^+ ({v^{-1}(A \cap U)})$, there is a sequence $(-t_n, p_n)_{n \in \Z_{\geq 0}}$ of ${v^{-1}(A \cap U)} \subseteq \R \times (A \cap U)$ converging to $(-t',p') \in \R \times \mathop{\mathrm{Sing}}(v)$. 
Lemma~\ref{lem:collapse} implies that $p' = \lim_{n \to \infty} v(-t_n, p_n)$. 
By $(-t_n, p_n) \in {v^{-1}(A \cap U)}$, we have that $v(-t_n, p_n) \in A \cap U$. 
The closedness of $A$ implies that $p' = \lim_{n \to \infty} v(-t_n, p_n) \in A$, which contradicts $p' \in S \setminus A$. 

From the openness of $U$ and definition of product topology, by $A \setminus U \subset \mathop{\mathrm{Sing}}(v)$, the inverse image $v^{-1}(A \setminus U) = \R \times (A \setminus U)$ is closed. 
By the previous claim, the inverse image $v^{-1}(A) = v^{-1}((A \setminus U) \sqcup (A \cap U)) = v^{-1}(A \setminus U) \sqcup v^{-1}(A \cap U) = (\R \times (A \setminus U)) \sqcup v^{-1}(A \cap U) = (\R \times (A \setminus U)) \cup \overline{v^{-1}(A \cap U)}$ is closed. 
%
\end{proof}

\begin{proof}[Proof of Lemma~\ref{lem:collapse_equiv}]
By construction, we obtain $S - \mathop{\mathrm{Sing}}(v) = S_{\mathrm{me}} - \mathop{\mathrm{Sing}}(v_{\mathrm{me}}) = S_{\mathrm{col}} - \mathop{\mathrm{Sing}}(v_{\mathrm{col}})$. 
Therefore $v = v_{\mathrm{me}} = v_{\mathrm{col}}$ on $S - \mathop{\mathrm{Sing}}(v) = S_{\mathrm{me}} - \mathop{\mathrm{Sing}}(v_{\mathrm{me}}) = S_{\mathrm{col}} - \mathop{\mathrm{Sing}}(v_{\mathrm{col}})$. 
Hence $v_{\mathrm{me}}$ and $v_{\mathrm{col}}$ are $\R$-actions. 
Since $S_{\mathrm{col}} - \mathop{\mathrm{Sing}}(v_{\mathrm{col}}) =S_{\mathrm{me}} - \mathop{\mathrm{Sing}}(v_{\mathrm{me}}) = S - \mathop{\mathrm{Sing}}(v)$ is an open subset, Lemma~\ref{lem:continuity_extension} implies that the $\R$-actions $v_{\mathrm{me}} \colon \R \times S_{\mathrm{me}} \to S_{\mathrm{me}}$ and $v_{\mathrm{col}} \colon \R \times S_{\mathrm{col}} \to S_{\mathrm{col}}$ are continuous. 
\end{proof}

\subsubsection{Charcterization of Q-sets}

We show the following statements.

\begin{lemma}\label{pre_q-}
For any positively recurrent point $y \in \mathrm{E}(v)$ and any point $x$ with $y \in \omega(x)$, the $\omega$-limit set $\omega(x) = \omega(y) = \overline{O(y)}$ is a transversely Cantor Q-set. 
\end{lemma}

The following proof is an analogous argument of the proof of \cite[Theorem~2.4.1]{nikolaev1999flows}. 

\begin{proof}[Proof of Lemma~\ref{pre_q-}]
By Lemma~\ref{lem:excptional_ch}, the Q-set $\overline{O(y)}$ is a transversely Cantor Q-set. 
%
Therefore, it suffices to show $\omega(x) = \overline{O(y)}$. 

Assume that there is a point $z \in \omega(x) - \overline{O(y)}$.
Then $O(x) \cap O(y) = \emptyset$. 
Since $S$ is a $T_3$-space, there is a closed disk $B$ whose interior contains $z$ such that $B \cap \overline{O(y)} = \emptyset$ and that the boundary $\gamma := \partial B$ is a simple closed curve. 
Therefore $(O(x) \cup B) \cap O(y) = \emptyset$. 
Since the point $y$ is positively recurrent, there are points $\alpha, \beta \in O^+(y)$, a transverse closed arc $J$ whose endpoints are $\alpha$ and $\beta$, and there is a sequence $(y_n)_{n \in \Z_{\geq 0}}$ of points in $J \cap O^+(y)$ converging to $\beta$ such that $y_{n+1} \in O^+(y_n)$ and $J \cap B = \emptyset$. 
For any $n \in \Z_{\geq 0}$, from $z \in \omega(x)$, $y_{2n+1} \in \omega(x) \cap J$ and $\beta = \lim_{m \to \infty} y_m$, there are points $x_n, \widetilde{x}_n \in O^+(x) \cap \gamma$ and open orbit arcs $C_n \subset O^+(x)$ as in Figure~\ref{figure:decomp_001} satisfying the following properties: 
\\
{\rm(1)} $\widetilde{x}_n \in O^+(x_n)$ and $x_{n+1} \in O^+(\widetilde{x}_n)$. 
\\
{\rm(2)} The endpoints of $C_n$ are $x_n$ and $\widetilde{x}_n$. 
\\
{\rm(3)} $C_n \cap B = \{x_n, \widetilde{x}_n\}$ and $C_n \cap J_{y_{2n}, y_{2n+2}} \neq \emptyset$. 
\\
Here $J_{a,b}$ is the closed sub-arc of $J$ whose endpoints are $a$ and $b$. 
Put $C := \gamma \cup \bigsqcup_{n \in \Z_{\geq 0}} C_n$. 
\begin{figure}
\begin{center}
\includegraphics[scale=0.5]{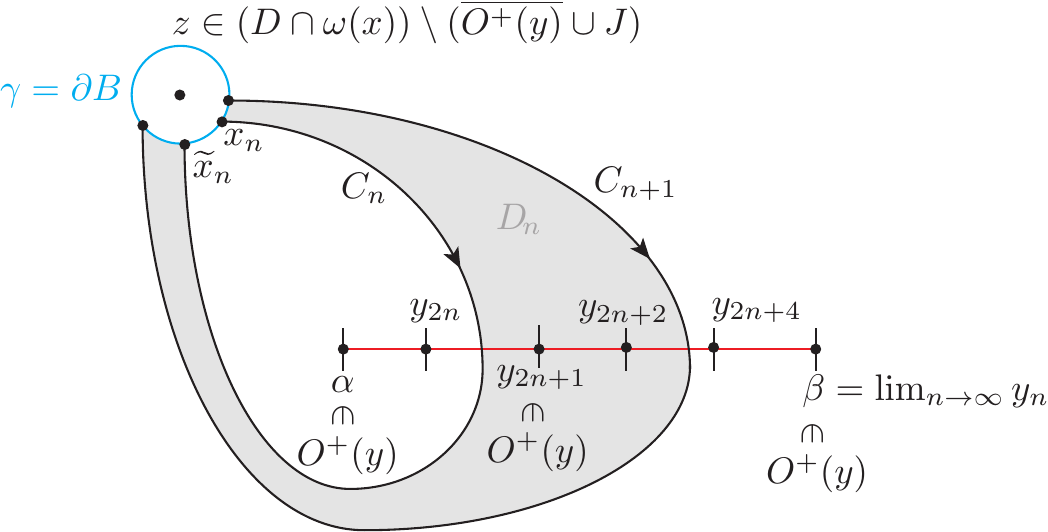}
\end{center}
\caption{The loop $\gamma = \partial B$, the orbits arc $C_n$, and the transverse closed arc $J$.}
\label{figure:decomp_001}
\end{figure}

\begin{claim}
We may assume that, for any $n \in \Z_{\geq 0}$, we have $C_n \cap J \subset J_{y_{2n}, y_{2n+2}}$ and there is the unique connected component $D_n$ of $S -C$ which is an open disk containing $y_{2n+1}$ and whose boundary is a loop consisting of $C_{n} \sqcup C_{n+1}$ and two sub-arcs in $\gamma$ such that $D_{n'}$ and $D_{n''}$ are disjoint for any $n' \neq n'' \in \Z_{\geq 0}$, by replacing $\alpha$ with some $y_{n_0} \in O^+(y_0)$ and taking subsequences of $(x_n)_{n \in \Z_{\geq 0}}$ and $(y_n)_{n \in \Z_{\geq 0}}$. 
\end{claim}
\begin{proof}
Denote by $D_n$ connected components of $S -C$ whose boundary contain $C_{n} \sqcup C_{n+1}$. 
Since $S$ is compact, the genus of $S$ is at most finite and the boundary $\partial S$ has at most finitely many connected components. 
Then there are at most finitely many connected components of $S -C$ that either have non-zero genus or have at least two boundary components. 
Therefore there are at most finitely many $n$ such that $D_n$ are not simply connected. 
By renumbering $(y_n)_{n \in \Z_{\geq 0}}$ if necessary, we may assume that any $D_n$ are simply connected and so pairwise disjoint open disks each of whose boundaries is a loop consisting of $C_{n} \sqcup C_{n+1}$ and two sub-arcs in $\gamma$. 
Then the union $D := \bigsqcup_{n \in \Z_{\geq 0}} D_n \sqcup C_{n+1}$ is an open disk. 
By construction, there is an open intercal $I_0 \subset J \cap D$ whose bounday consists of $\beta$ and a point $y_{n_0} \in O^+(y_0)$. 
Replacing $\alpha$ with $y_{n_0}$ and renumbering $(y_n)_{n \in \Z_{\geq 0}}$, we may assume that $\{ y_n \mid n \in \Z_{\geq 0} \} \subset I_0$ and that $\vert C_{n} \cap I_0 \vert \in \{0,1 \}$ for any $n \in \Z_{\geq 0}$. 
If $C_{n} \cap I_0$ is a singleton, denoted by $\{ x'_n \}$. 
By renumbering $(C_n)_{n \in \Z_{\geq 0}}$ if necessary, we may assume that $C_{n} \cap I_0 \neq \emptyset$ for any $n \in \Z_{\geq 0}$. 
From $x'_n, y_n \in I_0$ and $\lim_{n \to \infty} x'_n = \beta = \lim_{n \to \infty} y_n$, taking subsequences of $(x_n)_{n \in \Z_{\geq 0}}$ and $(y_n)_{n \in \Z_{\geq 0}}$, we may assume that $C_n \cap J \subset J_{y_{2n}, y_{2n+2}}$ for any $n \in \Z_{\geq 0}$. 
\end{proof}

By construction, the open disk $D_n$ contains $y_{2n+1} \in O^+(y)$. 
Since $O^+(y) \cap C = \emptyset$, we have $O^+(y) \cap \partial D_n = \emptyset$ and so $O^+(y) \subset D_n$. 
By $y_{2n+3} \in O^+(y) \cap D_{n+1}$, we have $y_{2n+3} \in D_n \cap D_{n+1} = \emptyset$, which is a contradiction. 
\end{proof}


\begin{lemma}\label{lem:ch_Q}
Let $v$ be a flow on a compact surface $S$. 
The $\omega$-limit set of a point is a Q-set if and only if it contains a non-closed recurrent orbit. 
\end{lemma}

\begin{proof}
Fix a point $x \in S$. 
If $\omega(x)$ is a Q-set, then the Q-set $\omega(x)$ contains a non-closed recurrent orbit by definition of Q-set. 
Conversely, suppose that $\omega(x)$ contains a non-closed recurrent orbit. 
If $x$ is positively recurrent, then $\omega(x)$ is a Q-set. 
Thus we may assume that $x$ is not positively recurrent (i.e. $x \notin \omega(x)$). 
Fix a non-closed recurrent orbit $O \subset \omega(x)$. 
If $O$ is locally dense,  then $O(x) \subseteq \overline{O} \subseteq \omega(x)$, which contradicts $x \notin \omega(x)$. 
Thus $O \subset \mathrm{E}(v)$. 
From \cite[Theorem~VI]{cherry1937topological}, there is a Poisson stable point $y \in \overline{O}$ with $\omega(y) = \overline{O(y)} = \overline{O}$. 
Since $y \in \omega(y) \cap \overline{O} \subseteq \omega(x)$, Lemma~\ref{pre_q-} implies that $\omega(x) = \omega(y) = \overline{O(y)} = \overline{O}$ is a transversely Cantor Q-set.
\end{proof}

We have the following statement. 

\begin{proposition}\label{prop:correspondence_q-set}
Let $v$ be a flow on a compact surface $S$ and $x \in S$ a point whose $\omega$-limit set is a Q-set. 
Then $\omega(x) = \overline{O(y)}$ for any non-closed recurrent points $y \in \omega(x)$. 
\end{proposition}

\begin{proof}
Fix a non-closed recurrent point $y \in \omega(x)$. 
Suppose $y \in \mathrm{E}(v)$. 
Lemma~\ref{lem:locally_dense_ch} and Lemma~\ref{lem:excptional_ch} imply that $\overline{O(y)}$ is a transversely Cantor Q-set and contains no locally dense orbits. 
From \cite[Theorem~VI]{cherry1937topological}, there is a Poisson stable point $z \in \overline{O(y)}$ with $\omega(z) = \overline{O(z)} = \overline{O(y)}$. 
By \cite[Proposition~2.2]{yokoyama2016topological}, we have $\hat{O}(z) = \overline{O(z)} \setminus (\Sv \sqcup \mathrm{P}(v)) = \overline{O(y)} \setminus (\Sv \sqcup \mathrm{P}(v)) \subseteq \mathrm{R}(v) - \mathrm{LD}(v) = \mathrm{E}(v)$. 
Since $z \in \omega(x) \cap \mathrm{E}(v)$ is positively recurrent, Lemma~\ref{pre_q-} implies that $\omega(x) = \omega(z) = \overline{O(z)} = \overline{O(y)}$ is a Q-set.

Suppose that $y \in \mathrm{R}(v) - \mathrm{E}(v) = \mathrm{LD}(v)$. 
Since $\emptyset \neq \mathop{\mathrm{int}} \overline{O(y)} \subseteq \omega(x)$, we have $O^+(x) \cap \mathop{\mathrm{int}} \overline{O(y)} \neq \emptyset$ and so $\omega(x) \subseteq \overline{O(y)}$. 
From $\overline{O(y)} \subseteq \omega(x)$, we obtain $\omega(x) = \overline{O(y)}$. 
\end{proof}

\subsubsection{Properties of connecting quasi-separatrices}

We have the following statement. 

\begin{lemma}\label{preQ02-}
A nontrivial quasi-Q-set consists of singular points and connecting quasi-separatrices. 
\end{lemma}

\begin{proof}
Let $\mathcal{M}$ be a nontrivial quasi-Q-set (i.e. quasi-Q-set that is not a Q-set). 
By time reversion if necessary, we may assume that $\mathcal{M}$ is the $\omega$-limit set of a non-positive-recurrent point. 
Then there is a non-positive-recurrent point $x \in S$ with $\omega(x) = \mathcal{M}$. 
If $\mathcal{M}$ contains a periodic point, then Lemma~\ref{lem:limit_cycle} implies that $\mathcal{M}$ is a semi-attracting limit cycle such that it intersects any essential closed transversal at most finitely many times, which contradicts the definition of quasi-Q-set.
Thus $\mathcal{M} \cap \Pv = \emptyset$. 
By Lemma~\ref{lq_q}, we have $\mathcal{M} \cap \mathrm{LD}(v) = \emptyset$. 

We claim that $\mathrm{E}(v) \cap \mathcal{M} = \emptyset$.
Indeed, assume that $\mathrm{E}(v) \cap \mathcal{M} \neq \emptyset$.
By \cite[Theorem~VI]{cherry1937topological}, there is a point $y \in \mathrm{E}(v) \cap \mathcal{M}$ with $\overline{O(y)} = \alpha(y) = \omega(y)$. 
Since $\omega(x) = \mathcal{M}$, Lemma~\ref{pre_q-} implies that the quasi-Q-set $\mathcal{M} = \omega(x) = \overline{O(y)}$ is a Q-set, which contradicts the hypothesis. 

Therefore $\mathcal{M} \subset \Sv \sqcup \mathrm{P}(v)$. 
Proposition~\ref{preQ03} implies that $\mathcal{M}$ consists of singular points and connecting quasi-separatrices. 
\end{proof}

We generalize the Poincar\'{e}-Bendixson theorem into one for a flow with arbitrarily many singular points on a compact surface and the Ma\v{i}er's description of recurrence as follows.

\subsection{Proof of Theorem~\ref{main:a}}
Proposition~\ref{prop:correspondence_q-set} implies assertion $(d)$. 
Let $v$ be a flow on a compact surface $S$ and $x$ a non-closed point. 
Lemma~\ref{lemma:a} implies that $\omega(x)$ is either a nowhere dense subset of singular points, a semi-attracting limit cycle, a quasi-semi-attracting limit quasi-circuit, a locally dense Q-set, or a quasi-Q-set that is not locally dense.
Proposition~\ref{preQ03} implies that any non-recurrent orbit in the $\omega$-limit set of a point is a connecting quasi-separatrix. 
This means that assertion $(b)$ holds. 
Corollary~\ref{cor:td} implies assertion $(c)$. 
By Lemma~\ref{lem:qc_ch}, a quasi-semi-attracting limit quasi-circuit either is the image of a circle or is not locally connected exclusively. 
From Corollary~\ref{preQ}, a quasi-Q-set that is not locally dense is either a transversely Cantor Q-set or a nontrivial quasi-Q-set. 
By Lemma~\ref{preQ02-}, a nontrivial quasi-Q-set corresponds to a quasi-Q-set that consists of singular points and connecting quasi-separatrices.
This completes the proof of Theorem~\ref{main:a}. 

\subsection{Proof of Theorem~\ref{main:rec}}

Let $v$ be a flow on a compact surface $S$ and a point $x \in \omega(z)$ for some point $z \in S$. 
Lemma~\ref{pos_rec} implies that if $\omega(x) \setminus \mathop{\mathrm{Cl}}(v) \neq \emptyset$ then $x$ is non-closed positively recurrent. 
Conversely, if $x$ is non-closed positively recurrent, then $x \in \omega(x) \setminus \mathop{\mathrm{Cl}}(v)$ and so $\omega(x) \setminus \mathop{\mathrm{Cl}}(v) \neq \emptyset$. 
This completes assertion {\rm (1)}. 

If $x$ is non-closed negatively recurrent, then $x \in \alpha(x) \setminus \mathop{\mathrm{Cl}}(v)$ and so $\alpha(x) \setminus \mathop{\mathrm{Cl}}(v) \neq \emptyset$. 
Conversely, suppose that $\alpha(x) \setminus \mathop{\mathrm{Cl}}(v) \neq \emptyset$. 
Then $x \notin \mathop{\mathrm{Cl}}(v)$ and so $x \in \mathrm{P}(v) \sqcup \mathrm{R}(v)$. 
%
We claim that $x$ is non-closed recurrent. 
Indeed, assume that $x$ is non-recurrent. 
Theorem~\ref{main:a} implies that $O(x)$ is a connecting quasi-separatrix and so that $\alpha(x) \subseteq \Sv$, which contradicts $\alpha(x) \setminus \mathop{\mathrm{Cl}}(v) \neq \emptyset$. 
Thus $x$ is non-closed recurrent. 
%
By \cite[Theorem~VI]{cherry1937topological}, the orbit class $\hat{O}(x)$ contains infinitely many Poisson stable orbits. 
Therefore there is a Poisson stable point $y \in \hat{O}(x)$ such that $x \in \alpha(y)$. 
The dual of Lemma~\ref{pos_rec} implies that $x$ is non-closed negatively recurrent. 
This completes assertion {\rm (2)}. 

Assertions {\rm (3)} and {\rm (4)} are followed from assertions {\rm (1)} and {\rm (2)}. 

\subsection{Proof of Corollary~\ref{main:rec02}}

%
%
By \cite[Theorem~VI]{cherry1937topological}, the orbit class of a non-closed recurrent point contains infinitely many Poisson stable orbits. 
Therefore Theorem~\ref{main:rec} and its dual statement imply that assertions {\rm(1)}--{\rm(3)} are pairwise equivalent. 

\subsection{Essential property of quasi-Q-sets}

We state the following essential property of quasi-Q-sets. 

\begin{lemma}\label{lem:ess_qqset}
Every quasi-Q-set of a flow on a compact surface is essential. 
\end{lemma}

\begin{proof}
Collapsing the boundary component of the surface $S$ into singletons if necessary,  
we may assume that $S$ is closed. 
Let $Q$ be a quasi-Q-set of a flow $v$ on $S$. 
Lemma~\ref{lem4-07a} implies that $Q$ is not a limit quasi-circuit. 
By the time reversing if necessary, we may assume that $Q$ is an $\omega$-limit set.
Fix a point $x \in S$ with $\omega(x) = Q$. 

Assume that $Q$ is inessential. 
Then there is an open disk $U \subset S$ which is a \nbd of $Q$. 
Since $\omega(x) = Q \subset U$, we have $O^+(v_T(x)) \cap U \neq \emptyset$ for any $T > 0$. 

\begin{claim}
There is a positive number $T>0$ with $O^+(v_T(x)) \subset U$. 
\end{claim}
\begin{proof} 
Assume $O^+(v_T(x)) \not\subset U$ for any $T > 0$. 
Then are strictly increasing sequence $(t_n)_{n \in \Z_{\geq 0}}$ of $\R_{>0}$ with $\lim_{n \to \infty} t_n = \infty$ such that $\{ v_{t_{n}}(x) \mid n \in \Z_{\geq 0} \} \subset S - U$. 
This implies that $\emptyset \neq \bigcap_{n =0}^\infty \overline{\{ v_{t_{n}}(x) \mid n \in \Z_{\geq 0} \}} \subset (S - U) \cap \bigcap_{n\in \mathbb{R}}\overline{\{v_t(x) \mid t > n\}} = (S - U) \cap \omega(x)$, which contradicts $\omega(x) = Q \subset U$. 
\end{proof} 

Replacing $x$ with a point in $O^+(x)$, we may assume that $O^+(x) \subset U$. 
Collapsing the boundary $\partial U$ into a singleton $p^*$, the resulting surface $U^*$ is a sphere. 
Considering the singleton as a singlar point, the resulting flow $v^*$ on the sphere $U^*$ from the restriction $v|_U$ contains $Q$ as the $\omega$-limit set of $x$ with respect to $v^*$ up to topological equivalent. 

\begin{claim}
$Q$ is not a limit quasi-circuit with respect to $v^*$.
\end{claim}
\begin{proof} 
Assume that $Q$ is a limit quasi-circuit with respect to $v^*$. 
There is a small collar $\A \subseteq U^*$ one of whose boundary components is $Q$. 
Since $U^*$ is a sphere, by $p^* \notin Q$, taking $\A$ small if necessary, we may assume that  $\A \subseteq U = U^* - \{ p^*\}$. 
From $\A \subseteq U \subset S$, we also obtain that $Q$ is a limit quasi-circuit with respect to $v$, which contradicts that  $Q$ is not a limit quasi-circuit with respect to $v$. 
\end{proof} 

Theorem~\ref{main:a} implies that $Q$ must be also quasi-Q-set with respect to $v^*$. 
This implies the existence of an essential closed transversal in the sphere $U^*$, which contradicts the simply connectivity of the sphere.
Thus $Q$ is essential.
\end{proof}

\subsection{Unboundedness of quasi-Q-sets}

We have the following unboundedness of quasi-Q-sets. 

\begin{lemma}\label{lem:unbdd_qqset}
Let $x \in S$ be a point whose $\omega$-limit set is a quasi-Q-set of a flow $v$ on a compact surface $S$ and $\pi \colon \widetilde{S} \to S$ the canonical unversal covering map. 
The following statements hold: 
\\
{\rm(1)} For any $y \in S$ with $\omega(y) = \omega(x)$, the positive orbit of $\widetilde{v}$ for any point in $\pi^{-1}(y)$ is unbounded, where $\widetilde{v}$ is the lift of $v$ on $\widetilde{S}$. 
\\
{\rm(2)} If there is the connected component $C_x$ of $S - \omega(x)$ containing $x$, then the boundary $\partial \widetilde{C_x} \subseteq \pi^{-1}(\omega(x))$ for any connected component $\widetilde{C_x}$ of $\pi^{-1}(C_x)$ is unbounded. 
\end{lemma}

\begin{proof}
Put $Q := \omega(x)$. 
Let $\widetilde{Q}$ be a lift of $Q$ on the universal cover $\widetilde{S}$ of $S$. 
Fix a closed transversal $T$ intersecting $Q$ infinitely many times. 
The closedness of $\omega(x)$ implies that the complement $S - \omega(x)$ is an open subset and so is an open surface whose lift to $\widetilde{S}$ is unbounded. 
Fix a connected component $C$ of the open surface $S - \omega(x)$. 

\begin{claim}
We may assume that $S$ is orientable. 
\end{claim}
\begin{proof} 
Since the vector field is lifted to the orientation double covering $S'$ of $S$, let $v'$ be the left of $v$ on $S'$. 
There is a lift $x'$ of $x$ whose $\omega$-limit set $\omega_{v'}(x')$ is a lift of $Q$ on $S'$ such that $\widetilde{Q}$ is the lift of $\omega_{v'}(x')$. 
Since any lift $T'$ on the orientation double covering $S'$ of the closed transversal $T$ of $v$ is also a closed transversal of $v'$, the lift $Q'$ is a quasi-Q-set. 
Because the point $x'$ is a lift of $x$, the connected component of $S' - \omega_{v'}(x')$ containing $x'$ is a lift of $C$. 
This means that we may assume that $S$ is orientable. 
\end{proof}

\begin{claim}\label{claim:Cx}
$\partial C \subseteq \omega(x)$. 
\end{claim}
\begin{proof} 
Since $C$ is the connected component of $S - \omega(x)$, we have that $\overline{C} \cap (S - \omega(x)) = C$ and so that $ (\overline{C} - C) \cap (S - \omega(x)) = \emptyset$. 
Then $\partial C = \overline{C} - C \subseteq \omega(x)$. 
\end{proof}

\begin{claim}\label{claim:y_in_omega}
For any positively recurrent point $y \in S$ whose positive orbit intersects the essential closed transversal $T$ infinitely many times, 
the positive orbit  $O^+_{\widetilde{v}}(\widetilde{y})$ is unbounded.
\end{claim}
\begin{proof} 
Fix a positively recurrent point $y \in S$ which intersects the essential closed transversal $T$ infinitely many times. 
Because $T$ is essential, the genus of the surface $S$ is positive and any connected components of the preimage $\pi^{-1}(T)$ are simple curves between ideal boundary points of $\widetilde{S}$. 
Since $\vert O^+(y) \cap T \vert = \infty$, the positive orbit $O^+_{\widetilde{v}}(\widetilde{y}) \subset \pi^{-1}(C)$ for any point $\widetilde{y} \in \pi^{-1}(y)$ is the lift of $O^+(y)$ which intersects infinitely many simple curves $\widetilde{T}_i$ between ideal boundary points of $\widetilde{S}$, where $\widetilde{T}_i$ are lifts of $T$.
Choose infinitely many such simple curves $\widetilde{T}_i$ and points $\widetilde{y}_i \in O^+_{\widetilde{v}}(\widetilde{y}) \cap \widetilde{T}_i$ with $O^+_{\widetilde{v}}(\widetilde{y}_i) \subsetneq O^+_{\widetilde{v}}(\widetilde{y}_{i+1})$. 
This means that the positive orbit  $O^+_{\widetilde{v}}(\widetilde{y})$ is unbounded. 
\end{proof}

Fix a point $y \in S$ with $\omega(y) = \omega(x)$. 
Then $\vert O^+(y) \cap T \vert = \infty$.

\begin{claim}\label{claim:unbd_positive}
The positive orbit of $\widetilde{v}$ for any point in $\pi^{-1}(y)$ is unbounded. 
\end{claim}
\begin{proof} 
Fix a point $\widetilde{y} \in \pi^{-1}(y)$. 
Claim~\ref{claim:y_in_omega} implies that we may assume that $y \notin \omega(y)$. 
Choose the connected component $\widetilde{C}$ of an open surface $\pi^{-1}(C)$ with $\widetilde{y} \in \widetilde{C}$. 
Since any connected components of $\pi^{-1}(C)$ are open subsets, we obtain $\overline{\widetilde{C}} \cap \pi^{-1}(C) = \widetilde{C}$ and so $\partial \widetilde{C} \cap \pi^{-1}(C) = (\overline{\widetilde{C}} - \widetilde{C}) \cap \pi^{-1}(C) = \emptyset$. 
Then $\pi(\partial \widetilde{C}) \cap C = \emptyset$. 
From $\partial \widetilde{C} \subseteq \pi^{-1}(\overline{C})$, we have $\pi(\partial \widetilde{C}) \subseteq \overline{C} - C = \partial C \subseteq \omega(x)$, because of Claim~\ref{claim:Cx}. 
Because $T$ is essential, the genus of the surface $S$ is positive and any connected components of the preimage $\pi^{-1}(T)$ are simple curves between the ideal boundary of $\widetilde{S}$. 
Since $\vert O^+(y) \cap T \vert = \infty$, the positive orbit $O^+_{\widetilde{v}}(\widetilde{y}) \subset \widetilde{C}$ is the lift of $O^+(y)$ which intersects infinitely many simple curves $\widetilde{T}_i$ between the ideal boundary of $\widetilde{S}$.
Choose infinitely many such simple curves $\widetilde{T}_i$ and points $\widetilde{y}_i \in O^+_{\widetilde{v}}(\widetilde{y}) \cap \widetilde{T}_i$ with $O^+_{\widetilde{v}}(\widetilde{y}_i) \subsetneq O^+_{\widetilde{v}}(\widetilde{y}_{i+1})$,  where $\widetilde{T}_i$ are lifts of $T$. 
This means that the positive orbit  $O^+_{\widetilde{v}}(\widetilde{y})$ is unbounded. 
\end{proof}

Suppose that there is the connected component $C_x$ of $S - \omega(x)$ containing $x$. 
Then $x \notin \omega(x)$ and $\vert O^+(x) \cap T \vert = \infty$. 
Since $x \notin \omega(x)$, the positive orbit of $\widetilde{v}$ for any point in $\pi^{-1}(x)$ is contained in $\pi^{-1}(C_x)$. 
Claim~\ref{claim:unbd_positive} implies that the positive orbit of $\widetilde{v}$ for any point in $\pi^{-1}(x)$ is unbounded. 
Since any connected component of $\pi^{-1}(C_x)$ intersects $\pi^{-1}(x)$ and the positive orbit of $\widetilde{v}$ for any point in $\pi^{-1}(x)$ is contained in $\pi^{-1}(C_x)$, any connected component of $\pi^{-1}(C_x)$ is unbounded and so is the boundary $\partial \widetilde{C_x}$ for any connected component $\widetilde{C_x}$ of $\pi^{-1}(C_x)$. 
\end{proof}

Note that the positive orbit in Lemma~\ref{lem:unbdd_qqset} (1) has asymptotic direction in the sense of \cite{aranson1996maier} (see the proof of \cite[Theorem~3.1]{aranson1996maier}).
However, the author does not know whether the subset of absolutes (see \cite{aranson1996maier} for the definition) which are contained in the limits of curves on a connected component of the preimage $\pi^{-1}(\omega(x)) \subset \widetilde{S}$ of the quasi-Q-set $\omega(x)$ in the previous lemma consists of exactly two elements. 

\subsection{Finiteness of quasi-Q-sets}

To state finiteness, recall the end completion as follows. 

\subsubsection{Direct system}

A binary relation $\leq$ on a set $P$ is a pre-order (or quasiorder) if it is reflexive (i.e. $a \leq a$ for any $a \in P$) and transitive (i.e. $a \leq c$ for any $a, b, c \in P$ with $a \leq b$ and $b \leq c$).
The pair $(P, \leq)$ is called a pre-ordered set. 
A pre-ordered set $(P, \leq)$ is a {\bf directed set} if for any elements $a,b \in P$ there is an element $c \in P$ with $a \leq c$ and $b \leq c$. 

For a directed set $(\Lambda, \leq)$, a family $\{K_\lambda\}_{\lambda \in \Lambda}$ of sets indexed by $\Lambda$, and a family $\{i_{\lambda, \lambda'} \colon K_{\lambda} \to K_{\lambda'} \mid \lambda, \lambda' \in \Lambda, \lambda \leq \lambda' \}$, a pair $(\{K_\lambda\}, \{i_{\lambda, \lambda'} \})$ is a  {\bf direct system} if $i_{\lambda, \lambda} = 1_{K_\lambda}$ and $i_{\lambda, \lambda''} = i_{\lambda, \lambda'} \circ i_{\lambda', \lambda''}$ for any $\lambda, \lambda', \lambda'' \in \Lambda$, where $1_{K_\lambda}$ is the identity map on $K_\lambda$. 

\subsubsection{End completion of a topological space}
For a topological space $Y$, consider a direct system $(\{K_\lambda\}, \{i_{\lambda, \lambda'}\})$ of compact subsets $K_\lambda$ of $Y$ and inclusion maps $i_{\lambda, \lambda'} \colon K_{\lambda} \to K_{\lambda'}$ such that the interiors of $K_\lambda$ cover $Y$.  
There is a corresponding inverse system $\{ \pi_0( Y - K_\lambda ) \}$, where $\pi_0(Z)$ denotes the set of connected components of a topological space $Z$. 
Then the {\bf set of ends} of $Y$ is defined to be the inverse limit of this inverse system. 
Notice that $Y$ has one end $x_{\mathcal{U}}$ for each sequence $\mathcal{U} := (U_i)_{i \in \mathbb{Z}_{>0}}$ with $U_i \supseteq U_{i+1}$ such that $U_i$ is a connected component of $Y - K_{\lambda_i}$ for some $\lambda_i$. 
Considering the disjoint union $Y_{\mathrm{end}}$ of $Y$ and  $\{ \pi_0( Y - K_\lambda ) \}$ as set, a subset $V$ of the union $Y_{\mathrm{end}}$ is an open \nbd of an end $x_{\mathcal{U}}$ if there is some $i \in \mathbb{Z}_{>0}$ such that $U_i \subseteq V$. 
Then the resulting topological space $Y_{\mathrm{end}}$ is called the {\bf end completion} (or end compactification) of $Y$. 
Note that the end completion is not compact in general. 
Moreover, the surface $S_{\mathrm{col}}$ is the end completion of $S - \mathop{\mathrm{Sing}}(v)$. 
From Theorem~3~\cite{richards1963classification}, all connected surfaces of finite genus and finitely many boundary components are homeomorphic to the resulting surfaces from compact surfaces by removing closed totally disconnected subsets. 
Therefore the end compactification $S_{\mathrm{end}}$ of a connected surface $S$ of finite genus and finitely many boundary components is a compact surface.

We generalize the Ma\v{i}er's work \cite{markley1969poincare} for Q-sets into quasi-Q-sets as follows. 

\begin{proposition}\label{prop:num_qqset}
The number of quasi-Q-sets of a flow on an orientable compact surface is at most the genus.
\end{proposition}

\begin{proof}
Let $v$ be a flow on an orientable compact surface $S$. 
Taking the double of the surface $S$ if necessary, we may assume that $S$ is closed. 
Denote by $g$ the genus of $S$. 
By induction for $g$, we show the assertion. 
Lemma~\ref{lem:ess_qqset} implies that any flow on the sphere has no quasi-Q-sets. 
Thus, we may assume that $g > 0$ and that $v$ has a quasi-Q-set. 
Let $Q$ be a quasi-Q-set. 
From the closedness of $Q$, the complement $S -Q$ is an open subset and so an orientable open surface. 
Since any connected components of $S -Q$ are orientable open surfaces, the end completions of the connected components of $S-Q$ are orientable closed surfaces. 
Let $S'$ be the disjoint union of such end completions. 
By the construction of the end completion, we have that $S' - \mathcal{E} = S - Q$, where $\mathcal{E}$ is the set of ends. 

We claim that the sum $g'$ of genera of connected components of $S'$ is less than the genus $g$ of $S$. 
Indeed, assume $g' = g$. 
By the construction of the end completion, there is a disjoint union $D' \subset S'$ of finitely many closed disks whose interior contains the set $\mathcal{E}$ of ends. 
Then $S' - D' \subset S - Q$. 
Since $S' - \mathcal{E} = S - Q$, there is a closed \nbd $D \subset S$ of $Q$ such that $D' - \mathcal{E} = D - Q$ and so that $S' - D' = S - D$. 
Then the end completion of $S' - D' = S - D$ has the genus $g' = g$. 
Since $\partial D' = \partial (S' - D') = \partial (S - D) = \partial D$ is a disjoint union of finitely many loops, any connected component of $S - \partial D$ intersecting the interior $\mathop{\mathrm{int}} D$ of $D$ is a simply connected open subset containing $Q$. 
Because simply connected open surfaces are open disks, the open surface $\mathop{\mathrm{int}} D$ is a disjoint union of finitely many open disks containing $Q$. 
This means that $Q$ is inessential, which contradicts the essential property of $Q$. 

By inductive hypothesis, the resulting flow on $S'$ has at most $g'$ quasi-Q-sets. 
This implies that the number of quasi-Q-sets of $v$ is at most $g' + 1 \leq g$. 
\end{proof}

We partially generalize Markley's work \cite{markley1970number} for Q-sets into quasi-Q-sets as follows. 

\begin{proposition}\label{prop:num_qqset_02}
The number of quasi-Q-sets of a flow on a nonorientable compact surface is at most $p-1$, where $p$ is the number of nonorientable genus.
\end{proposition}

\begin{proof}
Let $v$ be a flow on a nonorientable compact surface $S$. 
Taking the double of the surface $S$ if necessary, we may assume that $S$ is closed. 
By Gutierrez's smoothing theorem~\cite{gutierrez1986smoothing}, the flow $v$ is topologically equivalent to a $C^1$-flow and so is generated by an integrable continuous vector field on $S$ which is integrable. 
Since the vector field is lifted to the orientation double covering $S'$ of $S$, let $v'$ be the left of $v$ on $S'$. 
Then the genus of $S'$ is $p-1$, where $p$ is the number of nonorientable genus of $S$.
Proposition~\ref{prop:num_qqset} implies the assertion. 
 \end{proof}

\subsection{Topological characterizations of (non-trivial) quasi-Q-sets}

We have the following observation. 

\begin{lemma}\label{lem:iso_qc}
If there are a non-periodic point $x$ and a transverse open arc $I$ with $\vert I \cap \omega(x) \vert = 1$, then $\omega(x)$ is either a limit cycle or a limit quasi-circuit. 
\end{lemma}

\begin{proof}
Suppose that there is a transverse open arc $I$ with $\vert I \cap \omega(x) \vert = 1$. 
Then $I$ contains a transverse closed arc $J: [-1,0] \to I$ with $J(-1) \in O^+(x)$ and $\{ J(0) \} = J([-1,0]) \cap \omega(J(-1)) \subset \overline{J([-1,0]) \cap O^+(J(-1))}$. 
Lemma~\ref{lem-quasi-limit-cricuits} implies $\omega(x)$ is either a limit cycle or a limit quasi-circuit. 
\end{proof}

The $\omega$-limit set $Q$ is {\bf transversely Cantor set at a non-singular point} if there are a non-singular point in $Q$ and its open \nbd $U$ such that the intersection $U \cap Q$ is the product of a Cantor set and an open interval.
We topologically characterize a quasi-Q-set as follows. 

\begin{proposition}\label{prop:ch_qqset}
An $\omega$-limit set of a point of a flow on a compact surface is a quasi-Q-set if and only if it is either locally dense or transversely Cantor set at a non-singular point. 
In the above cases, the $\omega$-limit set is essential and contains non-closed orbits. 
\end{proposition}

\begin{proof}
Let $Q$ be an $\omega$-limit set of a point $x$ of a flow on a compact surface $S$. 
If $Q$ is a closed orbit, then $Q$ is neither a quasi-Q-set, nor a locally dense subset, nor a transversely Cantor set at a non-singular point. 
Thus we may assume that $Q$ is not a closed orbit. 
Then $x$ is not closed. 

We claim that we may assume that $Q$ is not locally dense. 
Indeed, if $Q$ is locally dense then Lemma~\ref{lem:equiv_locally_Q} and Corollary~\ref{preQ} imply that $Q$ is a quasi-Q-set. 
Conversely, if $Q$ is a locally dense quasi-Q-set, then Lemma~\ref{lem:ess_qqset} implies the essential property and Lemma~\ref{lem:equiv_locally_Q} implies the existence of non-closed orbits. 
Thus the assertion holds if $Q$ is locally dense. 

Suppose that $Q$ is a quasi-Q-set. 
Then there is a closed transversal $\gamma$ which intersects $Q$ infinitely many times.
Moreover, the intersection $\gamma \cap Q$ is closed and so has an accumulation point. 
Since any quasi-Q-sets are neither limit cycles nor limit quasi-circuits, Lemma~\ref{lem:iso_qc} implies that $\vert I \cap \omega(x) \vert = \infty$ for any transverse open arc $I$ containing a point in $\omega(x) = Q$.  
This means that any point in $Q \cap \gamma$ is an accumulation point of the intersection $Q \cap \gamma$. 
Therefore $Q \cap \gamma$ is perfect and totally disconnected because of the absence of local density of $Q$.  
Since a Cantor set is characterized as a compact metrizable perfect totally disconnected space, 
there is a small neighborhood $U$ of a non-singular point of $Q \cap \gamma$ such that $Q \cap U$ is a product of an open interval and a Cantor set. 
This means that $Q$ is a transversely Cantor set at a non-singular point. 
Lemma~\ref{lem:ess_qqset} implies the essential property of $Q$. 
Theorem~\ref{main:a} implies that $Q$ contains non-recurrent points and so non-closed orbits. 
%

Conversely, suppose that $Q$ is a transversely Cantor set at a non-singular point. 
By definition of transversely Cantor set at a non-singular point, the $\omega$-limit set $Q$ contains non-closed orbits. 
Therefore $Q$ is neither a subset of the singular point set nor a limit cycle. 
By Lemma~\ref{lem:acc_lc}, the existence of a transversely Cantor set at a non-singular point implies that $Q$ is not a quasi-circuit. 
If $Q$ is a Q-set, then Corollary~\ref{preQ} implies that $Q$ is a quasi-Q-set. 
Thus, we may assume that $Q$ is not a Q-set. 
Theorem~\ref{main:a} implies that $Q$ is a quasi-Q-set that consists of singular points and non-recurrent points.  
\end{proof}

We topologically characterize a non-trivial quasi-Q-set as follows.

\begin{proposition}\label{prop:ch_qqset_nontrivial}
An $\omega$-limit set $Q$ of a flow on a compact surface is a non-trivial quasi-Q-set if and only if it satisfies the following two conditions:
\\
{\rm(1)} The $\omega$-limit set $Q$ consists of singular points and non-recurrent orbits. 
\\
{\rm(2)} The $\omega$-limit set $Q$ is a transversely Cantor set at a non-singular point. 
\\
In the above cases, the $\omega$-limit set $Q$ is essential and contains non-recurrent orbits. 
\end{proposition}

\begin{proof}
Let $Q$ be an $\omega$-limit set of a point $x$ of a flow on a compact surface $S$. 
Suppose that $Q$ is a non-trivial quasi-Q-set. 
Theorem~\ref{main:a} implies that $Q$ consists of singular points and non-recurrent points.  
Proposition~\ref{prop:ch_qqset} implies that assertion (2) follows from the non-trivial quasi-Q-set property. 
Conversely, suppose that $Q$ consists of singular points and non-recurrent orbits and is a transversely Cantor set at a non-singular point. 
Proposition~\ref{prop:ch_qqset} implies that $Q$ is a quasi-Q-set. 
Since any Q-sets contain non-closed recurrent points, the non-existence of non-closed recurrent points implies that $Q$ is not a Q-set. 
%
\end{proof}

Theorem~\ref{main:a} and Proposition~\ref{prop:ch_qqset} imply the following statement on surfaces without genus. 

\begin{corollary}\label{cor:nonex_qQset}
The following statements hold for a flow with arbitrarily many singular points on a compact surface that is contained in a sphere or a projective plane: 
\\
{\rm(a)} The $\omega$-limit set of any non-closed orbit is one of the following exclusively:
\begin{quote}
\setlength{\leftskip}{-25pt}$(1)$ A nowhere dense subset of singular points.
\\
$(2)$ A semi-attracting limit cycle.
\\
$(3)$ A quasi-semi-attracting limit quasi-circuit that is the image of a circle.
\\
$(4)$ A quasi-semi-attracting limit quasi-circuit that is not locally connected. 
\end{quote}
{\rm(b)} Every non-closed orbit in the $\omega$-limit set of a point is a connecting quasi-separatrix. 
\\
{\rm(c)} If the singular point set is totally disconnected, then any non-recurrent orbits in the $\omega$-limit set of a point are connecting separatrices.  
\end{corollary}

The previous corollary shows the non-existence of quasi-Q-sets on a sphere and a projective plane.

\subsection{Existence of uncountably many intersections of non-trivial quasi-Q-sets and boundary components of the singular point set}

We have the following observation to show the existence of uncountably many intersections of non-trivial quasi-Q-sets and boundary components of the singular point set.

\begin{lemma}\label{lem:q-Q_02}
Let $v$ be a flow on a compact surface $S$ and $x \in S$ a point with $\omega(x) \subseteq \mathop{\mathrm{Sing}}(v) \sqcup \mathrm{P}(v)$. 
For any transverse closed arc $T$ with respect to $v$ and any point $x' \in T \cap \omega(x)$, the set of points $x'' \in T \cap \omega(x)$ with $\omega_{v_{\mathrm{col}}}(x') = \omega_{v_{\mathrm{col}}}(x'')$ and $\alpha_{v_{\mathrm{col}}}(x') =\alpha_{v_{\mathrm{col}}}(x'')$ is finite. 
\end{lemma}

\begin{proof}
Because the assertion holds for $v$ if one holds for $v_{\mathrm{col}}$, replacing $v$ with $v_{\mathrm{col}}$ if necessary, we may assume that $\mathop{\mathrm{Sing}}(v)$ is totally disconnected. 

Assume that there are a transverse closed arc $T$ and an infinite subset $T' \subseteq T \cap \omega(x)$ such that $\omega(x') = \omega(x'')$ and $\alpha(x') =\alpha(x'')$ for any points $x', x'' \in T' \subseteq T \cap \omega(x)$. 
Then $x$ is non-singular. 
Since $\omega(x) \subseteq \mathop{\mathrm{Sing}}(v) \sqcup \mathrm{P}(v)$, if $x \in \omega(x)$ then $x \in \mathrm{R}(v) \cap \omega(x) \subseteq \mathrm{R}(v) \cap (\mathop{\mathrm{Sing}}(v) \sqcup \mathrm{P}(v)) = \emptyset$, which is a contradiction. 
Thus the point $x$ is not positive-recurrent and so $O(x) \cap \omega(x) = \emptyset$. 
Put $\omega := \omega(x')$ and $\alpha :=\alpha(x')$ for any points $x' \in T' \subseteq T \cap \omega(x)$. 
By Theorem~\ref{main:a}, any non-recurrent orbits in $\omega(x)$ are connecting separatrices and so are the orbits $O(x')$ for any non-singular points $x' \in T' \subseteq T \cap \omega(x)$. 
Since $T$ is compact, there is a sequence $(x_i)_{i \in \Z_{\geq 0}}$ of points $x_i \in T'$ which converges to an accumulation point $y \in T \cap \omega(x)$. 
Since $O(x_i)$ are connecting separatrices, the intersection $O(x_i) \cap \{ x_j \mid j \in \Z_{\geq 0}\}$ is at most finite. 
Taking a subsequence of $(x_i)_{i \in \Z_{\geq 0}}$, we may assume that $O(x_i)$ and $O(x_j)$ are disjoint for any $i \neq j$.    
Put $O_i := O(x_i) \subset \omega(x)$. 
Write $\Gamma := \{ \alpha, \omega \} \sqcup \bigsqcup_{i \in \Z_{\geq 0}} O_i$. 
Since $S$ is compact and the genus is finite, by the connectivity of $\Gamma = \{ \alpha, \omega \} \sqcup \bigsqcup_{i \in \Z_{\geq 0}} O_i$, there are at most finitely many connected components of $S -\Gamma$ that either have non-zero genus or have at least two boundary components. 
Therefore, taking a subsequence, we may assume that the pair of two orbits $O_i$ and $O_j$ for any $i \neq j \in \Z_{\geq 0}$ is homotopic relative to $\{ \alpha, \omega \}$ to each other and that the union $\gamma_i := \{ \alpha, \omega \} \sqcup O_i \sqcup O_{i+1}$ for any $i \in \Z_{\geq 0}$ bounds an open disk $B_i$. 
Then the disjoint unions $B_i \sqcup O_{i+1} \sqcup B_{i+1}$ are \nbds of $O_{i+1}$. 
By construction, any orbits intersecting some $B_j$ are contained in $B_j$. 
Since $O_1, O_2 \subset \omega(x)$, we have that $O(x) \cap (B_1 \sqcup B_2) \neq \emptyset$ and $O(x) \cap (B_2 \sqcup B_3) \neq \emptyset$ and so that $O(x) \subset B_2$. 
By $\overline{B_2} \cap O_4 = \emptyset$, we obtain $\omega(x) \subseteq \overline{O(x)} \subseteq \overline{B_2} \subset S - O_4$. 
This implies that $O_4 \cap \omega(x) = \emptyset$, which contradicts $O_4 \subset \omega(x)$. 
\end{proof}

The previous lemma implies the following statement. 

\begin{lemma}\label{lem:q-Q_uncountable}
A non-trivial quasi-Q-set of a flow on a compact surface intersects uncountably many connected components of the singular point set, and contains uncountably many connecting quasi-separatrices. 
\end{lemma}

\begin{proof}
Because the assertion holds for $v$ if one holds for $v_{\mathrm{col}}$, replacing $v$ with $v_{\mathrm{col}}$ if necessary, we may assume that $\mathop{\mathrm{Sing}}(v)$ is totally disconnected. 
Let $Q$ be a non-trivial quasi-Q-set of a flow $v$ on $S$. 
By the time reversing if necessary, we may assume that $Q$ is an $\omega$-limit set of a point $x \in S$.
Proposition~\ref{prop:ch_qqset_nontrivial} implies that $Q$ consists of singular points and non-recurrent orbits and is a transversely Cantor set at a non-singular point. 
By definition of transversely Cantor set at a non-singular point, there is a transverse closed arc $T$ such that the intersection $T \cap \omega(x)$ is a Cantor set.
By Theorem~\ref{main:a}, any non-recurrent orbits in $\omega(x)$ are connecting separatrices. 
Since any connecting separatrices intersects any transverse closed arcs at most finitely many times, the $\omega$-llimit set $\omega(x)$ contains uncountably many connecting separatrices, and the intersection $T \cap \omega(x)$ contains a subset $T'$ which consists of uncountably many non-recurrent points such that $O(x') \neq O(x'')$ for any $x' \neq x'' \in T'$. 

Assume that $Q = \omega(x)$ contains at most countably many singular points.  
Since $T'$ contains uncountable points, the countable exsitence of singular points implies that there is singular point $\alpha \in Q$ such that $\{ x_\lambda \in T' \mid \alpha = \alpha(x_\lambda) \}$ is uncountable.  
Similarly, the countable exsitence of singular points implies that there is singular point $\omega \in Q$ such that $\{ x_\lambda \in T' \mid \alpha = \alpha(x_\lambda), \omega = \omega(x_\lambda) \}$ is uncountable, which contradicts the non-existence of such the transverse closed arc $T$ because of Lemma~\ref{lem:q-Q_02}.  
\end{proof}

\subsubsection{Characterization of non-triviality of quasi-Q-sets}

We have the following characterization of the non-triviality of quasi-Q-sets. 

\begin{proposition}\label{lem:ch_qqset_trivial}
A quasi-Q-set of a flow on a compact surface is non-trivial if and only if it contains no orbit whose closure is the quasi-Q-set.  
\end{proposition}

\begin{proof}
By Proposition~\ref{prop:correspondence_q-set}, each Q-set contains an orbit whose closure is the Q-set.  
Therefore any quasi-Q-set $Q$ of a flow on a compact surface which contains no orbit $O$ with $\overline{O} = Q$ is not a Q-set and so is non-trivial.

Conversely, Lemma~\ref{lem:q-Q_uncountable} implies that any non-trivial quasi-Q-set consists of singular points and uncountably many connecting quasi-separatrices. 
This implies that each non-trivial Q-set contains no orbit whose closure is the quasi-Q-set.  
\end{proof}

\section{Poincar\'{e}-Bendixson theorem for a non-compact surface and generalization of Ma\v{i}er's description of recurrence}

In this section, we generalize the Poincar\'{e}-Bendixson theorem to one for a flow with arbitrarily many singular points on a surface of finite genus and finitely many boundary components, which is homeomorphic to a closed surface with punctures. 
Moreover, we also characterize the recurrence, which are generalizations of the Ma\v{i}er's description of recurrence for such a flow. 

\subsection{Poincar\'{e}-Bendixson theorem for a flow with arbitrarily many singular points on a surface of finite genus and finitely many boundary components}

\subsubsection{Concepts for flows on {\rm(}possibly non-compact{\rm)} surfaces}

For a flow $v$ on a surface $S$ of finite genus and finitely many boundary components, considering ends to be singular points, we obtain the resulting flow $v_{\mathrm{end}}$ on a surface $S_{\mathrm{end}}$ which is a union of compact surfaces. 
A non-recurrent orbit on $S$ is a {\bf virtual quasi-separatrix} if it is a connecting quasi-separatrix on $S_{\mathrm{end}}$ with respect to $v_{\mathrm{end}}$. 
A non-recurrent orbit on $S$ is a {\bf virtual separatrix} if it is a connecting separatrix on $S_{\mathrm{end}}$ with respect to $v_{\mathrm{end}}$. 
An invariant subset on $S$ is a {\bf quasi-semi-attracting limit virtual quasi-circuit} if it is the resulting subset from a quasi-semi-attracting limit quasi-circuit on $S_{\mathrm{end}}$ with respect to $v_{\mathrm{end}}$ by removing all the ends. 

\subsubsection{Poincar\'{e}-Bendixson theorem for flows on {\rm(}possibly non-compact{\rm)} surfaces}

By taking end completions, Theorem~\ref{main:a} and Lemma~\ref{lem:q-Q_uncountable} imply the following Poincar\'{e}-Bendixson theorem for a flow with arbitrarily many singular points on a surface of finite genus and finitely many boundary components. 

\begin{theorem}\label{main:b}
The following statements hold for a flow with arbitrarily many singular points on a surface of finite genus and finitely many boundary components: 
\\
{\rm(a)} The $\omega$-limit set of any non-closed orbit is one of the following exclusively:
\begin{quote}
\setlength{\leftskip}{-25pt}
$(1)$ A {\rm(}possibly empty{\rm)} nowhere dense subset of singular points.
\\
$(2)$ A semi-attracting limit cycle.
\\
$(3)$ A quasi-semi-attracting limit virtual quasi-circuit. 
\\
$(4)$ A locally dense Q-set. 
\\
$(5)$ A transversely Cantor Q-set. 
\\
$(6)$ A quasi-Q-set that consists of uncountably many singular points and uncountably many non-recurrent points.  
\end{quote}
{\rm(b)} Any non-recurrent orbit in the $\omega$-limit set of a point is a virtual quasi-separatrix. 
\\
{\rm(c)} If the singular point set is totally disconnected, then any non-recurrent orbits in the $\omega$-limit set of a point are virtual separatrices.  
\\
{\rm(d)} If the $\omega$-limit set of a point is a Q-set, then the Q-set corresponds to the orbit closure of any non-closed recurrent point in the Q-set. 
\end{theorem}

In the previous theorem, notice that the $\omega$-limit set of a point is empty if and only if it is a nowhere dense subset of singular points. 
Theorem~\ref{main:b} and Corollary~\ref{cor:nonex_qQset}
imply the following statement on surfaces without genus. 

\begin{corollary}
The following statements hold for a flow with arbitrarily many singular points on a surface which is contained in a sphere or a projective plane: 
\\
{\rm(a)} The $\omega$-limit set of any non-closed orbit is one of the following (the possibilities are mutually exclusive):
\begin{quote}
\setlength{\leftskip}{-25pt}
$(1)$ A {\rm(}possibly empty{\rm)} nowhere dense subset of singular points.
\\
$(2)$ A semi-attracting limit cycle.
\\
$(3)$ A quasi-semi-attracting limit virtual quasi-circuit. 
\end{quote}
{\rm(b)} Any non-recurrent orbit in the $\omega$-limit set of a point is a virtual quasi-separatrix. 
\\
{\rm(c)} If the singular point set is totally disconnected, then any non-recurrent orbits in the $\omega$-limit set of a point are virtual separatrices.  
\end{corollary}

\subsection{Topological characterizations of non-closed recurrence for surfaces}

For a flow $v$ on a connected surface with finite genus and finitely many boundary components, the end completion $S_{\mathrm{end}}$ is a compact connected surface, and the resulting flow $v_{\mathrm{end}}$ on it can be obtained by adding exactly new singular points.  
This implies that Theorem~\ref{main:rec} holds for a flow with arbitrarily many singular points on a surface with finite genus and finitely many boundary components as follows. 

\begin{theorem}\label{main:c}
Let $v$ be a flow on a surface $S$ with finite genus and finitely many boundary components. 
The following statements hold for a point $x \in \omega(z)$ for some point $z \in S$: 
\\
{\rm (1)} $\omega(x) \setminus \mathop{\mathrm{Cl}}(v) \neq \emptyset$ if and only if $x$ is non-closed positively recurrent. 
\\
{\rm (2)} $\alpha(x) \setminus \mathop{\mathrm{Cl}}(v) \neq \emptyset$ if and only if $x$ is non-closed negatively recurrent. 
\\
{\rm (3)} $(\omega(x) \cup \alpha(x)) \setminus \mathop{\mathrm{Cl}}(v) \neq \emptyset$ if and only if $x$ is non-closed recurrent.
\\
{\rm (4)} $\omega(x) \setminus \mathop{\mathrm{Cl}}(v) \neq \emptyset$ and $\alpha(x) \setminus \mathop{\mathrm{Cl}}(v) \neq \emptyset$ if and only if $x$ is non-closed Poisson stable.
\end{theorem} 

\begin{proof}
As mentioned above, taking the end completion $S_{\mathrm{end}}$ of $S$, the resulting surface $S_{\mathrm{end}}$ is a compact connected surface and the resulting flow $v_{\mathrm{end}}$ is obtained by adding singular points. 
Therefore the conditions that $\omega(x) \setminus \mathop{\mathrm{Cl}}(v) \neq \emptyset$ and $\alpha(x) \setminus \mathop{\mathrm{Cl}}(v) \neq \emptyset$ are invariant under taking end completions. 
This means that Theorem~\ref{main:rec} holds for a flow on a connected surface $S$ with finite genus and finitely many boundary components.  
\end{proof}

Corollary~\ref{main:rec02} implies the following topological characterizations of non-closed recurrence for a flow with arbitrarily many singular points on a surface of finite genus and finitely many boundary components. 

\begin{corollary}\label{cor:rec02}
Let $v$ be a flow on a surface $S$ with finite genus and finitely many boundary components. 
The following are equivalent for a point $x \in S$: 
\\
{\rm (1)} The point $x$ is non-closed positively recurrent. 
\\
{\rm (2)} $\omega(x) \setminus \mathop{\mathrm{Cl}}(v) \neq \emptyset$ and there is a point $z \in S$ with $x \in \omega(z)$. 
\\
{\rm (3)} $\omega(x) \setminus \mathop{\mathrm{Cl}}(v) \neq \emptyset$ and there is a point $z \in S$ with $x \in \alpha(z)$. 
\end{corollary}

\section{Reductions of quasi-Q-sets and quasi-circuits}

Let $v$ be a flow on a compact connected surface $S$. 

\subsection{Reductions of quasi-Q-sets into Q-set under countability of singular points}

Recall that every Q-set is a quasi-Q-set because of Corollary~\ref{preQ}. 
By Lemma~\ref{lem:q-Q_uncountable}, we show that a quasi-Q-set is a generalization of a Q-set.  

\begin{proposition}\label{q-Q}
Every quasi-Q-set of a flow with countably many singular points on a compact connected surface is a Q-set. 
\end{proposition}

\subsection{Reductions of quasi-circuits into circuits under total disconnectivity of singular points}

For a closed disk $W$ and disjoint transverse arcs $\mu',\mu'' \subset \partial W$, an orbit arc $I$ in $W$ has the orbit direction from $\mu'$ to $\mu''$ if $\vert I \cap \partial W \vert = 2$ and $I$ is an orbit arc from a point in $\mu'$ to a point in $\mu''$. 
We have the following observations. 

\begin{lemma}\label{lem:connecting_boundary}
Let $v$ be a flow with totally disconnected singular points on a surface $S$ and $x \in S$ a point whose $\omega$-limit set is a limit quasi-circuit. 
For any closed disk $W$ whose interior intersects $\omega(x)$ and whose boundary is transverse to $v$ except finitely many tangencies $p_1, \ldots , p_k \in S$ such that $x \notin W$ and $\omega(x) \setminus W \neq \emptyset$ and for any point $y \in \omega(x) \cap \mathop{\mathrm{int}}W$, there is a sequence $(I_n)_{n \in \Z_{\geq 0}}$ of connected components $I_n$ of $W \cap O^+(x)$ satisfying the following properties: 
\\
{\rm (1)} There are connected components $\mu', \mu''$ of the complement $\partial W - \{ p_1, \ldots , p_k \}$ of the tangencies of the loop  $\partial W$ such that the connected component $I_n$ are orbit arcs which have the orbit direction from $\mu'$ to $\mu''$ for any $n \in \Z_{\geq 0}$. 
\\
{\rm (2)} The set difference $\overline{\bigcup_{n=1}^\infty I_n}^W - \bigcup_{n=1}^\infty I_n \subset \omega(x)$ contains $y$, is connected, and is not a singleton, where $\overline{A}^W$ is the closure of a subset $A \subseteq W$ with respect to $W$. 
\\
{\rm (3)} 
Any connected components of $(\bigcup_{n=1}^\infty V_n - V_0) \cap O^+(x)$ are orbit arcs which have the orbit direction from $\mu'$ to $\mu''$, where $V_n$ is the connected component of $W - I_n$ not containing $y$ for any $n \in \Z_{\geq 0}$. 
\end{lemma}

\begin{proof}
By Theorem~\ref{main:a}, any non-recurrent orbit in the $\omega$-limit set of a point is a connecting quasi-separatrix. 
Therefore, the limit quasi-circuit $\omega(x)$ does not intersect $O(x)$ and so $\omega(x) \subseteq \overline{O^+(x)} \setminus O(x)$. 

\begin{claim}\label{claim:eq_01}
$\omega(x) = \overline{O^+(x)} - (O^+(x) \sqcup \{ x\})$. 
\end{claim}
\begin{proof}
For any $t \in \R_{>0}$, from $\overline{v(\R_{>0},x)} = \overline{v(\R_{>t},x)} \cup v([0,t],x)$, we have $\overline{O^+(x)} - (O^+(x) \sqcup \{ x\}) = \overline{v(\R_{>0},x)} - (O^+(x) \sqcup \{ x\}) = \overline{v(\R_{>t},x)} \setminus (O^+(x) \sqcup \{ x\})$.
By $\omega(x) \cap O(x) = \emptyset$, we obtain $\omega(x) = \omega(x) \setminus (O^+(x) \sqcup \{ x\}) = \bigcap_{t \in \R} \overline{v(\R_{>t},x)} \setminus (O^+(x) \sqcup \{ x\}) = \bigcap_{t >0} \overline{v(\R_{>t},x)} \setminus (O^+(x) \sqcup \{ x\}) = \overline{O^+(x)} - (O^+(x) \sqcup \{ x\})$. 
\end{proof}

By $\omega(x) \cap \mathop{\mathrm{int}}W \neq \emptyset$ and $\omega(x) \setminus W \neq \emptyset$, the intersectoin $O^+(x) \cap W$ contains a connected component $I$ of the intersection $W \cap O^+(x)$ which is an orbit arc with $\vert I \cap \partial W \vert = 2$. 
This implies that there are tangencies of the loop $\partial W$. 

Fix any point $y$ in $\omega(x) \cap \mathop{\mathrm{int}}W$. 
Let $p_1, \ldots , p_k \in S$ be the tangencies of the loop $\partial W$ for some $k \in \Z_{\geq 2}$. 
Denote by $\mu_1, \ldots , \mu_k$ the connected components of the complement $\partial W - \{ p_1, \ldots , p_k \}$ of the tangencies of $\partial W$.
Then $\mu_1, \ldots , \mu_k$ are transverse open arcs. 
For any point $x' \in O^+(x) \cap W$, denote by $I_{x'}$ the connected component of $O^+(x) \cap W$ containing $x'$. 
By $\omega(x) \cap \mathop{\mathrm{int}}W \neq \emptyset$, $\omega(x) \setminus W \neq \emptyset$, and $x \notin W$, any connected components of $O^+(x) \cap W$ are closed arcs between $\partial W$. 
From $y \in \omega(x) \cap \mathop{\mathrm{int}}W \subset W - O(x)$, there is a sequence $(x_{n})_{n \in \Z_{\geq 0}}$ of points $x_{n} \in O^+(x)$ converging to $y$ with $I_{x_n} \neq I_{x_m}$ for any integers $n\neq m \in \Z_{\geq 0}$. 
By the finiteness of tangencies $p_1, \ldots , p_k \in \partial W$ of the loop $\partial W$, taking a subsequence of of $x_{n}$, we may assume that $\vert I_n \cap \partial W \vert = 2$ for any $n  \in \Z_{\geq 0}$, and that there are numbers $i_{-}, i_{+} \in \{ 1, \ldots , k \}$ such that the connected component $I_n$ are orbit arcs which have the orbit direction from $\mu_{i_{-}}$ to $\mu_{i_{+}}$ for any $n  \in \Z_{\geq 0}$. 
Denote by $U_n$ (resp. $V_n$) the connected component of $W - I_n$ containing (resp. not containing) $y$. 
Then $U_{n+1} \subset U_n$ and $V_n \subset V_{n+1}$.  
Put $V_\infty := \bigcup_{n=1}^\infty V_n$ and $I_\infty := \overline{\bigcup_{n=1}^\infty I_n}^W - \bigcup_{n=1}^\infty I_n$, where $\overline{A}^W$ is the closure of a subset $A \subseteq W$ with respect to $W$. 

\begin{claim}
$y \in I_\infty \subseteq \omega(x)$. 
\end{claim}

\begin{proof}
Since $x$ is not positive recurrent, by $\bigcup_{n=1}^\infty I_n \subset O^+(x)$, we obtain that $\overline{\bigcup_{n=1}^\infty I_n}^W \cap (O^+(x) \sqcup \{ x\}) =\bigcup_{n=1}^\infty I_n$. 
By definition of $I_\infty$, the closedness of $W$ implies that $I_\infty = \overline{\bigcup_{n=1}^\infty I_n}^W  - \bigcup_{n=1}^\infty I_n = \overline{\bigcup_{n=1}^\infty I_n}^W \setminus (O^+(x) \sqcup \{ x\}) \subseteq \overline{O^+(x) \cap W} \setminus (O^+(x) \sqcup \{ x\}) \subseteq \overline{O^+(x)} - (O^+(x) \sqcup \{ x\}) = \omega(x)$, because of Claim~\ref{claim:eq_01}. 
Then $y \in \partial_W V_\infty = \partial_W (\bigcup_{n=1}^\infty V_n) = \overline{\bigcup_{n=1}^\infty V_n}^W - \bigcup_{n=1}^\infty V_n = \overline{\bigcup_{n=1}^\infty I_n}^W - \bigcup_{n=1}^\infty I_n  = I_\infty \subset \omega(x)$, where $\partial_W V_\infty$ is the boundary of $V_\infty \subset W$ with respect to $W$. 
\end{proof}

\begin{claim}
The set difference $I_\infty$ is connected and is not a singleton. 
\end{claim}
\begin{proof}
Considering the doubling $S_W$ of $W$ which is a sphere, the doubling $\widetilde{I_n}$ of $I_n$ in $S_W$ is a loop in $S_W$ and the doubling $\widetilde{V_\infty}$ of $V_\infty$ is an open disk in the sphere $S_W$.
By construction, the doubling $\widetilde{I_\infty}$ of $I_\infty$ is the boundary of the open disk $\widetilde{V_\infty}$. 
From \cite[Lemma~4]{sannami2021topological}, the boundary $\partial \widetilde{V_\infty} = \widetilde{I_\infty}$ is connected. 

Assume that $I_\infty$ is not connected. 
Then there are disjoint nonempty open subsets $U_y$ and $V_y$ with $I_\infty \subset U_y \sqcup V_y$. 
Then the doublings $\widetilde{U_y}$ and $\widetilde{V_y}$ of $U_y$ and $V_y$, respectively, are disjoint nonempty open subsets and form an open covering of $\widetilde{I_\infty}$ in $S_W$. 
This means that the doubling $\widetilde{I_\infty}$ is not connected, which contradicts the connectivity of $\widetilde{I_\infty}$. 
Thus $I_\infty$ is connected. 

Since $I_\infty$ intersects $\partial W$, by $y \in I_\infty \cap \mathop{\mathrm{int}}W$, the connectivity of $I_\infty$ implies that the closed subset $I_{\infty}$ is not a singleton. 
\end{proof}


\begin{claim}
Assertion {\rm(3)} holds.
\end{claim}
\begin{proof}
For any $n  \in \Z_{\geq 0}$, since $I_n = \partial_W (V_{n+1} - V_n) - I_{n+1}$ is a  closed orbit arc from a point in $\mu_{i_{-}}$ to a point in $\mu_{i_{+}}$, any connected components of $(V_{n+1} - V_n) \cap O^+(x)$ are orbit arcs from points in $\mu_{i_{-}}$ to points in $\mu_{i_{+}}$.
Therefore any connected components of $(V_\infty - V_0) \cap O^+(x)$ are orbit arcs from points in $\mu_{i_{-}}$ to points in $\mu_{i_{+}}$.
\end{proof}
This completes the proof. 
%
%
\end{proof}

\begin{lemma}\label{lem:non-q-lim-circuit_-1}
Let $v$ be a flow with totally disconnected singular points on a compact connected surface $S$. 
Then every limit quasi-circuit is locally connected. 
\end{lemma}

\begin{proof}
Fix any Riemannian metric on $S$ which induces the Riemannian distance.
Let $\gamma$ be a limit quasi-circuit. 
Since a limit quasi-circuit contains non-recurrent points, so does the limit quasi-circuit $\gamma$. 
By Theorem~\ref{main:a}, the total disconnectivity of $\Sv$ implies that any non-recurrent orbits $O(x)$ in $\gamma$ are connecting separatrices. 
Therefore, the limit quasi-circuit $\omega(x)$ does not intersect $O(x)$ and so $\omega(x) \subseteq \overline{O^+(x)} \setminus O(x)$. 
The invariance of $\gamma$ implies that each orbit closure in $\gamma$ is either a singular point, a closed interval, or a loop. 
By time reversion if necessary, we may assume that $\omega(x) =\gamma$ for a point $x$. 
Then the difference $\gamma \setminus \mathop{\mathrm{Sing}}(v)$ is a disjoint union of open intervals.

Assume that $\gamma$ is not locally connected.
There is a point $x_{\infty} \in \gamma$ at which $\gamma$ is not locally connected. 
Fix a small closed disk $W'$ which is a \nbd of $x_{\infty}$ with $\gamma \setminus W' \neq \emptyset$ and $x \notin W'$. 
By the total disconnectivity of $\Sv$, the complement $S - \Sv$ is a connected open surface.
Therefore, there are a loop $\mu \subset W' - \Sv = W' \cap (S - \Sv)$ and a closed disk $W \subseteq W'$ whose boundary is $\mu$ such that $x_{\infty} \in \mathop{\mathrm{int}}W$. 
Put $\gamma_W := \gamma \cap W$. 
From \cite[Lemma~3.1]{Kibkalo2022topological}, by modifying $W$, we may assume that $\mu$ is transverse to $v$ except for finitely many tangencies $p_1, \ldots , p_k \in \partial W$  of the loop $\partial W$. 
Moreover, by modifying $W$, we may assume that the length of $\mu$ is finite. 

\begin{claim}
There are tangencies of the loop $\mu$.
\end{claim}
\begin{proof}
Because $x_{\infty} \in \omega(x) = \gamma$, the intersection $O^+(x) \cap W$ contains pairwise disjoint infinitely many orbit arcs. 
Then there is a connected component of $O^+(x) \cap W$ which transversely intersects the loop $\mu = \partial W$ once in the inward direction and once in the outward direction of the orbit $O(x')$. 
Therefore, there are tangencies of the loop $\mu$. 
\end{proof}

By the previous claim, denote by $\mu_1, \ldots , \mu_k$ the connected components of the complement in $\partial W$ of the tangencies of the loop $\mu = \partial W$ for some $k \in \Z_{\geq 2}$. 
For any point $y \in \omega(x) \cap \mathop{\mathrm{int}}W$, applying Lemma~\ref{lem:connecting_boundary} to $y$, there is a sequence $(I_{y,m})_{m \in \Z_{\geq 0}}$ of the connected components $I_{y,m}$ of $W \cap O^+(x)$ with $I_{y,\infty}= \overline{\bigcup_{m=1}^\infty I_{y,m}}^W - \bigcup_{m=1}^\infty I_{y,m}$ such that $I_{y,\infty}$ contains $y$, is connected, and is not a singleton, and that any connected components $I_{y,m}$ have the orbit direction from $\mu_{i_{y,-}}$ to $\mu_{i_{y,+}}$ for some $i_{y,-}, i_{y,+} \in \{1, \ldots , k\}$. 

\begin{claim}
The subset $\gamma_W = \gamma \cap W$ has infinitely many connected components. 
\end{claim}

\begin{proof}
Assume that $\gamma_W$ has at most finitely many connected components. 
Denote by $J_{x_\infty}$ the connected component of $\gamma_W  = \gamma \cap W$ containing $x_\infty$. 
Since any connected components are closed, the closedness of $W$, the complement $\gamma_W - J_{x_\infty}$ is a finite disjoint union of connected components of $\gamma_W$ and is a closed subset of $S$. 
This means that $J_{x_\infty}$ is a connected open \nbd of $x_\infty$ with respect to the subspace $\gamma_W$.
Since $W$ is a \nbd of $x_\infty$ in $S$, the subset $_{x_\infty} \subset W$ is a connected \nbd of $x_\infty$ with respect to the subspace $\gamma$, which contradicts the absence of the local connectivity of $\gamma$ at $x_\infty$. 
\end{proof}

\begin{claim}
The subset $\gamma_W = \gamma \cap W$ has infinitely many connected components which contain no tangencies $p_1, \ldots , p_k$. 
\end{claim}

\begin{proof}
By the invariance of $\gamma = \omega(x)$, any connected component of $\gamma_W = \gamma \cap W = \omega(x) \cap W$ intersecting a point of the boundary $\partial W - \{ p_1, \ldots , p_k \}$ intersects $\mathop{\mathrm{int}}W$. 
By the finiteness of tangencies $\{ p_1, \ldots , p_k \}$, the number of the connected components of $\gamma_W$ which do not intersect $\mathop{\mathrm{int}}W$ is finite. 
Therefore, the assertion holds, because of the previous claim. 
\end{proof}

\begin{claim}
The length of $\mu$ is infinite.
\end{claim}

\begin{proof}
By the previous claim, let $(J_{n\infty})_{n \in \Z_{\geq 0}}$ be a sequence of pairwise disjoint connected components of $\gamma_W$ which contains no tangencies $p_1, \ldots , p_k$. 
Moreover, by Lemma~\ref{lem:connecting_boundary}, for any $n \in \Z_{\geq 0}$, there is a sequence $(I_{n,m})_{m \in \Z_{\geq 0}}$ of the connected components $I_{n,m}$ of $W \cap O^+(x)$ with $I_{n,\infty} := \overline{\bigcup_{m=1}^\infty I_{n,m}}^W - \bigcup_{m=1}^\infty I_{n,m} \subseteq J_{n,\infty}$ such that $I_{n,\infty}$ intersects $\omega(x)$, is connected, and is not a singleton, and that any connected components $I_{n,m}$ have the orbit direction from $\mu_{i_{n,-}}$ to $\mu_{i_{n,+}}$ for some $i_{n,-}, i_{n,+} \in \{1, \ldots , k\}$. 
Taking a subsequence of the sequence $(I_{n\infty})_{n \in \Z_{\geq 0}}$, we may assume that there are numbers $i_{\infty,-}, i_{\infty,+} \in \{ 1, \ldots , k\}$ such that each $I_{nm}$ has the orbit direction from $\mu_{i_{\infty,-}}$ to $\mu_{i_{\infty,+}}$ for any $n \in \Z_{\geq 1}$ and any $m \in \Z_{\geq 0}$. 
For any $n \in \Z_{\geq 1}$ and any $m \in \Z_{\geq 0}$, by $\vert I_{nm} \cap \mu_{i_{\infty,-}}\vert = 1$ and $\vert I_{nm} \cap \mu_{i_{\infty,+}}\vert = 1$, denote by $x_{nm,-}$ (resp. $x_{nm,+}$) the point in the singleton $I_{nm} \cap \mu_{i_{\infty,-}}$ (resp. $I_{nm} \cap \mu_{i_{\infty,+}}$). 
Taking a subsequence of the sequence $(I_{n\infty})_{n \in \Z_{\geq 0}}$, we may assume that the subsequence $(x_{nn,-})_{ \in \Z_{\geq 1}}$ (resp. $(x_{nn,+})_{ \in \Z_{\geq 1}}$ ) is monotonic in the open interval $\mu_{i_{\infty,-}}$ (resp. $\mu_{i_{\infty,+}}$). 

By the existence of a collar basin $\A$ of $\omega(x) = \gamma$, for any $n \in \Z_{\geq 1}$, there is a positive number $d_{\A}$ such that the lengths of the arcs in $\mu_{i_{\infty,-}} \subset \partial W$ connecting $I_{n\infty} \cap \mu_{i_{\infty,-}}$ and $I_{n+1 \infty} \cap \mu_{i_{\infty,-}}$ is more than $d_{\A}$ as in the proof of Lemma~\ref{lem:acc_lc}. 
The infinity of $(I_{n\infty})_{n \in \Z_{\geq 1}}$ implies that the length of $\partial W = \mu$ is infinite.
\end{proof}

The previous claim contradicts the finiteness of the length of $\mu$. 
\end{proof}

Lemma~\ref{lem:qc_ch} and the previous lemma imply the following observation.

\begin{lemma}\label{lem:non-q-lim-circuit}
Let $v$ be a flow with totally disconnected singular points on a compact connected surface $S$. 
Then a limit quasi-circuit is the image of a circle.
\end{lemma}

\subsubsection{Case of totally disconnected singular point set}

By Lemma~\ref{lem:q-Q_uncountable} and Lemma~\ref{lem:non-q-lim-circuit}, Theorem~\ref{main:a} can be reduced into the following statement, which is a refinement of  \cite[Theorem~3.1]{marzougui1996} (cf. \cite[Theorem~2.1]{marzougui2000flows})), if the singular point set is totally disconnected. 

\begin{corollary}\label{cor:PB_totally_disconn}
The following statements hold for a flow with totally disconnected singular point set on a compact surface: 
\\
{\rm(a)} The $\omega$-limit set of any non-closed orbit is one of the following exclusively:
\begin{quote}
\setlength{\leftskip}{-25pt}$(1)$ A singular point.
\\
$(2)$ A semi-attracting limit cycle.
\\
$(3)$ A quasi-semi-attracting limit quasi-circuit that is the image of a circle.
\\
$(4)$ A locally dense Q-set. 
\\
$(5)$ A transversely Cantor Q-set. 
\\
$(6)$ A quasi-Q-set that consists of singular points and non-recurrent points.  
\end{quote}
{\rm(b)} Every non-recurrent orbit in the $\omega$-limit set of a point is a connecting separatrix. 
\end{corollary}

The countability of singular points implies the fololwing statement. 

\begin{theorem}\label{cor:PB_countable}
For a flow with countably many singular point set on a compact surface, the $\omega$-limit set of any non-closed orbit is one of the following exclusively:
\begin{quote}
\setlength{\leftskip}{-25pt}$(1)$ A singular point.
\\
$(2)$ A semi-attracting limit cycle.
\\
$(3)$ A quasi-semi-attracting limit quasi-circuit that is the image of a circle and consists of singular points and connecting separatrices.
\\
$(4)$ A locally dense Q-set. 
\\
$(5)$ A transversely Cantor Q-set. 
\end{quote}
\end{theorem}

\subsubsection{Reduction under finiteness of singular points}

We show that a limit quasi-circuit is a generalization of a limit circuit.  

\begin{proposition}\label{non-q-lim-circuit}
Let $v$ be a flow with finitely many singular points on a compact connected surface $S$. 
Then a limit quasi-circuit is a semi-attracting or semi-repelling limit non-periodic circuit, which is a continuous image of a circle. 
\end{proposition}

\begin{proof}
Let $\gamma$ be a limit quasi-circuit. 
By Corollary~\ref{cor:PB_totally_disconn} and its dual statement, the limit quasi-circuit  $\gamma$ is a continuous image of a circle that consists of singular points and connecting separatrices. 
The finiteness of $\Sv$ implies that $\gamma$ is a non-periodic circuit.

%

We claim that $\gamma$ is semi-attracting.
Indeed, by definition of limit quasi-circuit, there is a small closed collar $\A$ that is either positive invariant or negative invariant such that $\gamma$ is a boundary component of its collar $\A$ and that $\partial \A - \gamma$ is a loop consisting of a closed orbit arc and a transverse closed interval. 
By time reversion if necessary, we may assume that $\omega(x) =\gamma$ for a point $x$. 
Then $\A$ is positive invariant. 
Since $\mathop{\mathrm{Sing}}(v)$ is finite, taking $\A$ small if necessary, we may assume that $\A$ contains no singular points. 
Because a limit quasi-circuit contains non-recurrent points, so does the limit quasi-circuit $\gamma$. 
Since $\A$ contains no singular point, by $\omega(x) =\gamma$, there is a transverse closed arc $T \subset \gamma \sqcup \A$ whose boundary intersects $\gamma$ such that the first return map on the interior $\mathop{\mathrm{int}}T$ is attracting. 
If $\A$ contains a periodic orbit $O$, then $O$ bounds an invariant closed disk $D \subset \A$ which contains a singular point because of Poincar\'{e}-Hopf theorem to the restriction $v\vert_D$, which contradicts the non-existence of singular points in $\A$. 
Thus $\A$ contains no periodic points and so $\A \subset \mathrm{P}(v)$.  
Taking $\A$ small, we may assume that the boundary component $\partial \A - \gamma$ consists of one orbit arc in $O(x)$ and one sub-arc in $T$ such that any connected components $B_i$ of $\A - (T \cup O(x))$ are flow boxes in $\mathrm{P}(v)$. 
Since any flow boxes $B_i$ can be considered as the restriction of a flow on a sphere, by Corollary~\ref{cor:nonex_qQset}, any flow boxes $B_i$ are trivial and so $\A \subset v(T)$ such that the domain of the first return map to $\mathop{\mathrm{int}}T$ is $\mathop{\mathrm{int}}T$. 
Because $\omega(x) = \gamma$ is a boundary component of $\partial \A$ which is a limit circuit, 
we have $\gamma = \omega(x) = \omega(y)$ for any $y \in \A$.
\end{proof}

\subsection{Reduction of the Poincar\'{e}-Bendixson theorem for flows with finitely many singular points}

The reductions of quasi-Q-sets (Proposition~\ref{q-Q}) and of quasi-circuits (Proposition~\ref{non-q-lim-circuit}) imply a proof of the following generalization of the Poincar\'{e}-Bendixson theorem for a flow with finitely many singular points (see for example \cite{nikolaev1999flows}). 

\begin{corollary}\label{PBthm}
Let $v$ be a flow with finitely many fixed points on a compact surface $S$.
Then the $\omega$-limit set of any non-closed orbit is one of the following exclusively:
\\
$(1)$ A singular point.
\\
$(2)$ A semi-attracting limit cycle.
\\
$(3)$ A semi-attracting limit non-periodic circuit.
\\
$(4)$ A locally dense Q-set.
\\
$(5)$ A transversely Cantor Q-set.
\end{corollary}

\subsubsection{Poincar\'{e}-Bendixson theorem for a flow with finitely many singular points on possibly non-compact surfaces}
Recall that a non-recurrent orbit on $S$ is a virtual separatrix if it is a connecting separatrix on $S_{\mathrm{end}}$ with respect to $v_{\mathrm{end}}$. 
An invariant subset on $S$ is a {\bf semi-attracting limit virtual circuit} if it is the resulting subset from a semi-attracting limit circuit on $S_{\mathrm{end}}$ with respect to $v_{\mathrm{end}}$ by removing all the ends. 

The previous corollary implies the following generalization of Poincar\'{e}-Bendixson theorem for a flow with finitely many singular points on possibly non-compact surfaces of finite genus and finitely many boundary components.

\begin{corollary}\label{PBthm_non_cpt}
The following statements hold for a flow with finitely many singular points on a surface of finite genus and finitely many boundary components: 
\\
{\rm(a)} The $\omega$-limit set of any non-closed orbit is one of the following exclusively:
\begin{quote}
\setlength{\leftskip}{-25pt}
$(0)$ The empty set.
\\
$(1)$ A singular point.
\\
$(2)$ A semi-attracting limit cycle.
\\
$(3)$ A semi-attracting limit non-periodic virtual circuit. 
\\
$(4)$ A locally dense Q-set. 
\\
$(5)$ A transversely Cantor Q-set. 
\end{quote}
{\rm(b)} Any non-recurrent orbit in the $\omega$-limit set of a point is a virtual separatrix. 
\\
{\rm(c)} If the $\omega$-limit set of a point is a Q-set, then the Q-set corresponds to the orbit closure of any non-closed recurrent point in the Q-set. 
\end{corollary}

\section{Construction of flow boxes with non-arcwise-connected invariant subsets}

In this section, we introduce an operation that makes $\omega$-limit sets not arcwise-connected by constructing flow boxes with non-arcwise-connected invariant subsets. 
To state the operation, we have the following statement. 

\begin{lemma}\label{lem:deformation_01}
For any flow $v$ on a surface $S$ of finite genus and finitely many boundary components with a non-singular point $x$, there is a trivial flow box $B$ containing $x$ such that the resulting flow $w$ by replacing $B$ with a flow box satisfies the following properties: 
\\
{\rm(1)} The restriction $v\vert_{S - O(x)}$ is topologically equivalent to the restriction $w\vert_{S - O(x)}$. 
\\
{\rm(2)} For any point $y \in S - O(x)$, we have that $\alpha(y) = \alpha_w(y)$ and $\omega(y) = \omega_w(y)$. 
\\
{\rm(3)} 
If $O(x)$ is periodic, then $O_w(x)$ is the disjoint union of one non-recurrent orbit $O_0$ and one singular point $x$ with $\omega_w(O_0) = x = \alpha_w(O_0)$. 
\\
{\rm(4)} If $O(x)$ is not periodic, then $O_w(x)$ is the disjoint union of two non-singular orbits $O_1, O_2$ and one singular point $x$ with $\omega_w(O_1) = x = \alpha_w(O_2)$ such that $\alpha(x) = \alpha_w(O_1)$ and $\omega(x) = \omega_w(O_2)$. 
\end{lemma}

We call $w$ in the previous lemma the {\bf resulting flow of $\bm{v}$ by replacing a non-singular point $\bm{x}$ with a singular point}, and denote by $\bm{v_x}$ the resulting flow $w$. 
Roughly speaking, the resulting flow by replacing a non-singular point with a singular point is a flow obtained by replacing a trivial flow box with a flow box as in Figure~\ref{figure:s_flow_box}. 

\begin{proof}[Proof of Lemma~\ref{lem:deformation_01}]
Since $S$ can be identified with a subset of a compact surface by the end completion of $S$, the flow $v$ can be identified with the restriction of the resulting flow of $v$ considering the ends as singular points. 
By Gutierrez's smoothing theorem~\cite{gutierrez1986smoothing}, we may assume that the flow $v$ is a $C^1$-flow generated by an integrable continuous vector field $X$ on $S$. 
Fix any open trivial flow box $B$ intersecting $O(x)$. 
Identifying $B$ with the square $[-1,1]^2$ such that $\{0 \} \times [-1,1]$ is an orbit arc of $O(x)$ and that $v\vert_B$ is generated by a vector field $X = (1,0)$. 
Take a $C^\infty$ bump function $\varphi \colon B = [-1,1]^2 \to [0,1]$ with $\varphi^{-1}(0) = \{ 0 \}$ such that $\varphi$ is one near the boundary $\partial [-1,1]^2$. 
Consider the flow box $B'$ on $B$ whose orbits arc are generated by the vector field $\varphi X = (\varphi, 0)$ as in Figure~\ref{figure:s_flow_box}. 
\begin{figure}
\begin{center}
\includegraphics[scale=0.5]{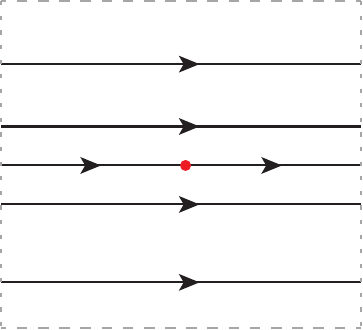}
\end{center}
\caption{A flow box with one singular point with respect to the vector field $\varphi X = (\varphi, 0)$.}
\label{figure:s_flow_box}
\end{figure}
Replacing the trivial flow box $B$ with the flow box $B'$, the resulting flow is desired. 
\end{proof}

\subsection{Resulting surface by collapsing a non-separating continuum into a singleton}

\subsubsection{Decompositions and upper semicontinuity}

By a {\bf decomposition}, we mean a family $\mathcal{F}$ of pairwise disjoint nonempty subsets of a set $X$ such that $X = \bigsqcup \mathcal{F}$, where $\bigsqcup$ denotes a disjoint union.
For a decomposition $\mathcal{F}$ on a set $X$, a subset $A \subseteq X$ is {\bf saturated} if $A$ is a union of elements of $\mathcal{F}$. 
or a decomposition $\mathcal{F}$ on a set $X$, the union of elements of $\mathcal{F}$ intersecting a subset $A \subseteq X$ is called the {\bf saturation} of $A$. 
We set $X/\mathcal{F} := X/\sim_{\mathcal{F}}$, where $p \sim_{\mathcal{F}} q$ if there is an element of $\mathcal{F}$ containing $p$ and $q$. 
A decomposition $\mathcal{F}$ of closed and compact elements on a topological space $X$ if {\bf upper semicontinuous} if for any element $L \in \mathcal{F}$ and for any open \nbd $U \subseteq X$ of $L$ there is a saturated \nbd of $L$ in $U$. 
We recall the following characterization of upper semicontinuty for a decomposition of connected compact elements of a locally compact Hausdorff space.

\begin{lemma}[Remark after Theorem~4.1\cite{Epstein1976}]\label{lem:eps}
The following statements are equivalent for a decomposition $\mathcal{F}$ of connected compact elements of a locally compact Hausdorff space $X$:
\\
{\rm(1)} The decomposition $\mathcal{F}$ is upper semicontinuous.
\\
{\rm(2)} The quotient space $X/\mathcal{F}$ is Hausdorff. 
\\
{\rm(3)} The canonical projection $p \colon X \to X/\mathcal{F}$ is closed {\rm(i.e.} the saturations of any closed subsets are closed{\rm)}.
\end{lemma}

\subsubsection{Continua and non-separating sets}

By a {\bf continuum}, we mean a nonempty compact connected metrizable space. 
A subset $C$ in a topological space $X$ is {\bf separating} if the complement $X - C$ is disconnected. 
A subset in a topological space is {\bf non-separating} if it is not separating.  

\subsubsection{Resulting surface by collapsing a continuum into a singleton}

We recall the following statement. 
 
 \begin{lemma}[Moore's theorem (cf. p.3\cite{daverman1986decompositions})]\label{lem:Moore}
 For an upper semicontinuous decomposition $\mathcal{F}$ into non-separating continua on a surface $S$ which is either a plane or a sphere, the quotient space $S/\mathcal{F}$ is homeomorphic to $S$ unless $\mathcal{F}$ is the singleton of the surface.
 \end{lemma}

We have the following tool. 

\begin{corollary}
Let $S$ be a surface and $C \subseteq S$ a non-separating continuum which is contained in an open disk in $S$. 
The quotient space $S/\mathcal{F}_C$ is homeomorphic to $S$, where $\mathcal{F}_C$ is a decomposition $\{ \{ x \} \mid x \in S - C \} \sqcup \{ C \}$.
\end{corollary}

Then the resulting surface $S/\mathcal{F}_C$ is called the {\bf resulting surface from $S$ by collapsing $C$ into a singleton}.

\begin{proof}
By definition of $\mathcal{F}_C$, since singletons are non-separating continua, the decomposition $\mathcal{F}_C$ consists of closed non-separating continua. 
Let $D$ be an open disk in $S$ containing $C$. 
Then the restriction $\mathcal{F}_C|_D$ is a decomposition on $D$. 
Since the saturation of any closed subset $A \subseteq S$ is either $A$ or $A \cup C$, the saturation of $A$ is closed. 
By Lemma~\ref{lem:eps}, the decomposition $\mathcal{F}_C$ and so the restriction $\mathcal{F}_C\vert_D$ to $D$ is upper semicontinuous. 
Moore's theorem (i.e. Lemma~\ref{lem:Moore}) implies that the quotient space $D/\mathcal{F}_C\vert_D$ is homeomorphic to $D$ and so that $S/\mathcal{F}_C$ is homeomorphic to $S$. 
\end{proof}

%
Recall that a flow $v \colon \R \times Z \to Z$ is {\bf topologically semi-conjugate} to a flow $w \colon \R \times Y \to Y$ via $h \colon Y \to Z$ if $h$ is a continuous surjection such that $v(t, h(y)) = h(w(t,y))$ for any $(t,y) \in \R \times Y$. 
In this section, we show the following statement. 

\begin{figure}
\begin{center}
\includegraphics[scale=0.35]{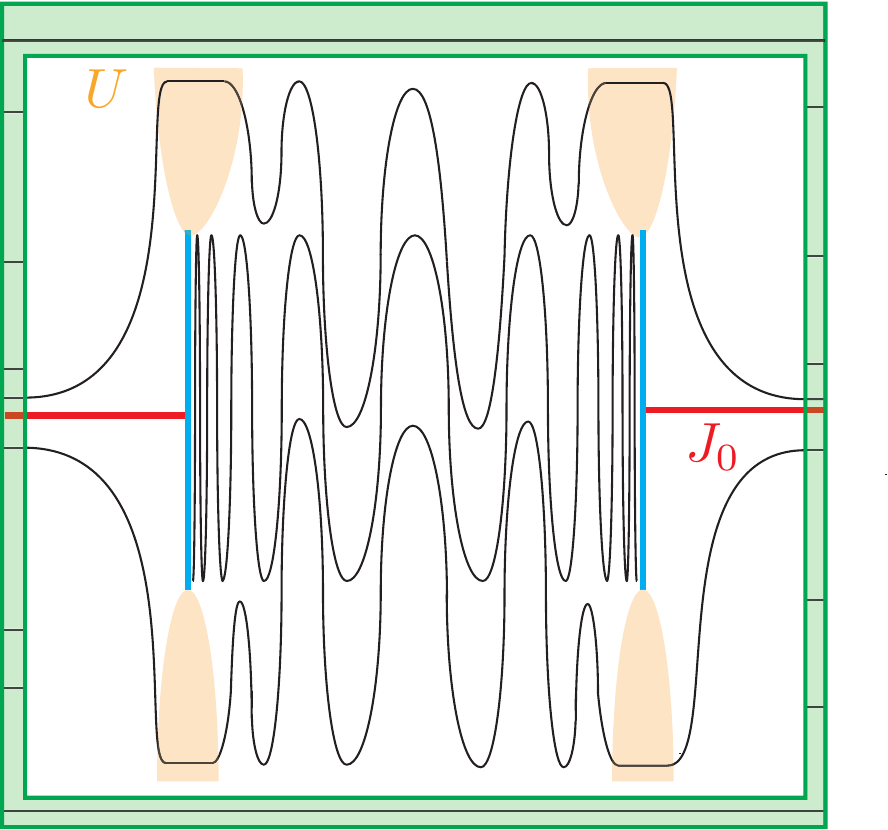}
\end{center}
\caption{The flow box $B$ with a non-arcwise-connected invariant subset with respect to the vector field $X$.}
\label{fig:vf_X}
\end{figure}

\begin{theorem}\label{thm:deformation_02}
Let $v$ be a flow on a surface $S$ with an $\omega$-limit set $\omega$ of a point containing non-singular point  $p_0$ and with a point $q_0 \in S - \omega$ satisfying $\omega(q_0) = \omega$. 
Then there is a trivial flow box $B_{p_0}$ containing $p_0$ such that the resulting flow $w$ by replacing $B_{p_0}$ with a flow box satisfies the following properties: 
\\
{\rm(1)} The $\omega$-limit set $\omega_w(q_0)$ is not arcwise-connected. 
\\
{\rm(2)} The restriction $v\vert_{S - \omega}$ to the complement $S - \omega$ is topologically equivalent to the restriction $w\vert_{S - \omega_w(q_0)}$. 
\\
{\rm(3)} The flow $v_{p_0}$ is topologically equivalent to a flow $v'$ which is topologically semi-conjugate to $w$, where $v_{p_0}$ is the resulting flow of $v$ by replacing $p_0$ with a singular point. 
\\
{\rm(4)} The topological semi-conjugacy from $w$ to $v'$ can be obtained by collapsing a closed invariant subset of $\omega_w(q_0)$ into a singleton. 
\end{theorem}
Roughly speaking, the resulting flow $w$ in the previous theorem can be obtained by replacing a trivial flow box with a flow box as in Figure~\ref{fig:vf_X}. 
This theorem implies Theorem~\ref{main:deformation_02}. 

\subsection{Construction of a flow box with a non-arcwise-connected connected invariant subset}\label{ex00}

Let $\mathbb{D} := [0,1] \times [-1,1]$ a closed square and $\partial_{\pitchfork} \mathbb{D} := \{0, 1\} \times [-1,1]$. 
Define a diffeomorphism $h_0 \colon (0,1) \to \R$ by $h_0(x) : = \tan \left( \dfrac{\pi(x-1)}{2} \right)$. 
Define a function $f \colon (0,1) \to [-1,1]$ as follows: 
\[
f (x) := \cos (h_0(x)) = \cos \left( \tan \left( \dfrac{\pi(x-1)}{2} \right)\right)
\]
Let $G := \{ (x,f (x) \mid x \in (0,1)\} \subset \mathbb{D}$ be the graph of $f$ and $C := G \sqcup \partial_{\pitchfork} \mathbb{D}$ the union. 
Consider a closed square $B := [0,1] \times [-2,2]$ containing $\mathbb{D}$. 
Put $\partial_{\pitchfork} B := \{0, 1\} \times [-2,2]$. 
Then the set difference $B - \partial B = (0, 1) \times (-2,2)$ is an open square $\mathop{\mathrm{int}} B$.
We have the following observation. 

\begin{lemma}\label{lem:complement_arc}
The set difference $\mathop{\mathrm{int}} B - G$ consists of two open disks. 
\end{lemma}

\begin{proof}
The closure $\overline{G}$ in $\R^2$ is the union $G \sqcup \partial_{\pitchfork} \mathbb{D} = C$. 
Let $\mathcal{F}$ be a decomposition of $\R^2$ by $\mathcal{F} := \{ \{ p \} \mid p \in \R^2 - \partial_{\pitchfork} \mathbb{D} \} \sqcup \{ \{0\} \times [-1,1], \{1\} \times [-1,1] \}$. 
By Moore's theorem (cf. p.3 in \cite{daverman1986decompositions})], the quotient space $\R^2/\mathcal{F}$ of the upper semicontinuous decomposition $\mathcal{F}$ into non-separating continua is homeomorphic to $\R^2$, where $\R^2/\mathcal{F}$ is the quotient space $\R^2/\sim_{\mathcal{F}}$ defined by $p \sim_{\mathcal{F}} q$ if there is an element of $\mathcal{F}$ containing $p$ and $q$. 
Let $p_{\mathcal{F}} \colon \R^2 \to \R^2/\mathcal{F}$ be the quotient map. 
Then the image $p_{\mathcal{F}}(C) = p_{\mathcal{F}}(G) \sqcup \{p_{\mathcal{F}}((0,0)), p_{\mathcal{F}}((1,0))\}$ is homeomorphic to a closed interval.  
The set difference $p_{\mathcal{F}}(B) - p_{\mathcal{F}}(C)$ consists of two disks that are homeomorphic to $[0,1] \times (0,2]$.
Therefore the set difference $p_{\mathcal{F}}(\mathop{\mathrm{int}} B) \setminus p_{\mathcal{F}}(G)$ consists of two open disks that are homeomorphic to $(0,1) \times (0,2)$.
Since the restriction $p_{\mathcal{F}}\vert_{\R^2 - \partial_{\pitchfork} \mathbb{D}}$ is identical, we have that $\mathop{\mathrm{int}} B - G = \mathop{\mathrm{int}} B - p_{\mathcal{F}}(G) = p_{\mathcal{F}}(\mathop{\mathrm{int}} B) \setminus p_{\mathcal{F}}(G)$ and so that 
the set difference $\mathop{\mathrm{int}} B - G$ consists of two open disks.
\end{proof}

\begin{figure}
\begin{center}
\includegraphics[scale=0.4]{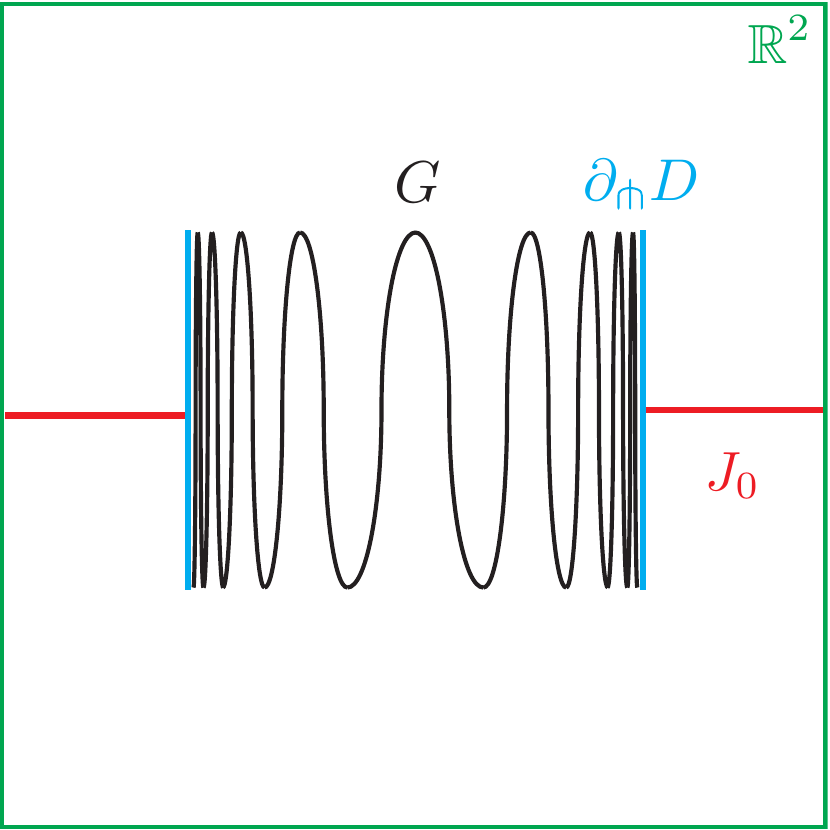}
\end{center}
\caption{An $\omega$-limit set which is a non-locally-connected quasi-circuit, which is neither the image of a circle nor a circuit.}
\label{nonac01}
\end{figure}

Let $\mathbb{S}^2$ be the one point compactification $\R^2 \sqcup \{ \infty \}$ of $\R^2$, $J_0 := (\R - [0,1]) \times \{ 0 \}$ the union of two intervals, and $\gamma := C \sqcup J_0 \sqcup \{ \infty \} = G \sqcup \partial_{\pitchfork} \mathbb{D} \sqcup J_0 \sqcup \{ \infty \}$. 
A subset is {\bf annular} if it is homeomorphic to an annulus. 
A continuum $Y$ contained in a surface $S$ is {\bf annular} if there is an open annular \nbd $A$ of $Y$ such that $A - Y$ consists of two open annuli. 
We have the following observations. 

\begin{lemma}
The union $\gamma$ is an annular continuum. 
\end{lemma}

\begin{proof}
Lemma~\ref{lem:complement_arc} implies that the complement $\mathbb{S}^2 - \gamma =\R^2 - C \sqcup J_0$  is the disjoint union of two open disks. 
Removing two points, we can obtain there is an open annular \nbd $A$ of $\gamma$ such that $A - \gamma$ consists of two open annuli. 
\end{proof}

\begin{lemma}\label{lem:collapsing}
Let $\R^2/\mathcal{F}_C$ be the resulting space collapsing the continuum $C = G \sqcup \partial_{\pitchfork} \mathbb{D}$ into a singleton and $p_{\mathcal{F}_C} \colon \R^2 \to \R^2/\mathcal{F}_C$ be the quotient map. 
Then the resulting space $\R^2/\mathcal{F}_C$ is a plane and the quotient map $p_{\mathcal{F}_C} $ is continuous.
\end{lemma}

\begin{proof}
Let $\mathcal{F}$ be a decomposition of $\R^2$ by $\mathcal{F} := \{ \{ p \} \mid p \in \R^2 - \partial_{\pitchfork} \mathbb{D} \} \sqcup \{ \{0\} \times [-1,1], \{1\} \times [-1,1] \}$ as in the proof of Lemma~\ref{lem:complement_arc}. 
By Moore's theorem (cf. p.3 in \cite{daverman1986decompositions})], the quotient space $\R^2/\mathcal{F}$ is homeomorphic to $\R^2$. 
Let $p_{\mathcal{F}} \colon \R^2 \to \R^2/\mathcal{F}$ be the quotient map. 
Then the image $p_{\mathcal{F}}(C) = p_{\mathcal{F}}(G) \sqcup \{p_{\mathcal{F}}((0,0)), p_{\mathcal{F}}((1,0))\}$ is homeomorphic to a closed interval.  
Let $\mathcal{F}_C$ be a decomposition of $\R^2$ by $\mathcal{F}_C := \{ \{ p \} \mid p \in \R^2 - C \} \sqcup \{ C \}$. 
Then the image $p_{\mathcal{F}}(\mathcal{F}_C)$ is a decomposition on a plane $\R^2/\mathcal{F}$ consisting of the closed interval $p_{\mathcal{F}}(C)$ and points. 
Define the quotient space $\R^2/\sim_{\mathcal{F}_C}$ by $p \sim_{\mathcal{F}_C} q$ if there is an element of $\mathcal{F}_C$ containing $p$ and $q$. 
Let $q \colon \R^2/\mathcal{F} \to \R^2/\mathcal{F}_C$ be the quotient map by collapsing the closed interval $p_{\mathcal{F}}(C)$ into a singleton. 
By Moore's theorem (cf. p.3 in \cite{daverman1986decompositions})], the quotient space $\R^2/\mathcal{F}_C$ is homeomorphic to $\R^2$ and the composition $p_{\mathcal{F}_C} : = q \circ p_{\mathcal{F}} \colon \R^2 \to \R^2/\mathcal{F}_C$  is continuous and is the quotient map. 
\end{proof}

Let $\mathbb{B} := [-1,2] \times [-5,5]$ be a closed square containing $B= [0,1] \times [-2,2]$, $I_- := [-1,0) \times \{ 0 \} \subset J_0$ an interval, and $I_+ := (1,2] \times \{ 0 \} \subset J_0$ an interval. 
Put $I := C \sqcup I_- \sqcup I_+$.
Then we have the following vector field. 

\begin{lemma}\label{lem:no_arc_conn_flowbox}
There is a $C^\infty$ vector field $X$ on $\R^2$ satisfying the following properties:  
\\
{\rm(1)} The square $\mathbb{B}$ is a flow box with respect to $X$. 
\\
{\rm(2)} The restriction $X\vert_{\R^2 - \mathbb{B}}$ is $(1,0)$. 
\\
{\rm(3)} Subsets $G$, $(- \infty,0) \times \{ 0 \}$ and $(1, \infty) \times \{ 0 \}$ of $\R^2$ are orbits of $X$. 
\\
{\rm(4)} The set $\mathop{\mathrm{Sing}}(X)$ of critical points of $X$ is $\partial_{\pitchfork} \mathbb{D}= \{0, 1\} \times [-1,1]$. 
\\
{\rm(5)} Each of the positive and negative orbits of any points in $\R^2 - (G \sqcup \partial_{\pitchfork} \mathbb{D} \sqcup J_0)$ is neither singular nor periodic but is unbounded and closed as subsets. 
\\
{\rm(6)} For any convergence sequence $(y_n)_{n \in \Z_{\geq 0}}$ of non-zero numbers $y_n$ tending to $0$ and for any point $p_g \in G$, there is a sequence $(t_n)_{n \in \Z_{\geq 0}}$ such that the sequence $v_X(t_n, (-1,y_n))_{n \in \Z_{\geq 0}}$ converges to the point $p_g$ in $G$, where $v_X$ is the flow generated by $X$. 
\end{lemma}

\begin{proof}
Define a $C^\infty$ bump function $\varphi \colon \R \to [0,1]$ with $\varphi^{-1}(0) = (- \infty, 1/3]$ and $\varphi^{-1}(1) = [2/3, \infty)$ such that $\varphi$ is increasing on $[1/3, 2/3]$ as in Figure~\ref{fig:bump_increasing}. 
\begin{figure}
\begin{center}
\includegraphics[scale=0.45]{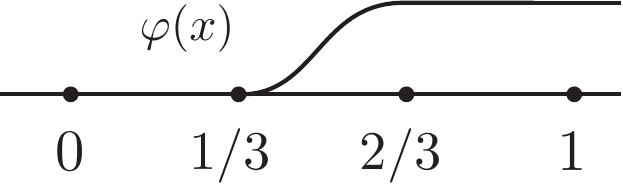}
\end{center}
\caption{A bump function.}
\label{fig:bump_increasing}
\end{figure}
Put $C_{\infty} = C'_{\infty} := \{ (x,\cos x) \mid x \in \R \}$ and $D_0 := \R \times (\R - [-3,3])$. 
Let $f_{0} \colon \R \to \{3 \}$ be a constant map and $g_0 \colon \R \to \{-3 \}$ a constant map and $f_\infty \colon \R \to [-1,1]$ a function defined by $f_\infty (x) := \cos x$. 

For any positive integer $n \in \Z_{>0}$, define a $C^\infty$ even function $f_n \colon \R \to [-1+1/n, 1+1/n]$ as follows:
\[
f_n(x) := 
\begin{cases}
1 + \dfrac{1}{n} & \text{for } x \in \R - [-2\pi n, 2\pi n]\\
\cos x + \dfrac{1}{n} & \text{for } x \in [-2\pi(n-1), 2\pi (n-1)]\\
\varphi \left( \dfrac{\vert x\vert  - \pi(2n -1)}{\pi} \right) \left(  1 - \cos x \right) + \cos x  + \dfrac{1}{n}  \hspace{-150pt} {}  & \\
&  \text{for } x \in (-2\pi n, -\pi(2n -1)) \sqcup (\pi(2n -1), 2\pi n)
\end{cases}
\] 
Denote by $C_n$ the graph $\{ (x, f_n(x)) \mid x \in \R \}$ of $f_n$ for any $n \in \Z_{\geq 0}$, and by $C_\infty$ the graph of $f_\infty = \cos$. 
For any point $x \in (-2\pi n, -\pi(2n -1)) \sqcup (\pi(2n -1), 2\pi n)$, we have the following relation: 
\[
\cos (x) + \dfrac{1}{n} \leq \varphi \left( \dfrac{\vert x\vert  - \pi(2n -1)}{\pi} \right) \left(  1 - \cos x \right) + \cos x  + \dfrac{1}{n} \leq 1+ \dfrac{1}{n} 
\]
Then the family $(C_n)_{n \in \Z_{\geq 0} \sqcup \{\infty\}}$ are pairwise disjoint. 
Let $D_n$ be the connected component of $\R^2 - \bigsqcup_{n \in \Z_{> 0} \sqcup \{\infty\}} C_n$ whose boundary is the union of $C_{n-1} \sqcup C_{n}$. 

Similarly, for any positive integer $n \in \Z_{>0}$, define a $C^\infty$ even function $g_n \colon \R \to [-1-1/n, 1-1/n]$ as follows:
\[
g_n(x) := 
\begin{cases}
-1 - \dfrac{1}{n} & \text{for } x \in \R - [-2\pi n, 2\pi n]\\
\cos x - \dfrac{1}{n} & \text{for } x \in [-2\pi(n-1), 2\pi(n -1)]\\
\varphi \left( \dfrac{\vert x\vert  - 2\pi(n -1)}{\pi} \right) \left(  -1 - \cos x \right) + \cos x  - \dfrac{1}{n}  \hspace{-150pt} {}  & \\
& {} \hspace{-35pt} {} \text{for } x \in (-2\pi n, -2\pi(n-1)) \sqcup (2\pi(n -1), 2\pi n)
\end{cases}
\] 
Denote by $C'_{n}$ the graph $\{ (x, g_n(x)) \mid x \in \R \}$ of $f_n$ for any $n \in \Z_{\geq 0}$. 
For any point $x \in (-2\pi n, -2\pi(n-1)) \sqcup (2\pi(n -1), 2\pi n)$, we have the following relation: 
\[
- 1- \dfrac{1}{n}  \leq \varphi \left( \dfrac{\vert x\vert  - 2\pi(n -1)}{\pi} \right) \left(  -1 - \cos x \right) + \cos x  - \dfrac{1}{n}  \leq \cos (x) - \dfrac{1}{n}
\]
Then the family $(C'_n)_{n \in \Z_{\geq 0}}$ are pairwise disjoint. 
Let $D_{-n}$ be the connected component of $\R^2 - \bigsqcup_{n \in \Z_{> 0}} C'_n$ whose boundary is the union of $C'_{n-1} \sqcup C'_{n}$. 
Then $\R^2 = \bigsqcup_{n \in \Z} D_n \sqcup C_\infty \sqcup \bigsqcup_{n \in \Z_{\geq 0}} C_n \sqcup C'_n$.

Define a non-singular vector field $X_0$ on $\R^2$ as follows: 
\[
X_0(x,y) := 
\begin{cases}
(1, 0) & \text{for } (x,y) \in D_0 \\ 
(1, f_\infty '(x)) = (1, - \sin(x)) & \text{for } (x,y) \in C_\infty \\ 
(1, f_n '(x)) & \text{for } (x,y) \in C_n \\ 
(1, g_n '(x)) & \text{for } (x,y) \in C'_n \\ 
\left(1, f_{n} '(x) + \varphi \left( \dfrac{y - f_{n}(x)}{f_{n-1}(x) - f_{n}(x)} \right)(f_{n-1} '(x) - f_{n} '(x))\right)   \hspace{-100pt} {}& \\
 & \text{for } (x,y) \in D_n \hspace{5pt} (n > 0) \\ 
\left(1, g_{n} '(x) + \varphi \left( \dfrac{g_{n}(x) - y}{g_{n}(x) - g_{n-1}(x)} \right)(g_{n-1} '(x) - g_{n} '(x))\right)   \hspace{-100pt} {}& \\
 & \text{for } (x,y) \in D_n \hspace{5pt} (n < 0) 
\end{cases}
\] 
By construction, for any point $p = (x,y) \in \bigsqcup_{n \in \Z_{\geq 0}} C_n$, we have a small \nbd $U_p$ of $p$ such that $X\vert_{U_p} = (1, f_n '(x))$. 
Similarly, for any point $p = (x,y) \in \bigsqcup_{n \in \Z_{\geq 0}} C'_n$, we have a small \nbd $U_p$ of $p$ such that $X_0\vert_{U_p} = (1, g_n '(x))$. 
This means that $X_0$ is an integrable continuous vector field on $\R^2$ such that $X_0$ is $C^\infty$ on  both $\R^2 - C_\infty$ and $C_\infty$. 
Take a diffeomorphism $h \colon (0,1) \times [-4,4] \to \R \times [-4,4]$ by $h(x,y) = (h_0(x), y) = \left(\tan (\pi(x-1/2)), y \right)$. 
Denote by $G$ the inverse image $h^{-1}(C_\infty)$. 
The pushforward $X_1 := (h^{-1})_*(X_0\vert_{\R \times [-4,4]})$ is a non-singular continuous  vector field on $B_0 := (0,1) \times [-4,4]$ such that $X_1$ is $C^\infty$ on $G$ and $B_0 - G$. 
By construction of $X_1$, the vector field $X_1$ generates a flow $v_{X_1}$ on $B_0$. 
%
Define $C^\infty$ functions $\phi_1, \phi_2 \colon \R \to [0,1]$ with 
\[\phi_1^{-1}(0) = \R - (0,1), \hspace{10pt} \phi_1^{-1}(1) = [1/3,2/3],\] 
\[\phi_2^{-1}(0) = \R - (-4,4), \text{ and } \phi_2^{-1}(1) = [-3,3]\]
such that $\phi_1$ (resp. $\phi_2$) is  increasing on $[0,1/3]$ (resp. $[-4,-3]$) and decreasing on $[2/3,1]$ (resp. $[3,4]$). 
Define an integrable continuous vector field $X_2$ on $\R^2 - \partial_{\pitchfork} \mathbb{D}$ as follows: 
\[
X_2(x,y) := 
\begin{cases}
\phi_1(x) \phi_2(y) X_1(x,y)/\vert X_1(x,y)\vert  & \text{for } (x,y) \in B_0 \\ 
0 &  \text{otherwise}\\
\end{cases}
\]
Then $X_2$ is $C^\infty$ on $\R^2 - C$. 
Define a $C^\infty$ function $\phi_3 \colon \R \to [0,1]$ with 
\[\phi_3^{-1}(0) = [0,1] \text{ and } \phi_3^{-1}(1) = \R - [-1/3,4/3]\] 
such that $\phi_3$ is decreasing on $[-1/3,0]$ and increasing on $[1,4/3]$. 
Write $B_{-1} := [-1/3,0) \times [-1,1]$ and $B_{1} := (1,4/3] \times [-1,1]$. 
Define $B'_{-1} := \{ (x,y) \mid x \in [-1/2,0), y \in [-1 - \phi_3(x),1+\phi_3(x)] \}$ and $B'_{1} := \{ (x,y) \mid x \in (1,3/2], y \in [-1 - \phi_3(x),1+\phi_3(x)]$. 
Then $B'_{-1}$ (resp. $B'_{1}$) is a closed \nbd of $B_{-1}$ (resp. $B_{1}$) on $\R^2 - \partial_{\pitchfork} \mathbb{D}$. 
Therefore there are a $C^\infty$ function $f_Y \colon \R^2 - \partial_{\pitchfork} \mathbb{D} \to [0,1]$ and a $C^\infty$ vector field $Y = (0, f_Y(x,y)y)$ on $\R^2 - \partial_{\pitchfork} \mathbb{D}$ such that 
\[Y\vert_{B_{-1} \sqcup B_{1}} = (0, \phi_3(x) y) \text{ and } Y^{-1}(0) = \R^2 - (\partial_{\pitchfork} \mathbb{D} \sqcup \mathop{\mathrm{int}}(B'_1 \sqcup B'_{-1})).\] 
There is a closed \nbd $U \subseteq h^{-1}(X_0^{-1}((1,0))) \cup ([-1,0] \sqcup [1,2]) \times ([-5,-1] \sqcup [1,5]) \subseteq \mathbb{B}$ of $\partial B_0 - \partial_{\pitchfork} \mathbb{D}$ on $\R^2 - \partial_{\pitchfork} \mathbb{D}$ with $U \cap (B'_1 \sqcup B'_{-1}) = \emptyset$ and there are a small positive number $\varepsilon \in (0, 1/3)$ and a $C^\infty$ function $\phi_4 \colon \R^2 \to [0,1]$ with 
\[((0,1)\times[-7/2,7/2]) \setminus \mathop{\mathrm{int}} U = \phi_4^{-1}(0) \text{ and }\]
\[\R^2 - ((-2/3,5/3) \times (-9/2,9/2)) \subset \phi_4^{-1}(1)\] 
such that $\phi_4\vert_{(-\varepsilon, 0) \times [-1,1]}(x,y) \leq -\phi_3(x)x$ on any $x \in (-\varepsilon, 0)$ and that $\phi_4\vert_{(0, \varepsilon) \times [-1,1]}(x,y) \leq \phi_3(x)(x-1)$ on any $x \in (0, \varepsilon)$. 
Define $C^\infty$ vector fields $Z$ and $X$ on $\R^2 - \partial_{\pitchfork} \mathbb{D}$ by $Z(x,y) = (\phi_4(x,y),0)$ and $X := X_2 + Y + Z$. 
Then 
\[
X(x,y) = 
\begin{cases}
X_2(x,y) = \phi_1(x) \phi_2(y) X_1(x,y)/\vert X_1(x,y)\vert  &\\
& {} \hspace{-155pt} \text{for } (x,y) \in B_0 \setminus U \\ 
X_2(x,y) + Z(x,y) = (\phi_1(x) \phi_2(y) + \phi_4(x,y), 0) &\\
& {} \hspace{-155pt}  \text{for } (x,y) \in (\R^2 -(B_0 \sqcup B'_{-1} \sqcup B'_{1} \sqcup \partial_{\pitchfork} \mathbb{D})) \cup U \\ 
Y(x,y) + Z(x,y) = (\phi_4(x,y), f_Y(x,y)y) &\\
& {} \hspace{-155pt}  \text{for } (x,y) \in B'_{-1} \sqcup B'_{1} \\ 
\end{cases}
\]
and $X(x,y)\vert_{B_{-1} \sqcup B_{1}}= (\phi_4(x,y), \phi_3(x) y)$ as in Figure~\ref{fig:support}.
\begin{figure}
\begin{center}
\includegraphics[scale=0.2]{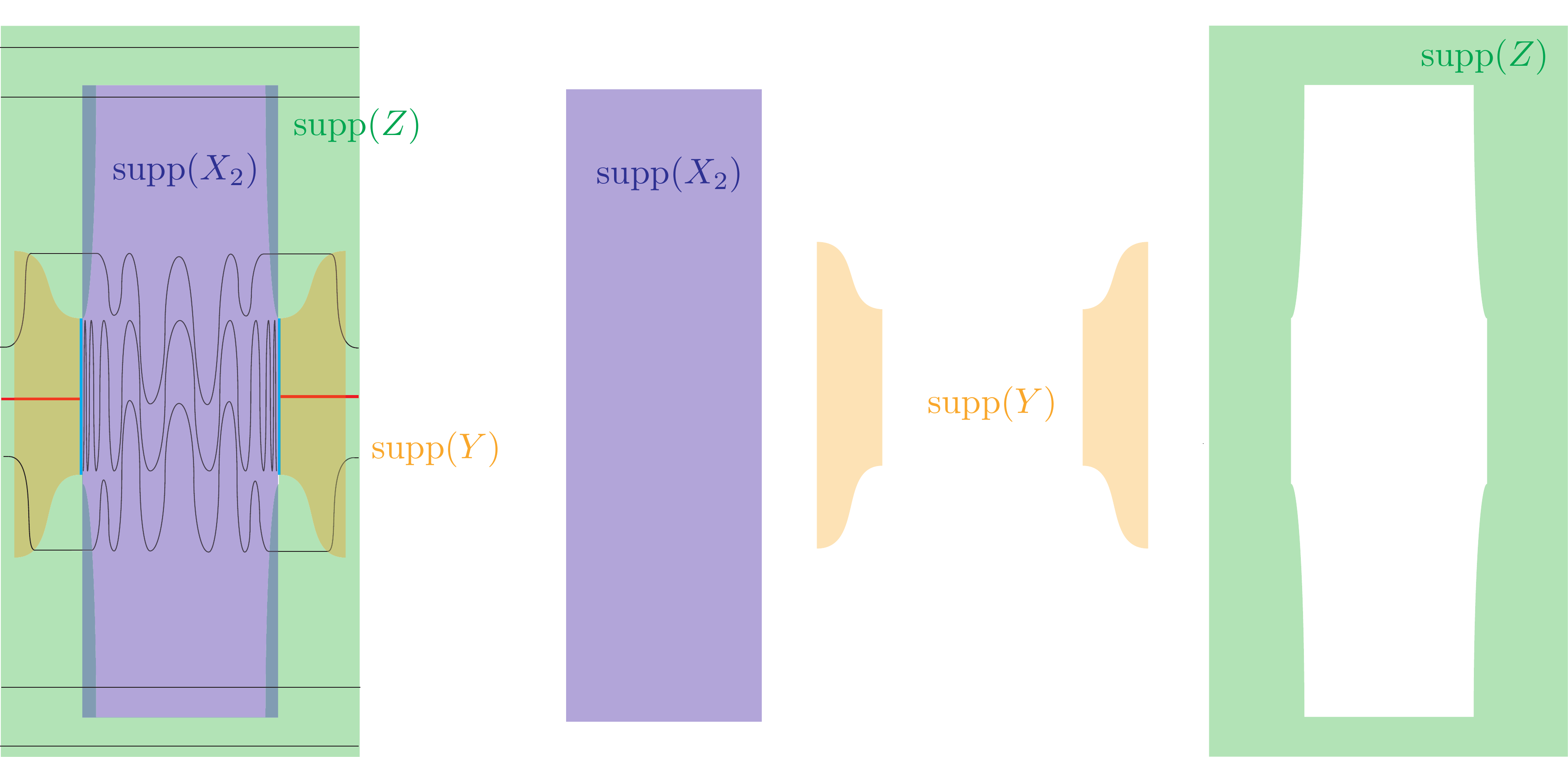}
\end{center}
\caption{Orbits with respect to $X$ and supports of vector fields $X_2$, $Y$, and $Z$.}
\label{fig:support}
\end{figure} 
Adding $\partial_{\pitchfork} \mathbb{D}$ as critical points, we extend $X$ into a vector field on $\R^2$, also denoted by $X$. 
Let $v_X$ be the $\R$-action generated by $X$ on $\R^2$.  
We will show that $X$ is as desired.

\begin{claim}
The flow $v_X$ is continuous on $\R^2$. 
\end{claim}
\begin{proof}
Since $\partial_{\pitchfork} \mathbb{D}$ is compact, for any $\varepsilon_0 >0$, there is a \nbd $V$ of $\partial_{\pitchfork} \mathbb{D}$ such that $\sup_{p \in V} \vert X(p)\vert  < \varepsilon_0$. 
This implies the continuity of $v_X$ at any points in the closed subset $\R \times \partial_{\pitchfork} \mathbb{D} \subseteq \R \times \Sv$. 
On the other hand, by the openness of $\R \times (\R^2 - \partial_{\pitchfork} \mathbb{D})$ and the invariance of $\R^2 - \partial_{\pitchfork} \mathbb{D}$, the $\R$-action $v_X$ is continuous at any points in the complement $\R \times (\R^2 - \partial_{\pitchfork} \mathbb{D})$. 
This means that $v_X$ is continuous on $\R^2$. 
\end{proof}

\begin{claim}
Assertion {\rm(4)} holds. 
\end{claim}
\begin{proof}
We have that $X(x,y) = X_2(x,y) \neq 0$ on $B_0 \setminus U$, $X(x,y) = X_2(x,y) + Z(x,y) \neq 0$ on $(\R^2 -(B_0 \sqcup B'_{-1} \sqcup B'_{1} \sqcup \partial_{\pitchfork} \mathbb{D})) \cup U$, and the first component of $X(x,y)$ is  $Z(x,y) = \phi_4(x,y) \neq 0$ on $B'_{-1} \sqcup B'_{1}$. 
This means that $X(x,y) \neq 0$ on $\R^2 - \partial_{\pitchfork} \mathbb{D}$. 
\end{proof}

\begin{claim}
Assertions {\rm(1)} and {\rm(2)} hold. 
\end{claim}
\begin{proof}
Since $U \subseteq \mathbb{B}= [-1,2] \times [-5,5]$, $B_0 = (0,1) \times [-4,4]$, $B_{-1} = [-1/3,0) \times [-1,1] \subset B'_{-1}$, and $B_{1} = (1,4/3] \times [-1,1] \subset B'_1$, we have $B_0 \cup (B'_{-1} \sqcup B'_1 \sqcup U) \subset \mathbb{B}$. 
Then $X(x,y) = Z(x,y) = (\phi_4(x,y), 0) = (1, 0)$ on $\R^2 - ((-2/3,5/3) \times (-9/2,9/2)) \subset \phi_4^{-1}(1)$. 
By $\R^2 - \mathop{\mathrm{int}} \mathbb{B} = \R^2 - ((-1,2) \times (-5,5)) \subset \R^2 - ((-2/3,5/3) \times (-9/2,9/2))$, the closed square $\mathbb{B}$ is a flow box with respect to $X$ such that the restriction $X\vert_{\R^2 - \mathbb{B}}$ is $(1,0)$. 
\end{proof}

Because $\R^2$ can be identified with a subset of the sphere by adding a point at infinty, the flow $v$ can be identified with the restriction of the resulting flow of $v$ considering the point at infinty as a singular point. 
Since $\mathbb{B}$ can be identified with a flow box in a compact surface, by Gutierrez's smoothing theorem~\cite{gutierrez1986smoothing}, we may assume that $v$ is $C^\infty$ and so is $X$. 

\begin{claim}
Assertion {\rm(3)} holds.
\end{claim}
\begin{proof}
Since $G$ is an orbit of $X_1$, from $X = X_2 = \phi_1(x) \phi_2(y) X_1(x,y)/\vert X_1(x,y)\vert $ on $B_0 \subset \mathbb{B} \setminus U$, the subset $G$ is the orbit of $X$. 
By $X_2 = Y = 0$ on $J_0 = (\R - [0,1]) \times \{ 0 \}$, we obtain that $X(x,y) = Z(x,y) = (\phi_4(x,y), 0)$ and $\phi_4(x,y)>0$ on $J_0$. 
This means that subsets $(- \infty,0) \times \{ 0 \}$ and $(1, \infty) \times \{ 0 \}$ are orbits of $X$ respectively. 
\end{proof}

We show assertion {\rm(5)}. 
On $[-1,2] \times ([-5,-4] \sqcup [4,5])$, the vector field $X = (\phi_4(x,y), 0)$ is non-singular and so the subset $[-1,2] \times \{y_0\}$ for any $y_0 \in [-5,-4] \sqcup [4,5]$ is an orbit arc for $X$.

\begin{claim}\label{claim:complete_intersection}
The positive orbit in any point in $\{-1\} \times ([-5, 5] - \{ 0 \}) \subset \partial \mathbb{B}$ intersects $\{0\} \times ([-5,5] - [-1,1])$. 
\end{claim}
\begin{proof}
Fix a point $p_0 = (x_0,y_0) \in \{-1\} \times ([-5, 5] - \{ 0 \})$. 
Suppose that $y_0 >1$. 
Then $[-1,0] \times [y_0,5]$ is compact and $Z(x,y) = (\phi_4(x,y),0)$ is non-singular on the domain $[-1,0] \times [y_0,5]$.  
By $\mathrm{supp}(Y) \subset ([-1,0] \sqcup [1,2]) \times [-2,2]$, $Y(x,y) = (0, f_Y(x,y)y)$, and $f_Y(x,y)y \geq 0$ on $[-1,0] \times [y_0,2]$, the positive orbit $O^+((x_1,y_0))$ for any $x_1 \in [-1,0)$ intersects $\{0\} \times [y_0,5] \subset \{0\} \times (1,5]$. 
In particular, the positive orbit $O^+(p_0)$ intersects $\{0\} \times [y_0,5] \subset \{0\} \times (1,5]$. 
Suppose that $y_0 = 1$. 
Then $X = (\phi_4(x,1), 0)$ on $x \in [-1,-1/2]$, $X = (\phi_4(x,1), f_Y(x,1))$ on $x \in [-1/2,0)$, and $f_Y(x,1) > 0$ on $x \in (-1/2,0)$. 
Therefore $O^+(p_0)$ intersects $[-1,0] \times (1,5]$ and so $\{0\} \times (1,5]$ because of the previous argument. 
Suppose that $y_0 \in (0,1]$. 
Then the positive orbit $O^+(p_0)$ intersects either $[-1,0] \times (1,5]$ or $(-\varepsilon, 0) \times (0,1)$. 
If $O^+(p_0)$ intersects $[-1,0] \times (1,5]$, then the previous argument implies that  the positive orbit $O^+(p_0)$ intersects $\{0\} \times [y_0,5] \subset \{0\} \times (1,5]$. 
Thus we may assume that $O^+(p_0)$ intersects $(-\varepsilon, 0) \times (0,1)$. 
By definition, we have that $Y(x,y) = (0, \phi_3(x) y)$ and $\phi_4(x,y) \leq -\phi_3(x)x$ on $(-\varepsilon, 0) \times (0,1) \subset B_{-1}$. 
Then $X(x,y) = (\phi_4(x,y), \phi_3(x) y) = \phi_3(x)(-x, y) + (\phi_4(x,y)+ \phi_3(x)x, 0)$ and $\phi_4(x,y) +\phi_3(x)x \leq 0$ on $(-\varepsilon, 0) \times (0,1) \subset B_{-1}$. 
Since the orbit of $p$ with respect to the vector field $Y'$ on $B_{-1}$ defined by $Y'(x,y) = (-x,y)$ intersects a point $(x',y')$ in the horizontal boundary $[-1/3,0) \times \{ 1 \}$, from $\phi_4(x,y) \geq 0$, the orbit $O^+(p_0)$ with respect to the vector field $X\vert_{(-\varepsilon, 0) \times (0,1)} = (\phi_4(x,y), \phi_3(x) y) = \phi_3(x)(-x, y) + (\phi_4(x,y)+ \phi_3(x)x, 0)$ intersects a point in $[-1/3, x'] \times \{1\}$. 
From the previous argument, the orbit $O^+(p_0)$ intersects $\{0\} \times (1,5]$. 
By symmetry, if $y_0 <0$, then the orbit $O^+(p_0)$ intersects $\{0\} \times [-5,-1)$. 
\end{proof}

By symmetry, the negative orbit in any point in $\{2\} \times ([-5, 5] - \{ 0 \}) \subset \partial \mathbb{B}$ intersects $\{1\} \times ([-5,5] - [-1,1])$. 
By construction, the positive (resp. negative) orbit of any point in $B_0 - G$ with respect to $X$ intersects $U$ and so $\{1\} \times ([-5,5] - [-1,1])$ (resp. $\{0\} \times ([-5,5] - [-1,1])$). 
By Claim~\ref{claim:complete_intersection} and its dual statement, each of the positive and negative orbits of any points in $\mathbb{B} - (G \sqcup \partial_{\pitchfork} \mathbb{D} \sqcup I_- \sqcup I_+)$ is neither singular nor periodic but is unbounded. 
This implies assertion {\rm(5)}. 
Finally, we show assertion {\rm(6)}.

\begin{claim}
Assertion {\rm(6)} holds.
\end{claim}
\begin{proof}
Fix a convergence sequence $(y_n)_{n \in \Z_{\geq 0}}$ of non-zero numbers $y_n$ to $0$. 
Denote by $z_n$ the point with $\{ z_n \} = O^+(-1,y_n) \cap (\{0\} \times (\R - [-1,1]))$. 
By construction, the sequence $(z_n)_{n \in \Z_{\geq 0}}$ converges to either $(0,1)$ or $(0,-1)$. 
Denote by $w_n$ the point with $\{ w_n \} = O^+(z_n) \cap (\{1/2\} \times \R)$. 
By construction, the sequence $(w_n)_{n \in \Z_{\geq 0}}$ converges to a point $p_G$ in $G$. 
For any point $p_g \in G$, there is a number $t_{p_g} \in \R$ with $p_g = v_X(t_{p_g}, p_G) = \lim_{n \to \infty} v_X(t_{p_g},w_n)$. 
\end{proof}
Therefore $X$ is as desired. 
\end{proof}

Considering a $C^\infty$ bump function $\varphi_0 \colon \R^2 \to [0,1]$ with $\varphi_0^{-1}(0) = \{ 0 \}$ and $\R^2 - (-1/2,1/2)^2 \subset \varphi_0^{-1}(1)$, the flow generated by the vector field $(\varphi_0,0)$ on $\R^2$ is called the resulting flow of a unit vector field $(1,0)$ on $\R^2$ by replacing a non-singular point with a singular point. 
Lemma~\ref{lem:deformation_01}{\rm)}.

The previous lemma implies the following statement.

\begin{lemma}\label{lem:blowdown_box}
Let $v_X$ be the flow generated by the vector field $X$ as in Lemma~\ref{lem:no_arc_conn_flowbox} and $p_{\mathcal{F}_C} \colon \R^2 \to \R^2/\mathcal{F}_C$ be the quotient map as in Lemma~\ref{lem:collapsing}. 
Then the following statements hold: 
\\
{\rm(1)} The mapping $v \colon \R \times \R^2/\mathcal{F}_C \to \R^2/\mathcal{F}_C$ defined by 
\[
v(t,p) := p_{\mathcal{F}_C}(v_X(t,p_{\mathcal{F}_C} ^{-1}(p)))
\]
is well-defined and continuous. 
\\
{\rm(2)} The mapping $v$ is semi-conjugate to $v_X$ via $p_{\mathcal{F}_C}$. 
\\
{\rm(3)} The mapping $v$ is topologically equivalent to the resulting flow of a unit vector field $(1,0)$ on $\R^2$ by replacing a non-singular point with a singular point. 
\end{lemma}

\begin{proof}
First, we show the well-definedness of $v$. 

\begin{claim}
The mapping $v$ is well-defined. 
\end{claim}
\begin{proof}
Fix $y \in \R^2/\mathcal{F}_C$. 
Suppose that $y \in p_{\mathcal{F}_C}(C)$. 
Then $p_{\mathcal{F}_C} ^{-1}(y) = C$ and so $v_X(t,p_{\mathcal{F}_C} ^{-1}(y)) = v_X(t,C) = C$. 
Therefore we have $p_{\mathcal{F}_C}(v_X(t,p_{\mathcal{F}_C} ^{-1}(y))) = p_{\mathcal{F}_C}(C) = y$. 
This means that $y$ is a singular point of $v$. 
Suppose that $y \not\in p_{\mathcal{F}_C}(C)$. 
Since $p_{\mathcal{F}_C}^{-1}(y) =y$, we have $v(t,y) = p_{\mathcal{F}_C}(v_X(t,p_{\mathcal{F}_C} ^{-1}(y))) = p_{\mathcal{F}_C}(v_X(t,y)) =v_X(t,y)$. 
This means that $v$ is well-defined. 
\end{proof}

\begin{claim}
The mapping $v$ is continuous. 
\end{claim}
\begin{proof}
Since the quotient map $p_{\mathcal{F}_C}$ is continuous and closed, the map $1_{\R} \times p_{\mathcal{F}_C} \colon \R \times \R^2 \to \R \times \R^2/\mathcal{F}_C$ defined by $1_{\R} \times p_{\mathcal{F}_C}(t,x) = (t,p_{\mathcal{F}_C}(x))$ is a quotient map and so is closed. 
Then $v(t,x) = p_{\mathcal{F}_C} \circ v_X ( (1_{\R} \times p_{\mathcal{F}_C})^{-1}(t,x))$. 
For any closed subset $A \subseteq \R^2/\mathcal{F}_C$, the inverse image $v^{-1}(A) = 1_{\R} \times p_{\mathcal{F}_C}(v_X^{-1}(p_{\mathcal{F}_C}^{-1}(A)))$ is closed. 
This means that the $\R$-action $v$ is a flow. 
\end{proof}

By construction, for any $(t,x) \in \R \times \R^2$, we obtain $v(t, p_{\mathcal{F}_C}(x)) = p_{\mathcal{F}_C}(v_X(t,x))$. 
This implies the semi-conjugacy. 
Therefore assertion (2) holds.

\begin{claim}
The flow $v$ is topologically equivalent to the resulting flow of the unit vector field $(1,0)$ on $\R^2$ by replacing a non-singular point with a singular point. 
\end{claim}
\begin{proof}
Since the quotient space $\R^2/\mathcal{F}_C$ is homeomorphic to $\R^2$ and the composition $p_{\mathcal{F}_C} = q \circ p_{\mathcal{F}} \colon \R^2 \to \R^2/\mathcal{F}_C$ of the quotient maps is continuous and is the quotient map, the set difference $\R^2/\mathcal{F}_C - [C] = (\R^2 - C)/\mathcal{F}_C$ is homeomorphic to an open annulus $\R^2 - \{ 0 \}$.
Because $\R^2$ can be identified with a subset of the sphere by adding a point at infinty, the flow $v$ can be identified with the restriction of the resulting flow of $v$ considering the point at infinty as a singular point. 
By Gutierrez's smoothing theorem~\cite{gutierrez1986smoothing}, we may assume that the flow $v$ is topologically equivalent to a $C^1$-flow. 
Put $p_- := (-3,0)$ and $p_+ := (3,0)$. 
Since the restriction $X\vert_{\R^2 - \mathbb{B}}$ is $(1,0)$, by $\mathbb{B} := [-1,2] \times [-5,5]$, we have that $\omega_v(p_-) = \{ [C] \}$ and $\alpha_v(p_+) = \{ [C] \}$ and that the orbits $O_v([(-3,y)])$ are closed subsets in $\R^2/\mathcal{F}_C$ for any $y \neq 0 \in \R$. 
Then the union $L_0 := O_v([(-3,0)]) \sqcup \{[C]\} \sqcup O_v([(3,0)])$ is a piecewise $C^1$-line. 
Fix a Riemannian metric on the plane $\R^2/\mathcal{F}_C$. 
Define a homeomorphism $D_0 \colon L_0 \to \R$ as follows: 
The value $D_0([(x,y)])$ for any point $[(x,y)] \in L_0$ with $y \geq -3$ is the arc-length of the arc connecting $[-3,0]$ and $[(x,y)]$ in $L_0$, and the value $D_0([(x,y)])$ for any point $[(x,y)] \in L_0$ with $y \leq -3$ is the arc-length of the arc connecting $[-3,0]$ and $[(x,y)]$ in $L_0$ multiplied by minus. 
Similarly, for any $y \neq 0 \in \R$, denote by $L_y$ the orbit $O_v([(-3,y)])$ and define a homeomorphism $D_y \colon L_y \to \R$ as follows: 
The value $D_y([(x,y')])$ for any point $[(x,y')] \in L_y$ with $y' \geq -3$ is the arc-length of the arc connecting $[-3,0]$ and $[(x,y')]$ in $L_{y}$, and the value $D_y([(x,y')])$ for any point $[(x,y')] \in L_y$ with $y' \leq -3$ is the arc-length of the arc connecting $[-3,0]$ and $[(x,y')]$ in $L_y$ multiplied by minus. 
By construction, the mapping $h \colon \R^2/\mathcal{F}_C \to \R^2$ defined by $h([(x,y')]) := (D_y([(x,y')]), y)$ if $[(x,y')] \in O_v([(-3,y)])$ is a continuous bijection. 
Moreover, the continuous mapping $h$ can be continuously extend to the spheres which are the one-point compactifications of $ \R^2/\mathcal{F}_C$ and $\R^2$ respectively. 
Since any continuous bijection from a compact space to a Hausdorff space is homeomorphic, the extension of $h$ is homeomorphic and so is the restriction $h$. 
Every orbit of the induced flow $v_h \colon \R \times \R^2 \to \R^2$ defined by $v_h(t,x,y) := h(v(t,h^{-1}(x,y)))$ is either the origin $\{0 \}$, a negative half of $x$-axis $\{ 0 \} \times \R_{<0}$, a positive half of $x$-axis $\{ 0 \} \times \R_{>0}$, or a horizontal line $\{ y \} \times \R$ for some $y \neq 0 \in \R$. 
Therefore $v_h$ is the resulting flow on $\R^2/\mathcal{F}_C$ of the unit vector field $(1,0)$ on $\R^2$ by replacing the non-singular point $h([C])$ with a singular point. 
\end{proof}
This completes the proof. 
\end{proof}

We demonstrate Theorem~\ref{thm:deformation_02} as follows. 

\begin{proof}[Proof of Theorem~\ref{thm:deformation_02}]
Let $v$ be a flow on a surface $S$ with an $\omega$-limit set $\omega$ containing non-singular point  $p_0$ and with a point $q_0 \in S - \omega$ satisfying $\omega(q_0) = \omega$. 

Take any closed trivial flow box $B_{p_0}$ with $p_0 \in \mathop{\mathrm{int}} B_{p_0}$ and $q_0 \notin B_{p_0}$. 
Identify $B_{p_0}$ with $[-1,2] \times [-5,5]$, $p_0$ with $0$, the set of orbit arcs in $B_{p_0}$ with $\{ [-1,2] \times \{ y \} \mid y \in [-5,5] \}$, and the connected component of $O(p_0) \cap B_{p_0}$ containing $p_0$ with $[-1,2] \times \{ 0 \}$. 
Then any connected components of $\omega(q_0) \cap B_{p_0}$ are of form $[-1,2] \times \{ y \}$ for some $y \in [-5,5]$. 
Replacing $B_{p_0}$ with the flow box $\mathbb{B}$ constructed in Lemma~\ref{lem:no_arc_conn_flowbox}, denote by $w$ the resulting flow.  
Lemma~\ref{lem:blowdown_box} implies that the resulting flow $v_{p_0}$ of $v$ by replacing $p_0$ with a singular point is topologically equivalent via $h$ to a flow $w_C$ on $S/\mathcal{F}_C$ which is topologically semi-conjugate to $w$, where $h \colon S \to S/\mathcal{F}_C$ is the homeomorphism constructed in the proof of Lemma~\ref{lem:blowdown_box} and the quotient space $S/\mathcal{F}_C$ is the resulting surface of $S$ by collapsing the closed invariant subset $C$ of $\omega_w(q_0)$ into a singleton. 
Then the restriction $v\vert_{S - \omega}$ is topologically equivalent to the restriction $w\vert_{S - \omega_w(q_0)}$. 
These mean that assertions (2)--(4) hold. 

Finally, we show the absence of arcwise-connectivity of $\omega_w(q_0)$. 
By Lemma~\ref{lem:no_arc_conn_flowbox}{\rm(6)}, the $\omega$-limit set $\omega_w(q_0)$ contains $C$.  
Therefore the disjoint union $([-1,0) \times \{ 0 \}) \sqcup C \sqcup ((1,2] \times \{ 0 \})$ is contained in a connected component of $\omega_w(q_0) \cap \mathbb{B}$. 
Since $h \colon S \to S/\mathcal{F}_C$ is the homeomorphism, we may assume that the restriction $h|_{S - \{ p_0 \}} \colon S - \{ p_0 \} \to (S -C)/\mathcal{F}_C$ is identical. 
Let  $p_C \colon S \to S/\mathcal{F}_C$ be the quotient map collasping $C$ into the singleton $[C]$ as in Lemma~\ref{lem:collapsing}. 
Then the restriction $p_C|_{S - C} \colon S -C \to (S - C)/\mathcal{F}_C$ is identical.

\begin{claim}\label{claim:tot_disconnected}
The intersection $(\{-1 \} \times [-5, 5]) \cap \omega_w(q_0) = (\{-1 \} \times [-5, 5]) \cap \omega(q_0)$ is totally disconnected. 
\end{claim}
\begin{proof}
Since $q_0 \in S - \omega = S - \omega(q_0)$, we have $q_0 \notin \omega(q_0) = \omega$ and so $O(q_0) \cap \omega(q_0) = \emptyset$. 
By definition of $v_{p_0}$, we obtain that $O_v^+(q_0) = O_{v_{p_0}}^+(q_0)$, and so that $\omega_v(q_0) = \omega = \omega_{v_{p_0}}(q_0)$ and $\overline{O_{v_{p_0}}^+(q_0)} = O_{v_{p_0}}^+(q_0) \sqcup \omega = \overline{O^+(q_0)}$. 
Since the restriction $v\vert_{S - \omega}$ is topologically equivalent to the restriction $w\vert_{S - \omega_w(q_0)}$, we have that $\overline{O^+(q_0)} \cap (S - \omega) = O^+(q_0)  = O_w^+(q_0) = \overline{O_w^+(q_0)} \cap (S - \omega_w(q_0))$ and so that $O_w^+(q_0) \cap \omega_w(q_0) = \emptyset$. 
Then $\overline{O_w^+(q_0)} \subseteq \omega_w(q_0) \sqcup O_w^+(q_0) \sqcup \{ q_0 \}$. 
By $\omega_w(q_0) \sqcup O_w^+(q_0) \sqcup \{ q_0 \} \subseteq \overline{O_w^+(q_0)}$, we have $\overline{O_w^+(q_0)} = O_w^+(q_0) \sqcup \{ q_0 \} \sqcup \omega_w(q_0)$. 
If $\omega_w(q_0)$ is locally dense, then $\omega_w(q_0)$ is a \nbd of a point of $\omega_w(q_0)$ and so $O_w^+(q_0) \cap \omega_w(q_0) \neq \emptyset$, which contradicts $O_w^+(q_0) \cap \omega_w(q_0) = \emptyset$. 
Thus $\omega_w(q_0)$ is not locally dense. 
Then the proof of the claim is completed. 
\end{proof}

Put $C' := [-1,0) \times \{ 0 \}) \sqcup C \sqcup ((1,2] \times \{ 0 \}$. 

\begin{claim}\label{claim:cc}
The disjoint union $C' \subset \mathbb{B} = B_{p_0} = [-1,2] \times [-5,5]$ is a connected component of $\omega_w(q_0) \cap \mathbb{B}$. 
\end{claim}
\begin{proof}
Assume that there is a connected component $K$ of $\omega_w(q_0) \cap \mathbb{B}$ with $C' \subsetneq K$. 
By Claim~\ref{claim:tot_disconnected}, the set difference $K - C'$ consists of closed intervals connecting pairs of points $(-1,y)$ and $(2,y)$ for some $y \in [-5,5]$. 
Therefore, the image $p_C(K - C')$ contains at least one interval and so the image $p_C(K)$ contains at least two intervals. 
%
By constructions of $C$ and $p_C$, the image $p_C(C')$ is a closed interval in $\mathbb{B}$ between points $(-1, 0), (2,0) \in \partial \mathbb{B}$.
Because $p_C|_{S - C}$ is homeomorphic and $p_C|_{\partial \mathbb{B}}$ is identical, the image $p_C(K - C')$ consists of closed intervals. 
Since $K$ is connected, the image $p_C(K)$ is connected and consists of pairwise disjoint closed intervals whose interiors are contained in $\mathop{\mathrm{int}}\mathbb{B}$ and which connect $\partial \mathbb{B}$. 
Since every connected component of $p_C(K)$ is a closed interval in $\mathbb{B}$ connecting a pair of points $(-1,y)$ and $(2,y)$ for some $y \in [-5,5]$ and since $(\{-1 \} \times [-5, 5]) \cap K) = (\{-1 \}$ is totally disconnected, the image $p_C(K)$ is a closed interval, which contradicts that $p_C(K)$ contains at least two intervals. 
Thus $C'$ is a connected component of $\omega_w(q_0) \cap \mathbb{B}$. 
\end{proof}

Since the restrictions $p_C|_{S - C}$ and $h\vert_{\mathbb{B} - C}$ can be identified with the identical maps on $S-C$, by $O^+_w(q_0) \cap C = \emptyset$, we have that $h^{-1} \circ p_C(O_w^+(q_0)) = h^{-1}(O_{w_C}^+([q_0])) = O_{v_{p_0}}^+(q_0) = O_{v}^+(q_0)$. 
Since the restriction $p_{\mathcal{F}}\vert_{\mathbb{B} - C}$ is identical, the set difference $\mathbb{B} - C$ are homeomorphic to $\mathbb{B} - p_{\mathcal{F}}(C)$. 
Because $C$ is not arcwise-connected, from $C \subset \omega_w(q_0)$, by Claim~\ref{claim:cc}, neither is the $\omega$-limit set $\omega_w(q_0)$. 
\end{proof}

\section{Examples}

We describe some kinds of $\omega$-limit sets that appear in Theorem~\ref{main:a}.  

\subsection{Nontrivial quasi-Q-sets}\label{ex01}

We show that there is a toral flow with a non-locally-dense nontrivial quasi-Q-set as follows.

\begin{lemma}
There are a toral flow $v_{\varphi}$ and a point $z$ whose $\omega$-limit set is a non-locally-dense nontrivial quasi-Q-set such that $\overline{\mathop{\mathrm{Cl}}(v_{\varphi})} \neq \Omega(v_{\varphi})$, where $\Omega (v_{\varphi})$ is the non-wandering set of $v_{\varphi}$.
\end{lemma}

\begin{proof}
Consider a Denjoy diffeomorphism $f : \mathbb{S}^1 \to \mathbb{S}^1$ with an exceptional minimal set $\mathcal{C}$.
Let $v_f$ be the suspension of $f$ on the torus $\T^2 := (\mathbb{S}^1 \times \R)/(x, r) \sim (f(x), r+1)$ and $\mathcal{M}$ the minimal set of $v_f$.

We will replace the minimal set $\mathcal{M}$ of $v_f$ with a union of singular points and separatrices of the resulting flow $v_{\varphi}$ as follows. 
Fix a bump function $\varphi : \T^2 \to \R_{\geq 0}$ with $\varphi^{-1}(0) = \hat{\mathcal{C}}$, where $\hat{\mathcal{C}} := \mathcal{M} \cap (\mathbb{S}^1 \times \{ 1/2 \})$ is a lift of $\mathcal{M}$. 
Let $X$ be the continuous vector field generating $v_f$ on the mapping torus $\mathbb{T}^2$. 
Since $\hat{\mathcal{C}}$ are covered by finitely many trivial flow boxes, there is an open \nbd $U$ of $\hat{\mathcal{C}}$ such that the restriction $X\vert_U$ can be considered as the restriction of an integrable continuous vector field $Y$ on a sphere. 
By Gutierrez's smoothing theorem~\cite{gutierrez1986smoothing}, we may assume that the vector field $Y$ is $C^\infty$ and so does $X\vert_U$. 
Since every closed subset of any paracompact $C^\infty$ manifold is a zero set of some $C^\infty$ function on it, take a $C^\infty$ bump function $\varphi \colon \mathbb{T}^2 \to [0,1]$ with $\varphi^{-1}(0) = \hat{\mathcal{C}}$ and $\varphi\vert_{\mathbb{T}^2 - U } = 1$. 
Define a continuous vector field $Z$ by $Z(p) := \varphi(p) X(p)$. 
Then the restriction $Z\vert_{\mathbb{T}^2 - U} = X\vert_{\mathbb{T}^2 - U}$ is non-singular and the restriction $Z\vert_{U} = X\vert_{U}$ is $C^\infty$. 
Therefore $Z$ is locally Lipschitz continuous and so generates a flow $v_{\varphi}$ with $\mathcal{M} = \mathop{\mathrm{Sing}(v_{\varphi})} \sqcup \{ \text{separatrix of } v_{\varphi} \}$, $O_{v_{\varphi}}(p) = O_{v_f}(p)$, and $\omega_v(p) = \omega_{v_{\varphi}}(p) =  \mathcal{M}$ for any point $p \in \T^2 - \mathcal{M}$.
Then $\T^2 =  \mathop{\mathrm{Sing}}(v_{\varphi}) \sqcup \mathrm{P}(v_{\varphi})$ and $\overline{\mathop{\mathrm{Cl}}}(v_{\varphi}) =  \mathop{\mathrm{Sing}}(v_{\varphi}) = \hat{\mathcal{C}}  \neq \mathcal{M} = \Omega (v_{\varphi})$,
where $\mathrm{P}(v_{\varphi})$ is the union of non-recurrent orbits of $v_{\varphi}$.
\end{proof}

\subsection{Quasi-circuit that is not a circuit}\label{ex02}

\begin{figure}
\begin{center}
\includegraphics[scale=0.17]{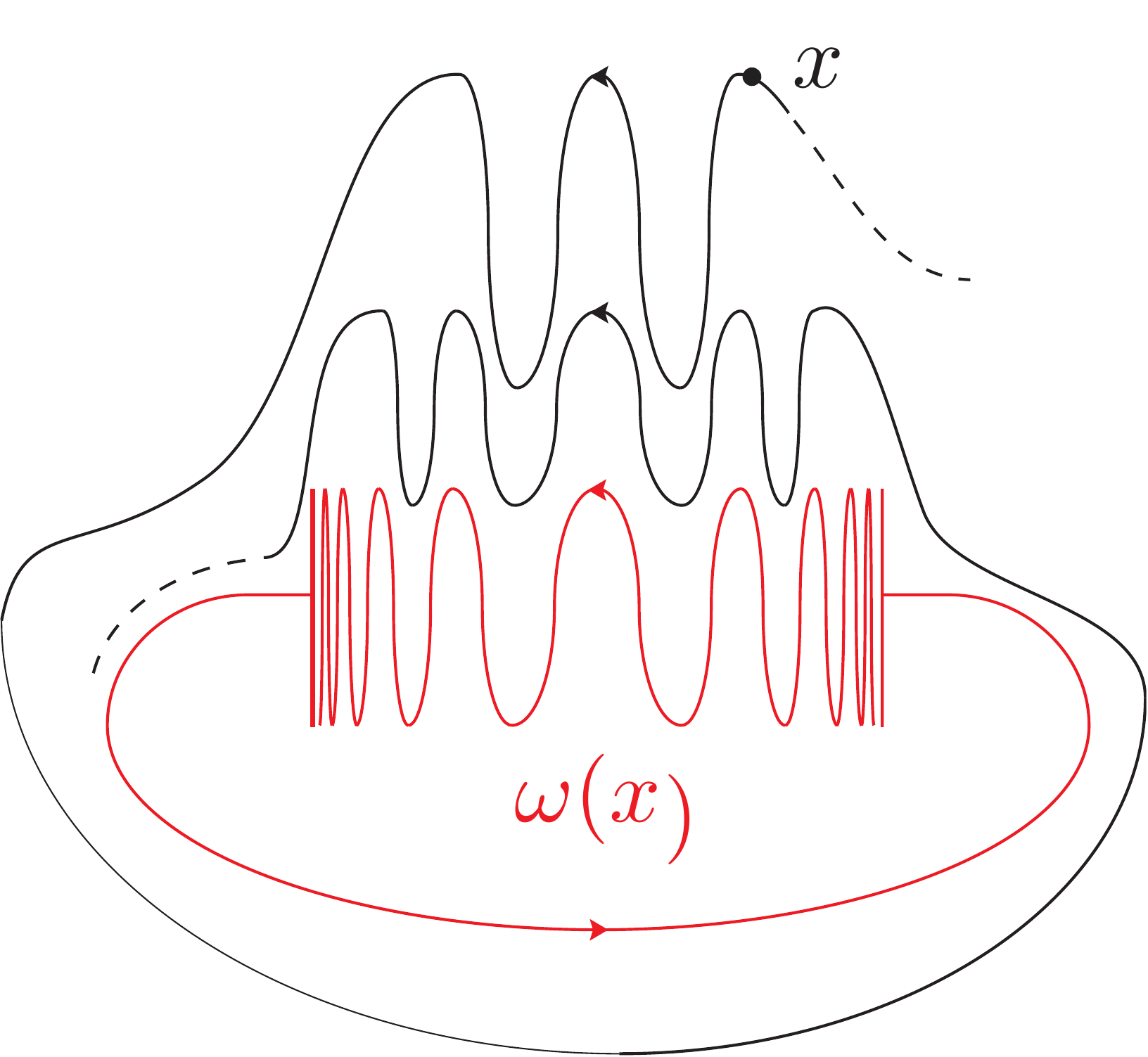}
\end{center}
\caption{An $\omega$-limit set which is a non-locally-connected quasi-circuit, which is neither the image of a circle nor a circuit.}
\label{nonac}
\end{figure}

Using the flow box in Lemma~\ref{lem:no_arc_conn_flowbox}, we can construct a flow with a quasi-circuit that is not a circuit such that it consists of two non-recurrent orbits and two closed intervals contained in the singular point set as in Figure~\ref{nonac}. 
In particular, the quasi-circuit is homeomorphic to the union of $G$ and a curve from a point in $G$ and to a point in $G$. 

\begin{figure}
\begin{center}
\includegraphics[scale=0.3]{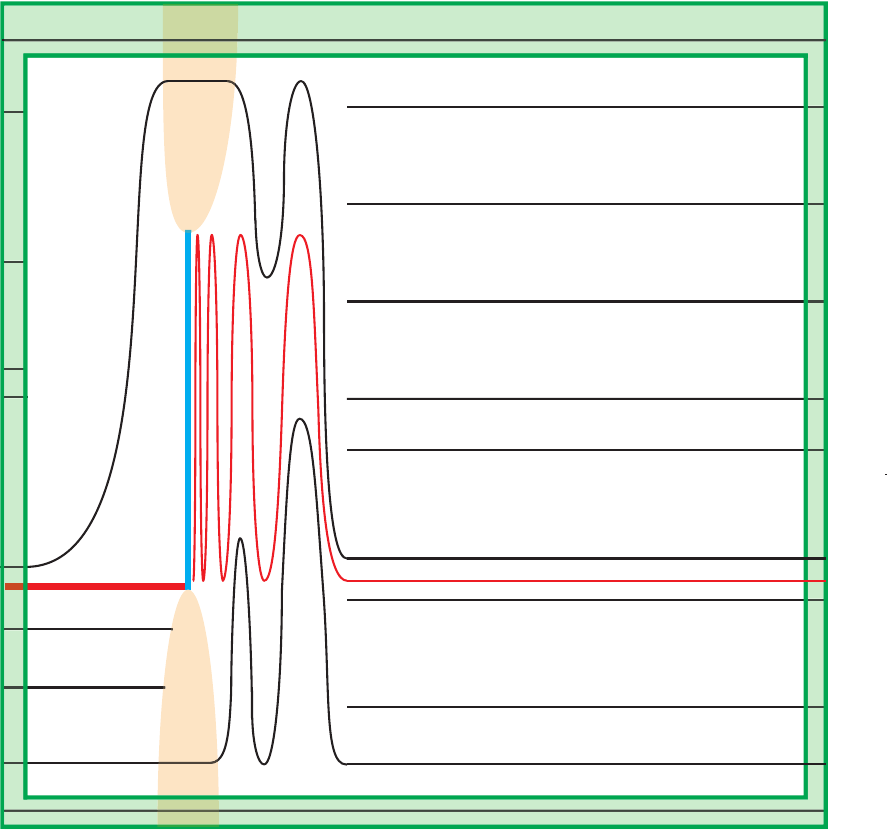}
\end{center}
\caption{A flow box with a non-arcwise-connected subset.}
\label{fig:vf_X2}
\end{figure}

Note that  Hastings constructed an attractor of a flow on $\R^2$ which is homeomorphic to a Warsaw circle (i.e. the disjoint union of the graph of the function $f \colon (0,1/\pi] \to [-1,1]$ by $f(x) = \sin 1/x$, the interval $\{ 0 \} \times [-1,1]$, and an open arc from $(0,-1)$ to $(1/\pi, 0)$) but is not an $\omega$-limit set in \cite[Example 3.3]{Hastings1979PB}.  
By a similar construction of the flow box in Lemma~\ref{lem:no_arc_conn_flowbox}, we can construct a flow box as in Figure~\ref{fig:vf_X2} and so a flow with a quasi-circuit that is homeomorphic to a Warsaw circle. 

\subsection{Non-locally connected subsets of singular points}

A continuum $\mathcal{M}$ on a surface is a {\bf transversely Cantor set} if there is a flow on a surface with a transversely Cantor Q-set which is homeomorphic to $\mathcal{M}$.  
We construct the following examples of flows. 

\begin{lemma}
There is a toral flow with an $\omega$-limit set consisting of singular points which is a transversely Cantor set.
\end{lemma}

\begin{proof}
Let be a non-recurrent orbit $O$ whose $\omega$-limit set is an exceptional minimal set $\mathcal{M}$ in a Denjoy flow on a torus $\mathbb{T}^2$. 
By Lemma~\ref{lem:deformation_01}, take the resulting flow of $v_x$ by replacing a non-singular point $x \in \mathcal{M}$ with a singular point. 
Then $\mathop{\mathrm{Sing}}(v_x) = \{ x \}$ and $\overline{O(y)} = \mathcal{M}$ for any non-singular point $y$.  
By Gutierrez's smoothing theorem~\cite{gutierrez1986smoothing}, we may assume that the flow $v_x$ is $C^\infty$. 
Let $X$ be the $C^\infty$ vector field generating $v_x$. 
Since every closed subset of any paracompact $C^\infty$ manifold is a zero set of some $C^\infty$ function on it by using $C^\infty$ bump functions and partitions of unity, take a $C^\infty$ bump function $\varphi \colon \mathbb{T}^2 \to [0,1]$ with $\varphi^{-1}(0) = \mathcal{M}$. 
Then the resulting flow $v$ generated by the vector field $\varphi X$ is a $C^\infty$ flow  
such that the $\omega$-limit set of a non-singular point is $\mathcal{M} = \Sv$. 
This means that $v$ is
a toral flow with an $\omega$-limit set consisting of singular points which is a transversely Cantor set.
\end{proof}


\begin{lemma}
There is a toral $C^\infty$ flow $w$ with an $\omega$-limit set of a point which is a quasi-semi-attracting limit quasi-circuit that is not locally connected and consists of two non-recurrent orbits and a subset of $\mathop{\mathrm{Sing}}(w)$ which is a transversely Cantor set as a set. 
Moreover, the flow can contain a non-recurrent point in the $\omega$-limit set whose orbit closure is not arcwise-connected.  
\end{lemma}

\begin{proof}
Let $X_1 := \varphi X$ be the vector field in the previous example. 
Then the complement $\mathbb{T}^2 - \mathcal{M}$ is a trivial flow box. 
Identify the flow box with a square $B := \R \times (-1/2,\pi+1/2)$ on a chart by the embedding $f \colon B \to \mathbb{T}^2 - \mathcal{M}$. 
Fix the Euclidean norm on $B \subset \R^2$ and the norm on $\T^2$ induced by the canonical quotient map $\R^2 \to \R^2/\Z^2$. 
Replacing the norm on $B$ with the norm on $B$ obtained by multiplying the norm by a positive scalar if necessary, we may assume that the norms of the vector at each point of any vector fields on $B$ are no more than the norms of the vector at each point of the push-forwards of the vector fields by $f$ on the square $\mathbb{T}^2 - \mathcal{M}$. 
Define $C^\infty$ bump functions $h, \psi \colon \R \to [0,1]$ with $\psi^{-1}(0) = \R - (-1/2,\pi+1/2)$ and $\psi^{-1}(1) = [0,1]$ such that $h$ is an even function which is strictly decreasing on $(0, \infty)$ and $\lim_{x \to \pm \infty} h(x) = 0$. 

Define the following continuous vector field $Y_0$ on $(-\pi/2, \pi/2) \times [-1/2,\pi+1/2]$ as follows: 
\[
Y_0(\theta,y) := 
\begin{cases}
(-\cos \theta, 0) & \text{for } y \in (\pi, \pi+1/2] \\ 
(\cos \theta, 0) & \text{for } y \in [-1/2, 0) \\ 
(\cos \theta \cos y, \sin \theta \sin y) & \text{for } (\theta,y) \in (-\pi/2, \pi/2) \times [0,\pi] \\ 
\end{cases}
\] 
Notice that the restriction $Y_0\vert_{(-\pi/2, \pi/2) \times [0,\pi]}$ is a Taylor-Green vortex. 
Define the vector field $Y_1$ on the square $B$ by $Y_1(x,y) := h(x)\psi(y)Y_0(\tan^{-1}(x),y)$. 
Then $Y_1(x,0) = (\psi(0)\cos (\tan^{-1}(x)), 0) = \psi(0)\cos (\tan^{-1}(x)) (1, 0)$ is non-singular on the line $\R \times\{0\}$ and $Y_1(x,\pi) = (-\psi(\pi)\cos (\tan^{-1}(x)), 0) = \psi(\pi)\cos (\tan^{-1}(x)) (-1, 0)$ is non-singular on the line $\R \times\{\pi\}$. 
Moreover, we have that $Y_1(x,-1/2) = 0$ and $Y_1(x,\pi + 1/2) = 0$. 
Since the restrictions $Y_1\vert_{\R \times ( [-1/2, 0) \sqcup (\pi, \pi+1/2] )}$ and $Y_1\vert_{\R \times [0,\pi]}$ are $C^\infty$, the vector field $Y_1$ is locally Lipschitz continuous. 
Then the push-forward $Y_2 := f_* Y_1$ on $\mathbb{T}^2 - \mathcal{M}$ by $f$ of the vector field $Y_1$ generates an $\R$-action $v_{Y_2}$ on $B$ as Figure~\ref{bufb03}. 
Since $\vert Y_1(x,y)\vert  \geq \vert Y_2(f(x,y))\vert $ for any point $(x,y) \in B$, by $\lim_{x \to \pm \infty}\max_{y \in [-1/2,\pi+1/2]} \vert Y_1(x,y)\vert  = 0$, the induced vector field $Y_2$ can be extended to a continuous vector field $Y_3$ on $\mathbb{T}^2$ by $Y_3\vert_{\mathcal{M}} = 0$. 
Since the closed subset $\mathbb{T}^2 - B$ consists of singular points of the $\R$-action $v_{Y_3} \colon \R \times \mathbb{T}^2 \to \mathbb{T}^2$ generated by $Y_3$, Lemma~\ref{lem:continuity_extension} implies that the $\R$-action $v_{Y_3} \colon \R \times \mathbb{T}^2 \to \mathbb{T}^2$ is a flow. 
By Gutierrez's smoothing theorem~\cite{gutierrez1986smoothing}, we may assume that the flow $v_{Y_3}$ is a desired $C^\infty$ flow. 
\begin{figure}
\begin{center}
\includegraphics[scale=0.3]{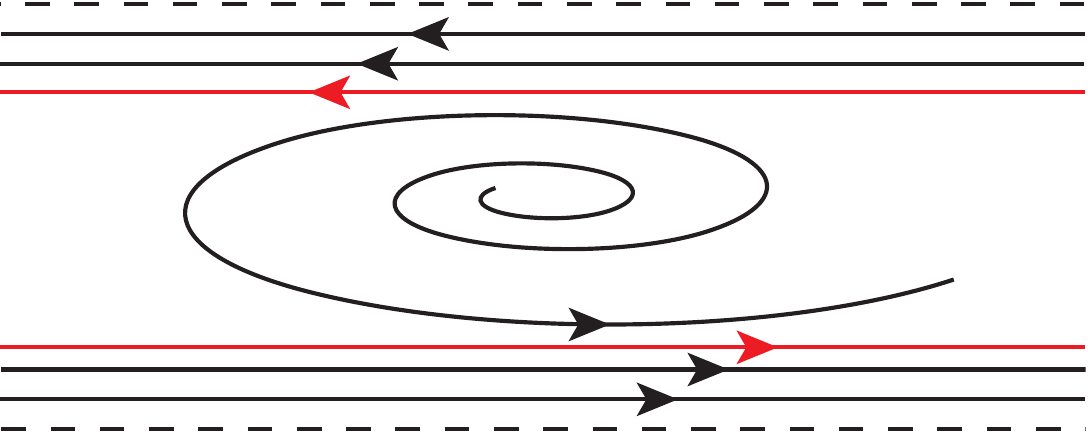}
\end{center}
\caption{A square with a flow.}
\label{bufb03}
\end{figure}
\end{proof}

Though the orbit closure of a non-recurrent point in the $\omega$-limit set in the above proof is not arcwise-connected and is the disjoint union of a non-recurrent orbit and a transversely Cantor set in the singular point set, notice that the orbit closure of a non-recurrent point in the $\omega$-limit set of a point for a flow with totally disconnected singular points on a compact surface is a closed arc because of Theorem~\ref{main:a}.




\bibliographystyle{abbrv}
\bibliography{../yt20210901}
\end{document}